\documentclass[12pt]{article}
\usepackage[left=3cm,right=3cm,top=3cm,bottom=3cm]{geometry}
\usepackage[french]{babel}

\usepackage{amsmath, amsfonts, amssymb,latexsym, bm}
\usepackage{wasysym,graphicx,stmaryrd,tikz,bclogo}
\usepackage{pgfplots}

\addto\captionsfrench{}

\usepackage{amsmath, amsfonts, amssymb,latexsym}

\input epsf.sty

\newtheorem{Theorem}{Theorem}[section]
\newtheorem{Lemma}{Lemma}[section]
\newtheorem{Proposition}[Lemma]{Proposition}

\newtheorem{Definition}[Lemma]{Definition}

\newcommand{\BEQ}{\begin{equation}}     
\newcommand{\BEA}{\begin{eqnarray}}
\newcommand{\BD}{\begin{displaymath}}
\newcommand{\EEQ}{\end{equation}}       
\newcommand{\EEA}{\end{eqnarray}}
\newcommand{\ED}{\end{displaymath}}

\newcommand{\del}{\delta}
\newcommand{\Del}{\Delta}
\newcommand{\eps}{\varepsilon}          

\newcommand{\Tr}{{\mathrm{Tr}}}

\newcommand{\Z}{\mathbb{Z}}
\newcommand{\N}{\mathbb{N}}

\newcommand{\Id}{{\mathrm{Id}}}

\def\T{{\mathbb{T}}}

\newcommand{\eop}{\hfill $\Box$}        
\newcommand{\Medskip}{\medskip\noindent}
\newcommand{\Bigskip}{\bigskip\noindent}

\newcommand{\half}{{1\over 2}}          

\catcode`\@=11
\def\numberbysection{\@addtoreset{equation}{section}
        \def\theequation{\thesection.\arabic{equation}}}
\numberbysection

\begin{document}
\renewcommand{\contentsname}{Table of contents}
\renewcommand{\bibname}{References}

\title{\bf Optimal multi-time-scale estimates for diluted autocatalytic
chemical networks. (1) Introduction and $\sigma^*$-dominant case.}
\vskip -2cm
\date{}
\vskip -2cm
\maketitle

\vspace{2mm}
\begin{center}
{\bf  J\'er\'emie Unterberger\footnote{IECL (Institut Elie Cartan de Lorraine), Université de Lorraine, Nancy, France}}
\end{center}

\vspace{2mm}
\begin{quote}

\renewcommand{\baselinestretch}{1.0}
\footnotesize
{Autocatalytic chemical networks are  dynamical systems whose linearization around zero
 has a positive Lyapunov exponent; this exponent gives
the growth rate of the system in the diluted regime, i.e. for
near-zero concentrations. The generator of the dynamics in the kinetic limit is then a Perron-Frobenius matrix,
suggesting the use of Markov chain techniques to get long-time asymptotics. 

This series of works introduces a new, general  procedure providing precise quantitative information about such asymptotics, based on estimates for the Lyapunov eigenvalue and
eigenvector.  The
algorithm, inspired from Wilson's renormalization group method in quantum field theory, is based
on a downward recursion on kinetic scales, starting from the fastest, and terminating with the
slowest rates.   Estimates  take on the form of simple rational functions of kinetic rates. They are accurate under a separation of scales hypothesis, loosely stating that kinetic rates span many orders of magnitude.   

We provide here a brief general motivation and  introduction to the method, present some simple examples, and 
derive a number of preliminary results, in particular the estimation of Lyapunov data for a subclass of 
so-called {\em $\sigma^*$-dominant} graphs.
    }

\end{quote}
\vspace{4mm}
\noindent


\smallskip
\noindent

\tableofcontents


\section{Growth rate of chemical networks} 


Large, complex chemical reaction networks are pervasive in the recent chemical literature; they are found in several
fields of interest, such as atmospheric chemistry  \cite{Jen1}, \cite{Jen2}, \cite{Top}, combustion
phenomena \cite{Boi}, both of interest in particular for the study of pollution; food processing,
thinking in particular about a long-studied set of browning  reactions very relevant for
flavor, called 
Maillard reaction \cite{Hod}, but also present in prebiotic reactions identified by astrochemists in pre-accretional organic analog formed from dense molecular ice analogs \cite{Dan}; and  bioinspired carbon conversion cycles, which are laboratory alternatives to the classical metabolic cycles that may run
in the absence of enzymes, and are of theoretical interest in the quest for the origin of life \cite{Much}. Mathematically, they take on the form of large-dimensional coupled ODEs (ordinary
differential equations), $dX/dt = F(X)$, where $X=X(t)$ is a vector, called {\em species distribution}, representing the concentration
of each chemical species, and $F$ is a  vector field parametrized by coefficients $(k_{\rho})_{\rho}$ called {\em kinetic rates}, with $\rho$ indexing reactions. We place ourselves here directly in the kinetic limit, where the evolution is
deterministic instead of stochastic, see \cite{FalEsp}, and assume so-called "mass action" rates, implying that $F$ is polynomial.

\Medskip  The general aim is to  obtain quantitatively
correct dynamical predictions for the species distribution. These face two main theoretical difficulties; the first one is the lack
of kinetic data, see e.g. \cite{vanBoe} for the Maillard reaction. The second one is the complexity of the networks; the number of compounds and reactions, in 
many experimental sets in the field,
is virtually infinite.

\Medskip In view of these difficulties, the aim of this series of papers is to provide a {\em   generic analysis of the long-term dynamics of complex chemical networks}, based upon the sole hypothesis of {\em scale
separation}. Though our methods are only semi-quantitative, the relative error in the predictions is in principle
of the same order as that in the kinetic rates. Furthermore, our predictions take on the form of simple
mathematical expressions in terms of the kinetic rates viewed as parameters, a key point for the inference of  kinetic rates from experimental data, which is a long-term project based on the present work, connecting theory with experiments.


\subsection{General context} \label{subsection:general-context}


We refer the reader e.g. to \cite{Fei} or \cite{Len} for a general discussion  of chemical reaction
networks.
Linearizing the ODEs around 0 yields a system of the general form $dX/dt = AX$, where $A$
is a Perron-Frobenius matrix associated to a set of reactions forming a hypergraph connecting reactants to products. A reaction
is of the general form $\rho: X_{\sigma} \overset{k_{\rho}}{\to} X_{\sigma'_1}+\cdots+X_{\sigma'_n}   $, where $\sigma\not= \sigma'_1,\ldots,\sigma'_n$ belong to a fixed  index set $\Sigma$,    the species $X_{\sigma}$ is 
called the  {\em reactant}, and species $X_{\sigma'_1},\ldots,X_{\sigma'_n}$ are the {\em products}. This reaction
contributes additively to the coefficients of the generator, $+k_{\rho}$ to off-diagonal coefficients  $A_{\sigma'_i,\sigma_i}$, $i=1,\ldots,n$,  $-k_{\rho}$ to the diagonal coefficient
$A_{\sigma,\sigma}$, so that the coefficients of $A$ are obtained by summing all these
 contributions for all reactions.  Note that reactions
with several reactants do not participate in $A$ since they contribute terms of order 
$\ge 2$ to the vector field $F(X)$.  Linearizing turns the hypergraph into a graph, whose
 edges (called: {\em split reactions}) are obtained from the list $\sigma\to\sigma'_i, i=1,\ldots,n$ with $\rho$ ranging in
 the set of reactions.  When $n=0$, $\rho$ is called a {\em degradation reaction}. When $n
 \ge 1$, 
$\rho$ is called a {\em 1-$n$ reaction}; in practice, one can concentrate on $1-1$ and $1-2$ reactions only, since reactions with $n\ge 3$ are never considered in chemistry. If only 1-1 reactions are present, $A$ is the adjoint generator
of a continuous-time Markov chain, see e.g. \cite{MeyTwe} for a general introduction
to the subject, and $k_{\rho}$ are transition rates.    As well-known  from standard Markov
chain theory, such matrices have a real top eigenvalue $\lambda^*$; assuming e.g. that $A$ is 
irreducible, the corresponding eigenspace is one-dimensional and generated by a vector $v^*$ with
positive components.  The fixed point 0 is unstable if and only if $\lambda^*$ is positive; such
chemical networks are called {\em autocatalytic}, and believed to play a role in the possibility for natural evolution starting from elementary biochemical bricks.  
 The foremost example of autocatalytic network
in chemistry is the formose reaction, a reaction network involving oses (sugars) in aqueous solution discovered by Butlerov
\cite{But} in 1861.

\Medskip Until very recently, there was no general
theory of autocatalytic networks. During the last
decade or so, however, interest shifted towards a systems 
chemistry point of view.
In 2020, Blokhuis et al.  \cite{BloLacNgh} gave a {\em stoechiometric} definition of autocatalysis for chemical reaction networks, stating the existence of a combination of reactions such that the balance for all autocatalytic species is strictly positive, and investigated minimal autocatalytic networks, called {\em autocatalytic cores}. In the same article,
 a complete list of minimal autocatalytic networks (so-called autocatalytic
  cores, or motifs) was established, and algorithms for
 automatic detection of such motifs in available reaction databases have been
 developed    \cite{Arya}, \cite{Hor}.
 
 \Medskip  By contrast,
  autocatalysis in our sense -- namely, exponential
 amplification of all species internal to a reaction network, starting from a {\em diluted regime}, i.e. low
 concentrations -- is a {\em dynamical} property. Networks with only 1-1 reactions have 
 $\lambda^*=0$ by Markov theory; not surprisingly, degradation, resp. $1-n$ reactions with $n\ge 2$,
 decrease, resp. increase $\lambda^*$. 
In a recent paper \cite{NgheUnt1}, we proved that, in absence of any degradation reaction, 
autocatalysis was a very general property for reaction networks,   and 
very closely related to the property of  stoechiometric autocatalysis.  
 Related works \cite{DesLacUnt1}, \cite{DesLacUnt2} study on chemically motivated examples the {\em autocatalystic threshold}, namely, the maximum value of the rate of a given
degradation reaction for which autocatalysis holds.


\subsection{Aim of the work} \label{subsection:aim}


\Bigskip {\em Separation of scales.} Our general aim in this series of articles  is to give accurate estimates of 
the long-term behavior of the linearized system $dX/dt = AX$  in great generality. 
The only hypothesis -- which may look innocent, but is not necessarily
well observed in practice -- is that of {\em separation of scales}.
Namely, we assume that {\em reaction rates do not 'accumulate' around any
given time-scale.} This can be made quantitative
by defining  the {\em scale} of a (reaction with)
kinetic rate $k$ to be the integer part 
 of the base $b$ logarithm
of $k$, 
\BEQ n(k) = \lfloor \log_{b} (k) \rfloor \label{eq:scale}
\EEQ
In practice, choosing $b\sim 10$ is very reasonable for experiments, given measurement errors when
attempting to determine kinetic rates from experiments. 
The choice of sign in (\ref{eq:scale}) implies that 
the larger $n(k)$, the faster the reaction. {\em Separation of 
scales thus means that for each integer $n$, there are 'not too many'
reactions with kinetic rate $k$ such that $n(k)=n$.}  A somewhat
more formal
statement would be that our estimates deteriorate, slowly but in general
irremediably, with the sizes of the connected components of the subgraphs obtained by
considering fixed
 scale  edges; in fact, even this condition is too drastic, as emphasized in our upcoming
 article \cite{Unt2}. 
Then, our estimates are {\em generic}, which means they are mathematically valid away from a {\em resonance set} obtained when some integer linear combinations of reaction scales vanish. 

\Bigskip Hypothesizing separation of scales, and leaving out the
delicate discussion of resonances, our results may be put bluntly
in this way.  {\em Assume the reaction network   is known, and
so are all reaction scales $n(k)$. Then we describe an explicit
algorithm   giving estimates of the Lyapunov exponent and eigenvector components in the simple form  
$\frac{P(k)}{Q(k)}$, where $k=(k_{\sigma})_{\sigma}$ is the vector of kinetic rates, and 
$P, Q$ are monomials (products of integer powers of the rates); equivalently, their logarithms are integer linear
combinations of the scales $n(k)$.  } In turn, this approach gives accurate estimates of 
the time behavior of the linearized system depending on explicitly described time regimes.

\vskip 2cm In the rest of the section, we loosely describe some
of the main mathematical tools, and say some words on the multi-scale
method.


\subsection{General setting and mathematical tools}  \label{subsection:maths}


We  assume for realism that {\em no reaction has more than two products}. This simplifies some
of the constants in the computations, but our results remain valid as soon as the number of products of
any reaction remains bounded by a constant (Theorem \ref{thm:homogeneous} below holds even for an infinite
network, so this very lenient assumption is not void). After linearization, we are thus left with
only three types of reactions:  degradation reactions $X_{\sigma} \overset{\beta_{\sigma}}{\to} 
\emptyset$;   1-1 reactions,
$X_{\sigma}\overset{k_{\sigma\to\sigma'}}{\to} X_{\sigma'}$; and  1-2  reactions, $X_{\sigma} \to 
X_{\sigma'} + X_{\sigma''}$.  

\Medskip  In absence of degradation reactions and one-to-many reactions, the linearized generator $A$ is the adjoint of the generator of a Markov chain, hence
$\lambda^*$ vanishes (and the stationary probability measure is 
a normalized null eigenvector). In general, however, {\em deficiency
rates}, defined as the species-dependent sum of the coefficients of $A$
over columns, and measuring a defect of probability preservation, are
non-zero. Degradation reactions decrease 
$\lambda^*$, while one-to-many reactions increase it, as can be
expected from purely stoechiometric considerations. If $A$ is irreducible, then there is 
a unique Lyapunov eigenvector $v^*$ up to normalization, which can be chosen with positive
coefficients.

\Medskip The {\em Resolvent Formula} is a well-known mathematical trick
in the study of Markov processes. Let $\alpha>0$ be a parameter (called: {\em threshold parameter}).  Formally, (assuming $\lambda^*>0$),
$A-\alpha$ has negative spectrum (more precisely, the real part
of eigenvalues is strictly negative) if $\alpha>\lambda^*$. Thus,
the inverse operator $(A-\alpha)^{-1}$ (known as resolvent in 
spectral theory) may be computed as the time-integral $\int_0^{+\infty} e^{t(A-\alpha)} \, dt$. Taking advantage of  the underlying Markov-like structure, matrix elements of $(A-\alpha)^{-1}$ may be computed as sum of (positive) weights of time trajectories of    
the associated  discrete-time Markov chain; these weights are products of dimensionless rate ratios $\frac{k_{\sigma\to\sigma'}}{\sum_{\sigma'} k_{\sigma\to\sigma'}}$, called {\em transition weights} (used 
in the Gillespie algorithm for the simulation of chemical networks).
 Then such sums diverge precisely
when $\alpha\to \alpha_{thr} \equiv \lambda^*$ (whence the name: threshold), yielding a  
{\bf self-consistent equation} for $\lambda^*$. The latter may be solved 
to a good approximation in whole generality. 
Very naively, we replace in general any sum of positive weights $W_1+W_2+\ldots$
by its maximum $\max(W_1,W_2,\ldots)$.  This is where our separation-of-scales hypothesis for the rates is important : the replacement of a sum by the max
is good if few 'dominant' $W_i$ have same order of magnitude, which is 
roughly tantamount to our hypothesis.  Then our solution is (save resonances) optimally good, in the sense that the relative error is
comparable to that in the kinetic parameters.


\subsection{Multi-scale method} \label{subsection:multi-scale-method}


The main line of thought in this work is directly borrowed from
 rigorous  renormalization group methods in 
quantum field theory (QFT), see e.g. \cite{FVT}, chap. 1 for a quick, 
self-contained introduction. Separation of scales in QFT is ensured by the smallness of a  coupling constant. Perturbative expansions for
the free energy or Green functions yield sums over Feynman graphs. 
Splitting the latter into different energy scales, it can be seen
that dominant contributions come in general from single-scale
graphs. Subleading corrections due to the insertion of higher-scale
subgraphs may be evaluated by a recursive, scale-dependent resummation procedure. A renormalization of the parameters of the theory keeps track scale after scale of 
the main, potentially divergent, corrections, and ensures one always
remains in a perturbative regime around a free theory. The renormalization procedure is best described pictorially as 
collapsing high-scale so-called {\em divergent subgraphs}.

\Medskip The above description translates almost verbatim to
the present context: Feynman graphs are just graphs; energy
scales are now rate scales; the resolvent 
plays the role of Green functions, and is evaluated (see \S 
\ref{subsection:maths}) by resumming time-integrals in the form of a sum of 
weights of paths along the graph of the reaction network. Renormalization 
 determines $\lambda^*$ to a large
extent by modifying the self-consistent equation discussed in \S \ref{subsection:maths}  above. 
Our representation convention for multi-scale graphs is to let scale indices decrease from top
to bottom (see \S \ref{subsection:example1-multi-scale-splitting}). Thus 
the film of chemical evolution shows gradually as time
increases (in a logarithmic scale) by scrolling down multi-scale
graphs.  Our general algorithm, leading to {\em hierarchical 
graphs}, a DAG (directed acyclic graph) structure much simplified compared to the original network,  will be described in our forthcoming
publication \cite{Unt2}. Because we do not tackle the general case, we shall not use
the full strength of the multi-scale structure  in
this article. Yet elementary multi-scale computations become central starting from Section 
\ref{section:multi-scale-splitting}, and DAGs appear in Section \ref{section:sigma*}.


\subsection{Important notations} \label{subsection:notations}


In the multi-scale setting, the notations below play an important
role.  Recall that scales are logarithms:
\BEQ n(k)= \lfloor \log_{b} k\rfloor \EEQ
is the scale of the (one-to-one or one-to-many) reaction  rate $k$. The base constant $b>1$ is fixed;
for applications, one would typically choose $b=10$.   Generally speaking, we order the
different transition scales $n_1,\ldots,n_{i_{max}}$ of the theory 'from top to bottom',
\BEQ n_1\succ n_2\succ\ldots\succ n_{i_{max}} \EEQ
and draw graphs as {\em multi-scale diagrams} (see Section \ref{section:multi-scale-splitting} and 
\S \ref{subsection:example1-multi-scale-splitting} for an illustration), where edges of scale $n_i$, $i<j$,
are placed above edges of scale $n_j$. 
 We sometimes define the scale of 
a discrete-time transition rate $w$ (see \S \ref{subsection:resolvent-formula}), in practice, a ratio of rates in $[0,1]$,
 as
\BEQ n(w) = - \lfloor \log_{b} w\rfloor \EEQ
the minus sign ensuring that $n\ge 0$. 
 Approximations are constantly used and are based on the
use of symbols $\prec,\succ,\preceq,\succeq,\sim,\simeq$.

\Medskip  First, scales; we  
write $n_1\prec n_2,\ {\mathrm{resp.}}\  n_1 \succ n_2,\  {\mathrm{resp.}}\  n_1\preceq n_2,\ {\mathrm{resp.}} \  n_1\succeq n_2$
if $n_1<n_2$ in the usual sense resp. $n_1>n_2$, $n_1\le n_2$, $n_1\ge n_2$. Then, about rates;  we write
\BEQ k\sim k_1 \EEQ
and say that {\em $k$ is approximately equal to $k_1$} 
if $k = \sum_{i=1}^j k_i$ with $n(k_1)\ge \max(n_{k_2},\ldots,n_{k_j})$.
Letting $b\to\infty$
in order to ensure separation of scales in a mathematically rigorous way, we would really have $k\sim_{b\to\infty} k_1$ for an (arbitrarily
large) but fixed network, provided $n(k_1)> \max(n_{k_2},\ldots,n_{k_j})+1$. However, it is difficult to bound the error involved in the symbol
$\sim$ for $b$ fixed. (Theorem \ref{thm:homogeneous} gives a rare instance where errors can be controlled
through explicit inequalities valid for an arbitrary size of the network).  In practice, one expects that 
$k/k_1$ is not necessarily close to 1, but there exists a constant $C<b$ or so  for
which $C^{-1}k_1<k<Ck_1$, which means in particular
that $\#\{i \ |\ k_i=k_1\}$ -- representing e.g. the number of
reactions $X\to \cdots$ of maximum rate scale for $X$ fixed --   should remain small.

\Medskip Similarly, if $k=\sum_{i=1}^j k_i$ with $n(k_i)<n$, $i=1,\ldots,j$, and 
$n(k')=n$,  we write $k\prec k'$
(equivalently, $k'\succ k$). 
 Finally,
$ (k\lesssim k')  \Leftrightarrow (k\prec k' \ {\mathrm{or}} \ k\sim k').$
When $k\sim k_1$, we also write
$n(k)\sim n(k_1),$
and expects in practice that $n(k)=n(k_1)-1,\ n(k)$ or $n(k)+1$.


\subsection{Overview of the article}


\noindent This series of articles is divided into three main parts, with increasing level of
complexity; only the first two  are included in the present article.

\Medskip
The first part (Sections \ref{section:resolvent} to \ref{section:example2}) is mainly based on our main formula: the {\bf Resolvent Formula} (Proposition \ref{prop:main-formula}), an implicit formula
for the Lyapunov exponent (more or less standard in probability theory),
based on the introduction of the {\em threshold parameter} $\alpha$,  which is the starting point for all subsequent computations. General mathematical tools (mostly coming from Markov chain theory)  are presented in Section \ref{section:resolvent}, in connection with the Resolvent Formula. Section \ref{section:easy} presents two readily usable results
on $\lambda^*$ : an elementary upper bound, and an estimate (with some indications of proof) based on the resolvent formula, valid when all 1-1 reactions have the same scale.  This is as far as one can go in general without
using multi-scale methods. Multi-scale splitting is introduced in Section \ref{section:multi-scale-splitting}. Splitting the reaction graph into multi-scale diagrams makes it possible to 'solve' easily the resolvent formula in some simple low-dimensional cases, that is, to get the scale of the Lyapunov exponent; we illustrate this pedestrian
 approach on Example 1 (a Type I, two-species network) and Example 2 (an expansion of Example 1 with
 one extra species), and validate our results with analytic/numerical data. The details of the arguments in  Example 2 make it quite clear that a more powerful approach is
 required to deal with more complicated networks.  
 
\Medskip  The second part (Sections \ref{section:sigma*} and \ref{section:formose})  shows how to 'solve' the Resolvent Formula for a specific 
class of networks for which the  {\em $\sigma^*$-dominant Hypothesis} holds. These are networks whose
autocatalytic cycles all pass through a fixed species $\sigma^*$. No recursion is needed in that case;
it suffices to estimate the weight of a finite generating set of loops through $\sigma^*$. Namely, 'bad' cycles (see below) are absent under the above hypothesis.  The scale
interval $\{n_{i+1},\ldots,n_i\}$ where $\alpha_{thr}$ is located is derived by cutting edges of scales $<n_{i+1}$ for $i=1,2,\ldots$ until the threshold is crossed.    At the same time, we
compute the Lyapunov eigenvector, see Lemma \ref{lem:Lya-weight} for
irreducible graphs.
We rederive in a few lines at the end of Section \ref{section:sigma*} the estimates obtained earlier for Examples 1 and 2, and
apply  in Section \ref{section:formose} the method of Section \ref{section:sigma*} to an idealized formose network, for which a more 
pedestrian approach would not be appropriate.

\Medskip Some of the  more technical results are postponed to Appendix (Section \ref{section:appendix}). A  notation index is given before the references.


\section{General notions and Resolvent formula}  \label{section:resolvent}


Because we are interested specifically in the linearized dynamics, 
the chemical network -- a hypergraph -- may be substituted by
a defective Markov graph $G=(V,E)$, where $V$ is the set 
of species, and $E$ a set of reactions $\sigma\to\sigma'$  with rates
$k_{\sigma\to\sigma'}$ encoding the linearized dynamics.   This is discussed in \S 
\ref{subsection:Markov}. The next subsection (\S \ref{subsection:resolvent-formula}) goes on to discuss the
Resolvent Formula, see Proposition \ref{prop:main-formula} (proved
in \cite{NgheUnt1}), which is our starting point for  estimating the
Lyapunov exponent.


\subsection{Markov theory}  \label{subsection:Markov}


For the whole subsection, and the beginning of the next one, we refer the non-expert reader to a 
good book on Markov chains, see e.g. \cite{MeyTwe}, or simply to the webpage \cite{Markov}  
 that gives a brief summary.

\Medskip {\em Hypergraph of reactions.} Let $V$ be the set of species, and  $R$ be the set of reactions of a chemical network with mass-action rates.
This set contains 1-1 reactions, one-to-many reactions {\em (in practice, we shall consider only
one-to-two reactions for the sake of realism, though extension to higher order reactions is 
mathematically straightforward)}, and
possibly, degradation reactions $\sigma\overset{\beta_{\sigma}}{\to}
\emptyset$, $\sigma\in V$ (considered as one-to-zero), which are
the only external flows (actually, outflows) of the system.
Denote by $R_p$, $p = 0,1,2$ the set of reactions with $p$ products.
These structures define a hypergraph Hyper$G=(V,R)$. 

\Medskip {\bf Adjoint (generalized Markov) generator.} Linearizing the dynamical system around 0,
we get $\frac{dX}{dt} = AX$, where $A$ is the {\em adjoint generator}. 
1-1 reactions $X_{\sigma}\overset{k}{\to}
X_{\sigma'}$ contribute $-k$, resp.
$k$ to the $(\sigma,\sigma)$, resp. $(\sigma',\sigma)$ entries of $A$. Summing up all 1-1 reactions gives
a matrix $A^{1-1}$ whose transpose $(A^{1-1})^t$ is a {\em Markov generator}, i.e. 
off-diagonal entries $(A^{1-1}_{\sigma',\sigma})_{\sigma\not=\sigma'}$ are positive, and $\sum_{\sigma'}A^{1-1}_{\sigma',\sigma}=0$ (probability conservation).  On the other hand, degradation reactions 
$X_{\sigma}\overset{\beta_{\sigma}}{\to} \emptyset$ only contribute 
a negative quantity $-\beta_{\sigma}$ to the $(\sigma,\sigma)$-entry; and one-to-many reactions $X_{\sigma}\overset{k}{\to} s_1 X_{\sigma'_1} +\ldots+s_n X_{\sigma'_n}$ ($s_1,\ldots,s_n\in\N^*$, 
$s_1+\ldots+s_n>1$) give both a negative
 $-k$ contribution to the $(\sigma,\sigma)$-entry, and  positive
 ones $s_i k$ to the $(\sigma'_i,\sigma)$-entries. In total, in all
 cases,
the contribution to $\sum_{\sigma'} A_{\sigma'\sigma}$ is
$(p-1)k$, where $p=s_1+\ldots+s_n$ $(p=0,1,2)$ is the {\em order} of the 
reaction. Thus {\em probability conservation does not hold for
$A^t$}. In fact, the {\em deficiency rate}
\BEQ \kappa_{\sigma} := \sum_{\rho\in R, \rho:X_{\sigma}\to\cdots}
(p(\rho)-1)k^{\rho}  \label{eq:deficiency-rate}
\EEQ
which is the sum of the $(\sigma,\sigma)$-contributions 
attached to all non 1-1 reactions with reactant $X_{\sigma}$,
is such that  $\sum_{\sigma'} A_{\sigma',\sigma} = \kappa_{\sigma}$. In absence of degradation reactions (all $\beta_{\sigma}=0$),
all deficiency rates $\kappa_{\sigma}$ are $\ge 0$; then (assuming $p\le 2$) $p(\rho)=2$ for
all reactions contributing to the sum in (\ref{eq:deficiency-rate}), whence {\em $\kappa_{\sigma}$ is
the sum of the rates of one-to-many reactions outgoing from $\sigma$}. Degradation
reactions contribute negative values.

\Medskip {\em Remark.} Negative $\kappa_{\sigma}$ are commonly interpreted as death rates in the probabilistic literature
\cite{Rev}; then $A^t$ is the generator of a Markov chain with death, and may give the average time evolution of a distribution of particles/species, with decreasing total mass. Positive $\kappa_{\sigma}$ 
could then be called birth rates, and the rate of variation of the total mass becomes unclear when these
are also included.  For lack of a standard terminology, we picked the word  "generalized".

\Medskip  {\bf   Lyapunov
data (exponent and eigenvector).} Matrices with non-negative off-diagonal coefficients (called: Perron-Frobenius, or Metzler matrices) such as $A$ are known by the Perron-Frobenius theorem to have a  {\em   Lyapunov
exponent}; namely, there is a unique eigenvalue $\lambda^*$ with maximal real part, it is  real, and admits a basis of 
eigenvectors with non-negative entries. If, furthermore, $A$ is {\em irreducible} (i.e. for every pair $(\sigma,\sigma')$,  there exists
a path $\sigma=\sigma_1\to\cdots\to\sigma_{\ell}=\sigma'$ with
$\sigma_i\not=\sigma_{i+1}$, $i=1,\ldots,\ell-1$, such that
$A_{\sigma_{i+1},\sigma_i}>0$), then the eigenspace of $\lambda^*$ is one-dimensional; 
the unique (up to normalization) eigenvector with non-negative entries, $v^*$,  is called {\bf Lyapunov eigenvector}.

\Medskip {\bf Associated Markov chain.} We let $\tilde{A}=(\tilde{A}_{\sigma,\sigma'})_{\sigma,\sigma'\in {\cal X}}$
be the matrix with  off-diagonal coefficients $\tilde{A}_{\sigma,\sigma'} =A_{\sigma,\sigma'}$ $(\sigma\not=\sigma')$ and negative diagonal coefficients
 
\BEQ \tilde{A}_{\sigma,\sigma}=-\sum_{\sigma'\not=\sigma} A_{\sigma',\sigma} = A_{\sigma,\sigma}-\kappa_{\sigma}.
\EEQ
 By construction, $\tilde{A}$ is an adjoint Markov generator, so
that its Lyapunov exponent (i.e. eigenvalue with largest real part)
is $0$. 

\Medskip {\em Transition rates.} Let  $k_{\sigma\to\sigma'} := A_{\sigma',\sigma}, \ \sigma\not=\sigma'$. 
We call $k_{\sigma\to\sigma'}$ the {\bf transition rate} from $\sigma$ to $\sigma'$. An interpretation
in terms of 'split' reactions is given in the next subsection.

\Medskip {\em Stationary measure.} If $\tilde{A}$ (equivalently, $A$) is 
{\em irreducible}, then there exists up to
normalization a unique stationary measure $\tilde{\mu}$. By definition,
$\tilde{A}\tilde{\mu}=0$; the constant vector ${\bf 1} = (1\ \cdots\ 1)^t$ satisfies the 
left eigenvalue equation ${\bf 1} \tilde{A}=0$ by probability conservation. If $A= \tilde{A}$ (only 1-1 reactions, no
degradation), then $\lambda^*=0$ and $v^*=\tilde{\mu}$ is an 
associated Lyapunov eigenvector. In some sense, all our results are
obtained by perturbing around this simple case (in particular, they
may be used to estimate the stationary measure of a continuous-time
Markov chain).

\Medskip {\em Markov graph.} Define $G=(V,E)$ to be the 
graph with edge set $E=\{(\sigma,\sigma')\ |\ \sigma\not=\sigma', \ 
\tilde{A}_{\sigma',\sigma}>0\}$. Then $\tilde{A}$ is irreducible if and only if  
$G$ is {\bf strongly connected} as a graph, i.e. if, for any pair of vertices 
$(v,v')\in V\times V$, there exists a path from $v$ to $v'$. In general, $G$ (if connected) may be partitioned into {\em strongly connected components} (maximal strongly connected subgraphs) which are 
connected by directed edges forming no cycle (the structure obtained by collapsing strongly connected
components is a directed acyclic graph).


\subsection{Resolvent Formula} \label{subsection:resolvent-formula}


The Resolvent Formula (proved in \cite{NgheUnt1}) is the key formula of this article; it allows a 
determination of the Lyapunov exponent as the solution of an implicit equation involving a discrete-time
generalized Markov chain. 

\Medskip We tacitly assume most of the time in this subsection (and in most
of the text) that $\tilde{A}$ is irreducible.

\Medskip {\bf Transition weights.} We let
\BEQ w_{\sigma\to\sigma'}= \frac{A_{\sigma',\sigma}}{|A_{\sigma,\sigma}|}, \qquad \tilde{w}_{\sigma\to\sigma'}= \frac{A_{\sigma',\sigma}}{|\tilde{A}_{\sigma,\sigma}|}, \qquad \sigma\not=\sigma'
\EEQ
These generalize the discrete-time transition weights used 
in the Gillespie simulation algorithm for the continuous-time
Markov chain; namely, starting from a vertex $\sigma$, wait for 
a random time following an exponential law  with parameter  $k_{\sigma}$ called 
{\bf outgoing rate}, 
\BEQ k_{\sigma} = \sum_{\sigma'\not=\sigma} k_{\sigma\to\sigma'},
\label{eq:k-sigma}
\EEQ
 and
then jump to $\sigma'$ with probability $w_{\sigma\to\sigma'}$.  Note that
\BEQ |\tilde{A}_{\sigma,\sigma}| = k_{\sigma}, \qquad 
\sum_{\sigma'\not=\sigma} \tilde{w}_{\sigma\to\sigma'} = 1
\EEQ 
by construction, so that $\tilde{w}_{\sigma\to\sigma'}$ are standard probabilistic weights, but  $\sum_{\sigma'\not=\sigma} w_{\sigma\to\sigma'} \not= 1$ in general.  Averaging the position $\sigma(t)$ at time $t$ over all random realizations, one gets
$\langle \del_{\sigma(t),\sigma'}\rangle = (e^{tA})_{\sigma',\sigma}$.
This gives the connection from Markov theory to the linearized 
dynamics.  The integral formula below for the resolvent may be seen
as a reformulation of the above.

\Medskip
More generally, if $\alpha\ge 0$ is a uniform degradation rate,
we let $A(\alpha):=A-\alpha\Id$, $\tilde{A}(\alpha):=\tilde{A}-\alpha\Id$, and denote by $w(\alpha)$, $\tilde{w}(\alpha)$ the associated 
transition weights,
\BEQ  w(\alpha)_{\sigma\to\sigma'}= \frac{A_{\sigma',\sigma}}{|A_{\sigma,\sigma}| + \alpha}, \qquad \tilde{w}(\alpha)_{\sigma\to\sigma'}= \frac{A_{\sigma',\sigma}}{|\tilde{A}_{\sigma,\sigma}| + \alpha}, \qquad \sigma\not=\sigma'
\label{eq:w(alpha)}
\EEQ
so that $w=w(0),\tilde{w}=\tilde{w}(0)$. 

\Medskip {\em Stationary measure, stationary weights.} A stationary measure
of the discrete-time Markov chain with transition weights
$\tilde{w}$ is a measure  $\tilde{\pi}$ such that 
\BEQ  \sum_{\sigma\not=\sigma'} \tilde{\pi}_{\sigma}
\tilde{w}_{\sigma\to\sigma'}  = \tilde{\pi}_{\sigma'}.
\label{eq:piwpi}
\EEQ
 Let  $\tilde{\mu}_{\sigma} := k_{\sigma}^{-1} \tilde{\pi}_{\sigma}$; 
plugging the relations $\tilde{w}_{\sigma\to\sigma'} = \frac{k_{\sigma\to\sigma'}}{|\tilde{A}_{\sigma,\sigma}|}$ and $|\tilde{A}_{\sigma',\sigma'}| =  k_{\sigma'}$  into (\ref{eq:piwpi}) 
yields $\sum_{\sigma\not=\sigma'} \tilde{\mu}_{\sigma} k_{\sigma\to\sigma'} = |\tilde{A}_{\sigma',\sigma'}|\, \tilde{\mu}_{\sigma'}$, from 
which $\tilde{A}\tilde{\mu}=0$. Thus $\tilde{\pi}$ is a stationary
measure of the discrete-time Markov chain if and only if $\tilde{\mu}$   is a stationary measure for $\tilde{A}$. 
Remember from the remark at the end of \S \ref{subsection:Markov} that
Lyapunov eigenvectors $v^*$  generalize stationary measures 
$\tilde{\mu}_{\sigma}$.  
We work most of the time with discrete-time transition weights, and
our results will provide estimates for the associated discrete-time
weights, called {\bf Lyapunov weights}, 
\BEQ \pi^*_{\sigma} := (|A_{\sigma,\sigma}|+\lambda^*) v^*_{\sigma}. \label{eq:pi*}
\EEQ   
Then $Av^* = \lambda^* v^*$ if and only if $\sum_{\sigma\not=\sigma'}
\pi^*_{\sigma} w^*_{\sigma\to\sigma'} = \pi^*_{\sigma'}$, where 
$w^*_{\sigma\to\sigma'} \equiv w(\lambda^*)_{\sigma\to\sigma'} = 
\frac{k_{\sigma\to\sigma'}}{|A_{\sigma,\sigma}|+\lambda^*}$.

\Medskip {\bf Path probability measure.} By construction,
$(\tilde{w}(0))_{\sigma\to\sigma'}$ are the transition probabilities of a true Markov chain. Therefore, they (formally) define
a probability measure 
\BEQ \tilde{w}(0)_{\gamma} := \prod_{i\ge 1}
\tilde{w}(0)_{\sigma_i\to\sigma_{i+1}}
\EEQ
  on the set of (infinite-length)  paths $\gamma:\sigma=\sigma_1\to \sigma_2
\to\cdots$ starting from some fixed state $\sigma$.

\Medskip {\bf Graph of split reactions.} Associate to each one-to-many reaction
$X\overset{k}{\to} s_1Y_1+\cdots+s_nY_n$ $(s_1+\ldots+s_n>1)$ the $n$ 1-1 'split' reactions $X\overset{s_i k}{\to}Y_i$ with
rates $s_i k$. Call $R_{split}$ the set of all 1-1
split reactions obtained in this way. Also (dealing with $R_0$), we add a 'cemetery' state $\emptyset$ to $V$ and represent the degradation of $X_{\sigma}$ as the
reaction $\sigma\to\emptyset$.  Then $R_0 \cup R_1 \cup R_{split}$
is the set of edges of a  multigraph, in the sense that two nodes
$X_{\sigma},X_{\sigma'}$ may be
connected by several edges $X_{\sigma}\to X_{\sigma'}$, each of one carrying a rate as upper index. 

\Medskip {\em  Markov graph.} We turn the above structure into a 
 Markov graph $G= (V,E), $
 with 1-1 rates denoted as $k_{\sigma\to\sigma'}$ by letting $E$ be the set of (oriented) pairs $(\sigma,\sigma')$
connected by one or several edges, and $k_{\sigma\to\sigma'}$ the
sum of all rates of 1-1 reactions or split reactions along the pair $(\sigma,\sigma')$.

\Medskip  We denote $A$ the associated  adjoint (generalized) Markov generator; 
its coefficients are
\BEQ A_{\sigma',\sigma} = k_{\sigma\to\sigma'}\qquad (\sigma\not=\sigma'), \qquad A_{\sigma,\sigma} = \kappa_{\sigma} 
- \sum_{\sigma'\not=\sigma} k_{\sigma\to\sigma'}.  \label{eq:A}
\EEQ
By construction, the linearization of the dynamical system around 0
is $\frac{dX}{dt}= AX$.

\Bigskip {\bf Resolvent.} Let $\alpha\ge 0$ be a positive parameter. The integral formula

 \BEQ R(\alpha)_{\sigma',\sigma}:=\int_0^{+\infty} dt\, 
(e^{tA(\alpha)})_{\sigma',\sigma}\in[0,+\infty]
\EEQ
 defines a matrix
with positive coefficients, which can be computed as a sum
over forward paths $\sigma= \sigma_1\to \sigma_2\to\cdots\to \sigma_{\ell}=\sigma'$ of arbitrary length $\ell\ge 0$,
\BEQ R(\alpha)_{\sigma',\sigma}=\sum_{\ell\ge 0} \sum_{\sigma_2,\ldots,\sigma_{\ell-1}\in V} \Big( \prod_{k=1}^{\ell-1}w(\alpha)_{\sigma_k\to \sigma_{k+1}}  \Big)
\, \times\, \frac{1}{|A_{\sigma',\sigma'}|+\alpha}.  \label{eq:traj}
\EEQ
 When finite, $R(\alpha)_{\sigma',\sigma}<\infty$
are the coefficients of the resolvent  $(\alpha-A)^{-1}=(-A(\alpha))^{-1}$; see
e.g.  \cite{Rev}, chap. III, or 
 \cite{Nor},  \S 4.2 for an introduction in connection to potential theory.
The connection to Lyapunov theory (see \cite{NgheUnt1}) is  that positivity of the Lyapunov exponent $\lambda^*$  of $A(\alpha)$
  is equivalent to having
\BEQ (R(\alpha))_{\sigma',\sigma}=+\infty  \label{eq:Ralpha=infty} \EEQ
for some (or all) $\sigma,\sigma'\in V$ and some $\alpha>0$.
Translating by $\alpha$, we thus have:

\BEQ \Big(\lambda^* \ge \alpha\Big) \Leftrightarrow 
\Big(\exists \sigma,\sigma'\in \Sigma, R(\alpha)_{\sigma',\sigma} = +\infty\Big)
 \Leftrightarrow 
\Big(\forall \sigma,\sigma'\in \Sigma, R(\alpha)_{\sigma',\sigma} = +\infty\Big)
\EEQ

\Medskip {\bf Excursions.} Let $\sigma,\sigma'\in V$. An {\em excursion} from $\sigma$ to $\sigma'$ is a path $\sigma=\sigma_1
\to \ldots \to \sigma_{\ell}=\sigma'$ (in particular, a loop if $\sigma'=\sigma$) which does not pass through
$\sigma$ at any intermediate step, i.e. such that $\sigma_2,\ldots,\sigma_{\ell-1}\not=\sigma$. Let $f(\alpha)_{\sigma\to\sigma'}\in 
[0,+\infty]$ be the total weight of
excursions from $\sigma$\ to $\sigma'$, namely,

\BEQ f(\alpha)_{\sigma\to \sigma'} =  
\sum_{\gamma:\sigma\to\sigma'} w(\alpha)_{\gamma} = \sum_{\ell\ge 1} \sum_{\sigma=\sigma_1\to \cdots
\to \sigma_{\ell}=\sigma', \ \sigma_2,\ldots,\sigma_{\ell-1}\not 
= \sigma} \prod_{k=1}^{\ell-1} w(\alpha)_{\sigma_k\to \sigma_{k+1}}
\EEQ

From (\ref{eq:w(alpha)}), it is apparent that $\alpha\mapsto f(\alpha)_{\sigma\to\sigma'}$ is a decreasing function. {\em From now on, we focus on the case when $\sigma'=\sigma$.}  It is actually
easy to prove that either $f(\alpha)_{\sigma\to\sigma}<+\infty$ for all $\alpha\ge 0$, or there exists $\alpha_0\ge 0$ such that $f(\alpha)_{\sigma\to\sigma}=+\infty$ if and only if $\alpha\le \alpha_0$; 
furthermore, the function is continuous and strictly decreasing on the interval where it is finite.
Then $R(\alpha)_{\sigma,\sigma}=\frac{1}{|(A(\alpha))_{\sigma,\sigma}|} \sum_{p\ge 0} (f(\alpha)_{\sigma\to \sigma})^p.$ 
Thus 
\BEQ \Big(R(\alpha)_{\sigma,\sigma}=+\infty\Big) 
\Leftrightarrow \Big( f(\alpha)_{\sigma,\sigma}\ge 1\Big)
\EEQ
In other words:

\begin{Proposition}[Resolvent formula] \label{prop:main-formula}
\BEQ \Big(\lambda^*\ge \alpha\Big) \Leftrightarrow 
\Big(\exists \sigma\in V,   f(\alpha)_{\sigma\to\sigma}\ge 1\Big)   \label{eq:main-formula}
\EEQ
or equivalently,
\BEQ \Big(\lambda^*\ge \alpha\Big) \Leftrightarrow 
\Big(\forall \sigma\in V,   f(\alpha)_{\sigma\to\sigma}\ge 1\Big)
\EEQ
\end{Proposition}

\Medskip {\bf Excursion probability measure.} Let  $\sigma\in V$. We can restrict
the above-defined path probability measure to excursions from
$\sigma$ to $\sigma$. Since finite-dimensional irreducible
Markov chains are recurrent, the total weight
$ \sum_{\gamma : \sigma\to\sigma} \tilde{w}(0)_{\gamma}$ of 
all excursions from $\sigma$ to $\sigma$ is equal to $1$. This
defines a probability measure $\langle \ \cdot \ \rangle_{\sigma}$ on the set of excursions from $\sigma$ to $\sigma$.


\section{Two easy results}   \label{section:easy}


 We  state in this section a formula and an estimate for $\lambda^*$ showcasing the interplay
between deficiency rates (equivalently, of  rates of one-to-many reactions, in absence of degradation)
and 1-1 reaction rates. 


\subsection{An upper bound}


Recall $\kappa_{\sigma} = \sum_{\rho \in R;
\rho:X_{\sigma}\to \cdots} 
(p(\rho)-1) k^{\rho}  = \sum_{\rho \in R_2; \rho:X_{\sigma\to\cdots}} k^{\rho}$ (assuming all one-to-many reactions are 1-2, see (\ref{eq:deficiency-rate})) is the deficiency rate of species $\sigma$. Linear
currents $J^{\rho}_{lin}$ are defined as
$J^{\rho}_{lin}(v) = k^{\rho} v_{\sigma(\rho)}$, where $\sigma(\rho)$ is the index of the
reactant of $\rho$. Let $({\mathbb S}_{\sigma,\rho})_{\sigma,\rho}$ be the stoechiometry matrix; by definition,
it is the rectangular matrix indexed by species and reactions with coefficients 
${\mathbb S}_{\sigma,\rho}= \begin{cases} -1, \qquad \sigma \ {\mathrm{reactant\ of\ }} \rho \\
+1, \qquad \sigma \ {\mathrm{product\ of\ }} \rho \\ 0 \qquad {\mathrm{else}} \end{cases}$. Products are counted with their multiplicities, e.g. ${\mathbb S}_{\sigma',\rho}=2$ if 
$\rho:X_{\sigma}\to X_{\sigma'}+X_{\sigma'}$, $\sigma'\not=\sigma$. 
   Note  that 
$\sum_{\sigma} {\mathbb S}_{\sigma,\rho} = p(\rho)-1$.  

\begin{Lemma}[deficiency rate identity]  \label{lem:1}
\begin{itemize}
\item[(i)] Let $v^*$ be a normalized  Lyapunov eigenvector, i.e.  $\sum_{\sigma} v^*_{\sigma}=1$. Then
\BEQ \lambda^* = \sum_{\sigma} \kappa_{\sigma}  v_{\sigma}^*     \label{id:1}
\EEQ

\item[(ii)] Assume  that the network $(\{X_{\sigma}\}, \{\rho\}_{\rho \in R_1})$
represents an irreducible, detailed-balanced Markov chain with
unique probability measure $(\pi^0_{\sigma})_{\sigma}$. 
Rescaling the rates of all other reactions, i.e. all
reactions of degree $\not=1$,    by a small coefficient $\eps$, i.e.
let $k^{\rho} \to \eps k^{\rho}$ for $\rho\in R_{\not=1}$. Then
\BEQ \lambda^* \sim_{\eps\to 0} \eps \sum_{\sigma} \kappa_{\sigma} \pi^0_{\sigma}    \label{id:1-eps}
\EEQ

\end{itemize}
\end{Lemma}

Thus, the growth rate (positive or negative) is a convex linear
combination of deficiency rates. In particular, the Lyapunov exponent
is {\em upper-bounded in absolute value by the largest deficiency rate}. For finite (not too large) networks, we
deduce the bound:
\BEQ |\lambda^*| = O(\max\{k_{\rho}, p(\rho)>1\}). \EEQ

\Medskip {\em Proof.}
\begin{itemize}
\item[(i)] If $v=v^*$ is the  Lyapunov eigenvector, then $Av=\lambda^* v$,
where  the adjoint generator $A$, see (\ref{eq:A}), is defined by 
$\sum_{\sigma'} A_{\sigma,\sigma'} v_{\sigma'}= \sum_{\rho} {\mathbb S}_{\sigma,\rho} J^{\rho}_{lin}(v)$. Taking 
the scalar product with some arbitrary vector $f = (f_{\sigma})_{\sigma\in V}$ in the space of species, we get the general
 formula:

\BEQ \lambda^* = \frac{\sum_{\sigma} f_{\sigma} \sum_{\rho}
{\mathbb S}_{\sigma,\rho} J^{\rho}_{lin}(v^*)}{\sum_{\sigma} f_{\sigma}
v^*_{\sigma}}   \label{eq:main}
\EEQ
 Choosing $f\equiv 1$ and permuting the two sums in the numerator leads to  (\ref{id:1}).

\item[(ii)] The unique
equilibrium probability measure $\pi^0=(\pi^0_{\sigma})_{\sigma}$ satisfies
$\pi^0_{\sigma'}k^0_{\sigma'\to\sigma}=\pi^0_{\sigma}k^0_{\sigma\to\sigma'}$. It is well-known and straightforward
to check that the matrix $\tilde{A}^0$ with coefficients
$\tilde{A}^0_{\sigma,\sigma'}= (\sqrt{\pi^0_{\sigma}})^{-1}
A^0_{\sigma,\sigma'} \sqrt{\pi^0_{\sigma'}}$ is symmetric.

\Medskip We now rescale all rate constants $k^{\rho}$ of one-to-many reactions
$\rho\in R_{>1}$ by a small parameter $\eps>0$, and
compute the asymptotics of $\lambda^*$ to leading order in $\eps$
as $\eps\to 0$.  Let $\tilde{A}_{\sigma,\sigma'} = (\sqrt{\pi^0_{\sigma}})^{-1}
A_{\sigma,\sigma'} \sqrt{\pi^0_{\sigma'}}$.   Then the matrix
$\tilde{A}=\tilde{A}^0 + \del A$ is the sum of a symmetric
matrix $\tilde{A}^0$ with top eigenvalue 0 and top eigenvector
$(\sqrt{\pi^0_{\sigma}})_{\sigma}$, and of a (in general non-symmetric) matrix $\del A$ with diagonal coefficients
 
\BEQ (\del A)_{\sigma,\sigma} = -\eps\sum_{\rho\in R_{>1}, 
\rho:X_{\sigma}\to \cdots} k^{\rho} 
\EEQ

and off-diagonal coefficients

\BEQ (\del A)_{\sigma,\sigma'} = \eps (\sqrt{\pi^0_{\sigma}})^{-1}
\Big(\sum_{\rho\in R>1, \rho:X_{\sigma'}\to X_{\sigma}+\cdots}
k^{\rho} \Big)  \sqrt{\pi^0_{\sigma'}}
 \qquad (\sigma\not=\sigma')  \label{eq:80}
 \EEQ
 
The standard Rayleigh perturbation theory (which applies
to non self-adjoint perturbations of self-adjoint operators)
yields to leading order
\BEA \lambda_{max}(\tilde{A}) &\sim& (\sqrt{\pi^0}, \del A \sqrt{\pi^0}) \nonumber\\
&=& -\eps\sum_{\sigma} \pi^0_{\sigma} \sum_{\rho\in R_{>1},
\rho:X_{\sigma}\to\cdots} k^{\rho} + \eps \sum_{\sigma\not=\sigma'} \pi^0_{\sigma'} \sum_{\rho\in R_{>2},
\rho:X_{\sigma'}\to X_{\sigma}+\cdots} k^{\rho}  \nonumber\\
&=& \eps \sum_{\sigma}  \Big( \sum_{\rho\in R_{>1},
\rho:X_{\sigma}\to\cdots} (p(\rho)-1)k_{\rho} \Big)\pi^0_{\sigma} = \eps \sum_{\sigma} \kappa_{\sigma} \pi^0_{\sigma}  
\EEA

\end{itemize}

\Medskip {\em Remark.} The validity of the perturbation argument implies directly that 
$v^*(\eps) \sim \pi^0$ in the limit $\eps\to 0$. Hence (\ref{id:1-eps}) may also be seen as a direct consequence of
(\ref{id:1}).


\subsection{Estimate in the homogeneous case} \label{subsection:homogeneous}


We state in this paragraph estimates which hold when the underlying
1-1 reaction network is homogeneous, i.e. when all non-zero
rates of 1-1 reactions $X_{\sigma}\to X_{\sigma'}$  are of the same scale. 

\Medskip Outgoing rates $k_{\sigma}$ are defined in (\ref{eq:k-sigma}). {\em We let $k_{\sigma}^{(0)} =
\sum_{\rho\in R_1; \rho: X_{\sigma\to\cdots}} k^{\rho}$
be the sum of the rates of 1-1 reactions outgoing from $\sigma$}. The Theorem below (proved in Appendix) is based on explicit, dimension- and network-independent,
inequalities, rather than on scale decomposition arguments as in \S \ref{subsection:notations};  as such, it is meaningful even for very large, or even infinite, networks. We write $a\lesssim b$  $(a,b>0)$ if there exists a constant $C>0$ such that $a\le Cb$, and 
$a\approx b$ (in words: $a$ and $b$ {\em are of the same order}) if $a\lesssim b$ and $b\lesssim a$. Statements are also true if one
replaces $a\lesssim b$ by $a\preceq b$, and $a\approx b$ by $a\sim b$, with the notations of \S
\ref{subsection:notations}.

\begin{Theorem}  \label{thm:homogeneous}
Consider a chemical network with zero degradation rates.
Assume all $k_{\sigma}^{(0)}$ are 
of the same order, i.e. there exists $k>0$ such that $\forall \sigma\in V, k^{(0)}_{\sigma}\approx k$. Then the following {\em three
scaling regimes} emerge:

\textbullet\ {\em (weakly autocatalytic regime)} if all $\kappa_{\sigma}\lesssim k$, i.e. if all autocatalytic rates are kinetically limiting,
then $\lambda^*\approx \max_{\sigma} \kappa_{\sigma}$ is of the same order as the leading autocatalytic
rate; 

\textbullet\ {\em (strongly autocatalytic regime)} if, on the
contrary, all $\kappa_{\sigma}\gtrsim k$, then 
 $\lambda^*\approx \min_{\sigma} \kappa_{\sigma}$ is of the same order as the smallest autocatalytic
rate; 

\textbullet\ {\em (intermediate regime)} if $\max_{\sigma}\kappa_{\sigma}\gtrsim k\gtrsim \min_{\sigma}\kappa_{\sigma}$, 
then  $\lambda^*\approx k$: 
transition rates of $1-1$ reactions are the limiting rates. 

\end{Theorem}

The applicability of the Theorem is not restricted to degradationless networks. 
Namely, adding degradation rates $\alpha=(\alpha_{\sigma})_{\sigma}$, 
i.e. letting $A_{\sigma,\sigma}\to A(\alpha)_{\sigma,\sigma}:=
A_{\sigma,\sigma}-\alpha_{\sigma}$, it was proved in \cite{NgheUnt1} that $ \lambda^*(A) - \max_{\sigma} (\alpha_{\sigma}) \le \lambda^*(A(\alpha)) \le 
 \lambda^*(A) - \min_{\sigma} (\alpha_{\sigma})
$.

\Medskip The proof relies on the resolvent formula (\ref{eq:main-formula}).  
Both generators $A,\tilde{A}$ have same off-diagonal coefficients,
but $|\tilde{A}_{\sigma,\sigma}|= k_{\sigma} =   |A_{\sigma,\sigma}|
+ \kappa_{\sigma}\ge |A_{\sigma,\sigma}|$. As follows from (\ref{eq:deficiency-rate}),  we actually have $|A_{\sigma,\sigma}| = k^{(0)}_{\sigma} + 
\kappa_{\sigma}$, $|\tilde{A}_{\sigma,\sigma}| = k^{(0)}_{\sigma} + 2\kappa_{\sigma}$; in particular, 
\BEQ |A_{\sigma,\sigma}|\le |\tilde{A}_{\sigma,\sigma}|\le 
2|A_{\sigma,\sigma}|  \label{eq:AAtilde2}
\EEQ
are of the same order. 
Let $\alpha\ge 0$. From $A(\alpha):=A-\alpha\Id$, $\tilde{A}(\alpha):=\tilde{A}-\alpha\Id$, one can define discrete-time generalized
Markov chains with transition kernels
\BEQ w(\alpha)_{\sigma\to\sigma'} = \frac{A_{\sigma',\sigma}}{|A_{\sigma,\sigma}|+\alpha}, \qquad 
\tilde{w}(\alpha)_{\sigma\to\sigma'} = \frac{A_{\sigma',\sigma}}{|\tilde{A}_{\sigma,\sigma}|+\alpha} \qquad \sigma'\not=\sigma
\EEQ
and then excursion weights $f(\alpha)_{\sigma\to\sigma'},\tilde{f}(\alpha)_{\sigma\to\sigma'}$. 
By definition, 
\BEQ w(\alpha)_{\sigma\to\sigma'}\ge \tilde{w}(\alpha)_{\sigma\to\sigma'}  
\label{eq:w>=wtilde}.
\EEQ 
By probability
conservation, $\tilde{f}(0)_{\sigma\to\sigma}=1$; then it is easy to prove by a first-step analysis
that the mapping $w\mapsto f(w)_{\sigma\to\sigma'}$, with $f = f(\alpha)$ or $\tilde{f}(\alpha)$, is
increasing, i.e. 
 $\Big(w'_{\sigma\to\sigma'}\ge w_{\sigma\to\sigma'},
\ \sigma\not=\sigma'\in V\Big)\Rightarrow \Big(f(w')_{\sigma\to\sigma'}
\ge f(w)_{\sigma\to\sigma'}, \ \sigma,\sigma'\in V\Big)$.  This is enough to get 
an explicit (network dimension independent) lower and an upper bound for $\lambda^*$, from which our Lemma follows.

\Bigskip We illustrate this result later on in Section \ref{section:example1} (Example 1). 
 For small networks, the conditions 
$k_{\sigma}^{(0)}\approx k$, $k$ fixed,  are equivalent to postulating simply that all rate constants of  1-1 reactions have
 scale $\sim n(k)$.


\section{Multi-scale splitting}  \label{section:multi-scale-splitting}


The Resolvent Formula
(\ref{eq:main-formula}) yields $\lambda^*$ as the solution
of the threshold equation, $\lambda^*=f_{\sigma\to\sigma}^{-1}(1)$, where $f_{\sigma\to\sigma}$ is the function $\alpha\mapsto
f(\alpha)_{\sigma\to\sigma}$; note that the solution is unique, since  $f_{\sigma\to\sigma}$ is 
continuous and strictly decreasing.  Unfortunately,
in most situations (except e.g. in the homogeneous case studied in
\S \ref{subsection:homogeneous}), it is difficult to extract information  directly from this implicit equation. 

\Medskip Multi-scale splitting  makes it possible to separate
the different contributions in $f(\alpha)$ according to their
scaling, and single out the most relevant ones, allowing to
determine the order of magnitude of $\lambda^*$.  We discuss multi-scale splitting at a formal level
in this section; see \S \ref{subsection:example1-multi-scale-splitting} for an example.

\Medskip  {\em Graph of split reactions.}  The starting point
is the graph of split reactions $G=(V,E)$
introduced in \S \ref{subsection:resolvent-formula}. By construction, 
$\sigma\to \sigma'$ is an edge in $E$ if (i) either there exists
a 1-1 reaction $X_{\sigma}\to X_{\sigma'}$;  or (ii) $X_{\sigma}
\to X_{\sigma'}$ is a split reaction coming from a one-to-many
reaction $X_{\sigma}\to X_{\sigma'}+\cdots$  There are no
multiple edges. The rate  $k_{\sigma\to\sigma'}$ attached to $(\sigma,\sigma')\in E$ is equal to the sum of the rates of reactions  (i) or
split reactions (ii) from $\sigma\to\sigma'$, counting
multiplicities.

\Medskip {\em Edge scales.} If $(\sigma,\sigma')\in E$,
we let 
\BEQ n_{\sigma\to\sigma'} := \lfloor \log_{b} (k_{\sigma\to
\sigma'}) \rfloor
\EEQ

\Medskip {\bf Dominant edges.} An edge
$\sigma\to\sigma'$ is dominant if 
\BEQ n_{\sigma\to\sigma'}  = \max_{\sigma''\in V}
n_{\sigma\to\sigma''}
\EEQ
 As a general rule, dominant edges appear as {\bf bold} lines on our representations of multi-scale graphs.

\Medskip {\bf Dominant paths/cycles.} A path $\gamma: X_{\sigma_1}\to\cdots\to X_{\sigma_{\ell}}$ is 
{\em dominant} if all the edges along $\gamma$ are dominant. A dominant path $\gamma$ is a dominant cycle
if $\sigma_{\ell}=\sigma_1$.

\Medskip {\bf Vertex scales.} A vertex $\sigma$ is an element
of $V$. The scale of $\sigma$ is 
\BEQ n_{\sigma}:= \max\{n_{\sigma\to\sigma'}, \sigma'\in V\}
\EEQ
i.e. the maximum scale of all edges $X_{\sigma}\to X_{\sigma'}$
 of $G$ outgoing from $X_{\sigma}$. 
 
\Medskip Note that an edge $\sigma\to\sigma'$ is dominant if 
$n_{\sigma\to\sigma'} = n_{\sigma}$.

\Medskip {\bf Total outgoing rate.}  Let
\BEQ k_{\sigma} := |\tilde{A}_{\sigma,\sigma}| = \sum_{\sigma'\not=\sigma} k_{\sigma\to\sigma'}.  \label{eq:k}
\EEQ
 By construction,  
 \BEQ \log_{b}  (k_{\sigma}) \sim n_{\sigma}. \EEQ

\Medskip {\em A short-hand notation for diagonal entries.} Let
\BEQ a_{\sigma}:= |A_{\sigma,\sigma}|. \label{eq:a} \EEQ
By construction, we always have $A_{\sigma,\sigma}< 0$, so that $a_{\sigma}=-A_{\sigma,\sigma}$.  
Similarly, $k_{\sigma} = -\tilde{A}_{\sigma,\sigma}$. If there is no degradation rate, then (from (\ref{eq:AAtilde2})), 
$|\tilde{A}_{\sigma,\sigma}|\simeq |A_{\sigma,\sigma}|$, hence $\log_{b} a_{\sigma} \sim n_{\sigma}$ too. 
These results still hold true if one
adds degradation rates $(\beta_{\sigma})_{\sigma}$ within our weak closedness hypothesis (see Assumption B in \S \ref{subsection:maths}), since $\beta_{\sigma}\ll |A_{\sigma,\sigma}|$.

\Medskip {\bf Deficiency weight.}  The deficiency
weight is the ratio
\BEQ \eps_{\sigma} : = \frac{\kappa_{\sigma}}{|A_{\sigma,\sigma}|} = \frac{\kappa_{\sigma}}{a_{\sigma}}.  \label{eq:deficiency-weight}
\EEQ

\noindent For example (see Example 1 in Section \ref{section:example1} below), a couple consisting of a  1-1 reaction $X_{\sigma}\overset{k_{off}}{\to} X_{\sigma'}$ and a 1-2 reaction
$X_{\sigma}\overset{\nu_+}{\to} 2X_{\sigma'}$ produces a deficiency weight $\eps_{\sigma} = \frac{\nu_+}{k_{off}+\nu_+}$. If 
$k_{off}\succ \nu_+$, then $\eps_{\sigma}\sim \frac{\nu_+}{k_{off}}$, 
otherwise $\eps_{\sigma}\sim 1$.  

\Medskip One always has $\kappa_{\sigma} \le |A_{\sigma,\sigma}|$, see (\ref{eq:AAtilde2}) and above, 
hence $\eps_{\sigma}\le 1$. We sometimes include degradation reactions $X_{\sigma}\overset{\beta_{\sigma}}{\to} \emptyset$ into the picture (a more systematic description of how this is done will be provided in our next publication).  If $\beta_{\sigma}=0$, then $\eps_{\sigma}\ge 0$;  as a consequence of the weak closedness assumption (B), see \S \ref{subsection:maths}, one always
have $\beta_{\sigma}\ll |A_{\sigma,\sigma}|$, hence, should it be the case that $\eps_{\sigma}<0$, one
would have: $|\eps_{\sigma}|\prec 1$.

\Medskip {\bf Scale ordering.} We let $(n_i)_{i=1,2,\ldots,n_{max}}$ be the
set of edge and vertex scales,  order them by decreasing order,
$n_1\succ n_2\succ\cdots$ and represent them from top to bottom, so 
that {\em scales decrease as one goes down a multi-scale graph}. By construction, there is 
a reaction $X_{\sigma}\to\cdots$ of scale $n_{\sigma}$ above any
reaction with reactant $\sigma$. In particular, there is at least 
one vertex of scale $n_1$.  

\Medskip {\bf Transition weights.}  Recall
$w_{\sigma\to\sigma'}= \frac{k_{\sigma\to\sigma'}}{|A_{\sigma,\sigma}|}$; hence,
\BEQ -\lfloor \log_{b}(w_{\sigma\to\sigma'}) \rfloor \sim n_{\sigma} -  n_{\sigma\to\sigma'}
\EEQ
Note that, if  $0\le \alpha\preceq k_{\sigma}$ (equivalently, $n_{\alpha}\preceq n_{\sigma}$), then 
$-\lfloor \log_{b}(w_{\sigma\to\sigma'}(\alpha)) \rfloor \sim n_{\sigma} -  n_{\sigma\to\sigma'}$ too.

\Medskip {\bf Path weights. }  The {\em weight} $w(\gamma)$ of a path $\gamma:X_{\sigma_1}\to\cdots\to X_{\sigma_{\ell}}$  is the
product of the transition probabilities $w$ along its edges,
$ w(\gamma)= \prod_{i=1}^{\ell-1} w_{\sigma_i\to\sigma_{i+1}}.$
From the above,
\BEQ -\lfloor \log_{b} w(\gamma) \rfloor \sim \sum_{i=1}^{\ell-1} 
\Big\{ n(\sigma_i) -  n(\sigma_i\to\sigma_{i+1})
 \Big\}  \label{eq:path-weight}
\EEQ
is a {\em sum of  differences of scales between nodes and edges}.  Paths from $X_{\sigma}$ to $X_{\sigma'}$ minimizing
the sum of differences of scales, i.e. the r.-h.s. of (\ref{eq:path-weight}), are called {\bf leading paths}. There
may be several ones, which all play an equivalent role.  {\em 
Dominant paths} are paths $\gamma$ such that  $\lfloor \log_{b} w(\gamma)\rfloor \sim 0$; a significant portion of time trajectories
flow through dominant paths.


\section{Example 1}  \label{section:example1}


  The simplest possible autocatalytic 
network (Type I in the classification of minimal autocatalytic cores in \cite{BloLacNgh}) is built upon a set of two species $X_0,X_1$,
\BEQ  X_0\overset{k_{on}}{\underset{k_{off}}{\rightleftarrows}}
X_1, \qquad X_1 \overset{\nu_+}{\to} 2X_0
\EEQ
It is a toy model e.g. of formose (see Example 3, Section \ref{section:formose}), summarizing
(looking from the distance) the successive addition of two abundant (chemostated) molecules with one carbon atom to $X_0$ (two carbons), followed by the splitting of $X_1$ (four carbons) into two molecules with
two carbons each. 

\Medskip The associated Markov graph is 
\begin{center}
\begin{tikzpicture}
\draw(0.2,0) node {$X_0$}; \draw(3,0) node {$X_1$};
\draw[->](0.5,0.2) arc(120:60:2 and 1);
\draw(1.5,0.8) node {$k_{on}$};
\draw(1.5,-0.8) node {$k_{off}+2\nu_+$};
\draw[->](2.5,-0.2) arc(-60:-120:2 and 1);
\end{tikzpicture} 
\end{center}

\noindent The Lyapunov exponent $\lambda^*$ is the maximum eigenvalue
of the matrix 
$A= \left(\begin{array}{cc}  -k_{on} & k_{off} + 2\nu_+ \\
k_{on} & -k_{off}-\nu_+ \end{array}\right).$
 Deficiency rates are $\kappa_0=0$ and $\kappa_1=\nu_+$. Assume that $k_{on}\approx k_{off}\approx k$ have same scale. Going back to Theorem \ref{thm:homogeneous},
we thus see that there are two cases: (i)  $\nu_+\preceq k$ (weakly autocatalytic regime); there, 
the autocatalytic rate $\nu_+$ is kinetically limiting, and $\lambda^*\simeq \nu_+$;  (ii) $\nu_+\succeq k$
(intermediate regime); there, 1-1 reactions are kinetically limiting, and $\lambda^*\simeq k$.  This is 
confirmed by (\ref{eq:lambda123}) below, from which we conclude that the above two estimates extend beyond the homogeneous regime; however, the formula in the kinetic regimes 2',3' where $k_{off} \gg k_{on}$, and $\nu_+$ is not
the highest rate, is different: it really follows a multi-scale logic, which is described at the end of 
the article.


\subsection{Analytic solution}  \label{subsection:example1-anal}


The matrix $A$ is of the form $A =\left[\begin{array}{cc} -a & b \\ c & -d \end{array}
\right]$ with $a,b,c,d>0$.  The determinant $\det(A)= - 
k_{on}\nu_+$ is $<0$. Thus 
the max eigenvalue (found by solving the second-order characteristic equation) 
 $\lambda^*=\half(-(a+d)+\sqrt{(a+d)^2-4\det(A)})
$
 is a strictly positive real number.
 Furthermore,   
 \BEQ \lambda^* \le  \frac{|\det(A)|}{|\Tr(A)|}=
\frac{|ad-bc|}{a+d} .
\EEQ
Namely, apply the inequality
$-x+\sqrt{x^2+y^2}=\frac{y^2}{x+\sqrt{x^2+y^2}}\le \frac{y^2}{2x}$ with $y^2=4|\det(A)| = 4k_{on}\nu_+$ and $x=|\Tr(A)| = k_{on} + k_{off} + \nu_+$. Note that 
the inequality is an equivalent, $-x+\sqrt{x^2+y^2}
\sim \frac{y^2}{2x}$ when $y/x\to 0$. Whatever the rates $k_{on},k_{off},\nu_+$, the above expressions
imply $y/x\le\sqrt{2}$. Thus we have
found the following upper bound for $\lambda^*$:
\BEQ \lambda^* \le \frac{k_{on}}{k_{on}+k_{off}+\nu_+}
\nu_+  \label{eq:TypeI-bd}
\EEQ 
which is actually an equivalent
in both limits $\nu_+\to 0,+\infty$. Furthermore, whatever the rates, we always have the estimate
\BEQ \lambda^* \approx \frac{k_{on}}{k_{on}+k_{off}+\nu_+}{\nu_+}
\label{eq:lambda-TypeI}
\EEQ

\Bigskip When $\nu_+\to +\infty$, we get the maximum
possible value: $\lambda^*\to k_{on}$; the 1-1 reaction
$X_0\to X_1$ is then the kinetically limiting step of 
the cycle. 
 
\Medskip When  $\nu_+\to 0$ instead (while $k_{on},k_{off}$ remain
constant), a physically interesting case in which the
doubling reaction $X_1\to 2X_0$ is a kinetically limiting step, 
we obtain the equivalent: 
\BEQ \lambda^* \sim \frac{k_{on}}{k_{on}+k_{off}} \nu_+. \EEQ  
  The idea is that (supposing
$k_{on},k_{off}$ are large) the system instantaneously equilibrates according to

\begin{center}
\begin{tikzpicture}
\draw(0.2,0) node {$X_0$}; \draw(3,0) node {$X_1$};
\draw[->](0.5,0.2) arc(120:60:2 and 1);
\draw(1.5,0.8) node {$k_{on}$};
\draw(1.5,-0.8) node {$k_{off}$};
\draw[->](2.5,-0.2) arc(-60:-120:2 and 1);
\end{tikzpicture} 
\end{center}

\noindent whence (dropping the doubling reaction 
$X_1 \overset{\nu_+}{\to} 2X_0$) it has probability $\mu^0_1:=\frac{X_1^{eq}}{X_0^{eq}+X_1^{eq}} = \frac{k_{on}}{k_{on}+k_{off}}$ of being in state $X_1$. Thus (in this limit) the system is
equivalent to the  1d system 
$X_1 \overset{\mu^0_1 \nu_+}{\to} 2X_1$, which has growth rate 
$\mu^0_1\nu_+ = \frac{k_{on}}{k_{on}+k_{off}} \nu_+$.  


\subsection{Multi-scale splitting}  \label{subsection:example1-multi-scale-splitting}


Eq. (\ref{eq:lambda-TypeI}) can also be obtained very simply from  the Resolvent Formula
 (Proposition \ref{prop:main-formula}). Instead of 
this single estimate, however, we get three different estimates depending on the ordering of the scales,
each of them  in the form $\frac{P(k)}{Q(k)} =$ (product of kinetic rates)/(product of kinetic rates), which are all equivalent
to  (\ref{eq:lambda-TypeI}).

\Medskip {\em Diagrammatic representation} (see Section \ref{section:multi-scale-splitting}).  We order the reaction scales
from highest $n_1$ to lowest $n_{i_{max}}$, i.e. $n_1>\ldots>n_{i_{max}}$.
In the diagrammatic representations, $n_1$ appears at the top,
and $n_{i_{max}}$ at the bottom, with scale indices decreasing as one
goes down.  A bullet (\textbullet) is attached to a vertex $\sigma$ (species) at
the highest scale at which there appears an outgoing reaction
$X_{\sigma}\to \cdots$, i.e., at scale $n_{\sigma}$.   Highest scale  edges $X_{\sigma}\to X_{\sigma'}$ outgoing from $\sigma$ 
are drawn as {\bf boldface lines}.  Let us draw the 6 possibilities:

\medskip
\begin{center}
\begin{tikzpicture}[scale = 0.81]  \label{fig:example1}
\draw(1.5,2) node {$1$};

\draw(-1,-1.5) rectangle(4,1.5);
\draw(2,1) node {\textbullet}; \draw(2.3,1) node {$X_1$};
\draw(3.5,1) node {$n_1$};
\draw(0.5,1) node {$X_0$};
\draw[dotted](2,1)--(2,-1);
\draw[dotted](0.8,1)--(0.8,-1);

\draw[dashed](0,0.5)--(3,0.5);

\draw[->, ultra thick](2,1)--(1,1); \draw[->, ultra thick](1,1)--(0.8,1);

\draw[->](2,0)--(0.8,0);
\draw(3.5,0) node {$n_2$};

\draw[dashed](0,-0.5)--(3,-0.5);
\draw[->, ultra thick](0.8,-1)--(2,-1);   \draw(0.8,-1) node {\textbullet};

\draw(3.5,-1) node {$n_3$};

\draw(1.5,-2) node {$n_1= \lfloor \log_{b}\nu_+\rfloor$};
\draw(1.5,-2.5) node {$n_2= \lfloor \log_{b}k_{off}\rfloor$};
\draw(1.5,-3) node {$n_3= \lfloor \log_{b}k_{on}\rfloor$};


\begin{scope}[shift={(6,0)}]
\draw(1.5,2) node {$1'$};

\draw(-1,-1.5) rectangle(4,1.5);
\draw(-0.75,-0.25) rectangle(3.75,1.25);
\draw(2,1) node {\textbullet}; \draw(2.3,1) node {$X_1$};
\draw(3.5,1) node {$n_1$};
\draw(0.5,1) node {$X_0$};
\draw[dotted](2,1)--(2,-1);
\draw[dotted](0.8,1)--(0.8,-1);

\draw[dashed](0,0.5)--(3,0.5);

\draw[->, ultra thick](2,1)--(1,1); \draw[->, ultra thick](1,1)--(0.8,1);

\draw[<-, ultra thick](2,0)--(0.8,0); \draw(0.8,0) node {\textbullet};
\draw(3.5,0) node {$n_2$};

\draw[dashed](0,-0.5)--(3,-0.5);
\draw[<-](0.8,-1)--(2,-1); 

\draw(3.5,-1) node {$n_3$};

\draw(1.5,-2) node {$n_1= \lfloor \log_{b}\nu_+\rfloor$};
\draw(1.5,-2.5) node {$n_2= \lfloor \log_{b}k_{on}\rfloor$};
\draw(1.5,-3) node {$n_3= \lfloor \log_{b}k_{off}\rfloor$};
\end{scope}


\begin{scope}[shift={(0,-6)}]
\draw(1.5,2) node {$2$};

\draw(-1,-1.5) rectangle(4,1.5);
\draw(-0.75,-0.25) rectangle(3.75,1.25);
\draw(0.8,1) node {\textbullet}; \draw(2.3,1) node {$X_1$};
\draw(3.5,1) node {$n_1$};
\draw(0.5,1) node {$X_0$};
\draw[dotted](2,1)--(2,-1);
\draw[dotted](0.8,1)--(0.8,-1);

\draw[dashed](0,0.5)--(3,0.5);

\draw[->, ultra thick](0.8,1)--(2,1); 

\draw[->, ultra thick](2,0)--(0.8,0); \draw(2,0) node {\textbullet};
\draw(3.5,0) node {$n_2$};

\draw[dashed](0,-0.5)--(3,-0.5);

\draw[->](2,-1)--(1,-1); \draw[->](1,-1)--(0.8,-1);

\draw(3.5,-1) node {$n_3$};

\draw(1.5,-2) node {$n_1= \lfloor \log_{b}k_{on}\rfloor$};
\draw(1.5,-2.5) node {$n_2= \lfloor \log_{b}k_{off}\rfloor$};
\draw(1.5,-3) node {$n_3= \lfloor \log_{b}\nu_+\rfloor$};
\end{scope}


\begin{scope}[shift={(6,-6)}]
\draw(1.5,2) node {$2'$};

\draw(-1,-1.5) rectangle(4,1.5);
\draw(-0.75,-0.25) rectangle(3.75,1.25);
\draw(2,1) node {\textbullet}; \draw(2.3,1) node {$X_1$};
\draw(3.5,1) node {$n_1$};
\draw(0.5,1) node {$X_0$};
\draw[dotted](2,1)--(2,-1);
\draw[dotted](0.8,1)--(0.8,-1);

\draw[dashed](0,0.5)--(3,0.5);

\draw[<-, ultra thick](0.8,1)--(2,1);  \draw(2,1) node {\textbullet};

\draw[<-, ultra thick](2,0)--(0.8,0); \draw(0.8,0) node {\textbullet};
\draw(3.5,0) node {$n_2$};

\draw[dashed](0,-0.5)--(3,-0.5);

\draw[->](2,-1)--(1,-1); \draw[->](1,-1)--(0.8,-1);

\draw(3.5,-1) node {$n_3$};

\draw(1.5,-2) node {$n_1= \lfloor \log_{b}k_{off}\rfloor$};
\draw(1.5,-2.5) node {$n_2= \lfloor \log_{b}k_{on}\rfloor$};
\draw(1.5,-3) node {$n_3= \lfloor \log_{b}\nu_+\rfloor$};
\end{scope}


\begin{scope}[shift={(0,-12)}]
\draw(1.5,2) node {$3$};
\draw(-0.75,-0.25) rectangle(3.75,1.25);
\draw(-1,-1.5) rectangle(4,1.5);
\draw(0.8,1) node {\textbullet}; \draw(2.3,1) node {$X_1$};
\draw(3.5,1) node {$n_1$};
\draw(0.5,1) node {$X_0$};
\draw[dotted](2,1)--(2,-1);
\draw[dotted](0.8,1)--(0.8,-1);

\draw[dashed](0,0.5)--(3,0.5);

\draw[<-, ultra thick](2,1)--(0.8,1);
\draw(3.5,0) node {$n_2$};
\draw[->, ultra thick](2,0)--(1,0); \draw[->, ultra thick](1,0)--(0.8,0); \draw[ultra thick](2,0) node {\textbullet};

\draw[dashed](0,-0.5)--(3,-0.5);
\draw[<-](0.8,-1)--(2,-1);   

\draw(3.5,-1) node {$n_3$};

\draw(1.5,-2) node {$n_1= \lfloor \log_{b}k_{on}\rfloor$};
\draw(1.5,-2.5) node {$n_2= \lfloor \log_{b}\nu_+\rfloor$};
\draw(1.5,-3) node {$n_3= \lfloor \log_{b}k_{off}\rfloor$};
\end{scope}


\begin{scope}[shift={(6,-12)}]
\draw(1.5,2) node {$3'$};

\draw(-1,-1.5) rectangle(4,1.5);
\draw(2,1) node {\textbullet}; \draw(2.3,1) node {$X_1$};
\draw(3.5,1) node {$n_1$};
\draw(0.5,1) node {$X_0$};
\draw[dotted](2,1)--(2,-1);
\draw[dotted](0.8,1)--(0.8,-1);

\draw[dashed](0,0.5)--(3,0.5);

\draw[->, ultra thick](2,1)--(0.8,1);
\draw(3.5,0) node {$n_2$};
\draw[->](2,0)--(1,0); \draw[->](1,0)--(0.8,0);

\draw[dashed](0,-0.5)--(3,-0.5);
\draw[->, ultra thick](0.8,-1)--(2,-1);   \draw(0.8,-1) node {\textbullet};

\draw(3.5,-1) node {$n_3$};

\draw(1.5,-2) node {$n_1= \lfloor \log_{b}k_{off}\rfloor$};
\draw(1.5,-2.5) node {$n_2= \lfloor \log_{b}\nu_+\rfloor$};
\draw(1.5,-3) node {$n_3= \lfloor \log_{b}k_{on}\rfloor$};
\end{scope}

\end{tikzpicture}

\medskip {\bf Fig. \ref{fig:example1}}
\end{center}

\noindent For the sake of the reader, the 1-1 reactions $X_1
\overset{k_{off}}{\to} X_0$ and the split reaction 
$X_1 \overset{2\nu_+}{\to} X_0$ coming from the doubling
reaction $X_1\overset{\nu_+}{\to} 2X_0$ have been represented separately, though they should be understood as a single reaction $X_1\overset{k_{1\to 0}}{\to} X_0$ 
with rate $k_{1\to 0}=k_{off}+2\nu_+$.

\Medskip (\ref{eq:lambda-TypeI}) can be rewritten
\BEQ \log_{b}(\lambda^*)\simeq \log_{b}(k_{on}) + 
\log_{b}(\nu_+) - n_1  \EEQ
where $n_1 = \max(\lfloor\log_{b}(k_{on}) \rfloor,\lfloor\log_{b}(k_{off}) \rfloor,
\lfloor\log_{b}(\nu_+) \rfloor)$ is the maximum scale.  In other words,
\BEQ \lambda^* \approx \begin{cases} k_{on} \qquad ({\mathrm{cases}}\ 1,1') \\
\nu_+\qquad ({\mathrm{cases}}\ 2,3) \\ \frac{k_{on}}{k_{off}}\nu_+
\qquad ({\mathrm{cases}}\ 2',3') \end{cases}  \label{eq:lambda123}
\EEQ

\Medskip Let us see briefly how to retrieve this result
using our Resolvent Formula. Choose $\sigma=1$, then 

\BEQ (f(\alpha))_{1\to 1}= w(\alpha)_{1\to 0} w(\alpha)_{0\to 1}
= \frac{k_{off}+2\nu_+}{k_{off}+\nu_++\alpha} 
\ \times\ \frac{k_{on}}{k_{on}+\alpha}. 
\label{eq:f(alpha)-TypeI}
\EEQ
 In cases 1,1', i.e.
when $n_1=\lfloor\log_{b}(\nu_+)\rfloor$, the largest
kinetic rate is $\nu_+$, therefore the most susceptible factor in
(\ref{eq:f(alpha)-TypeI}) is $\frac{k_{on}}{k_{on}+\alpha}$;   letting $\alpha\approx k_{on}$
yields  $(f(\alpha))_{1\to 1}\approx 2\times \half = 1$. 
In cases 2,3,  i.e. when $n_1=\lfloor\log_{b}(k_{on})\rfloor$,
  the largest
kinetic rate is $k_{on}$,  therefore the most susceptible factor in
(\ref{eq:f(alpha)-TypeI}) is $\frac{k_{off}+2\nu_+}{k_{off}+\nu_++\alpha}$; letting $\alpha\approx \nu_+$ yields 
 $(f(\alpha))_{1\to 1}\approx 1\times 1 = 1$.   In cases $2'$,3', i.e.
when $n_1=\lfloor\log_{b}(k_{off})\rfloor$, we take out relevant terms in
(\ref{eq:f(alpha)-TypeI}) by Taylor expanding to order 1; letting
$\eps = \frac{\nu_+}{k_{off}+\nu_+} \sim \frac{\nu_+}{k_{off}}$, 
we get
$(f(\alpha))_{1\to 1} = (1+\eps) \frac{k_{off}+\nu_+}{k_{off}+\nu_++\alpha} \ \times\ \frac{k_{on}}{k_{on}+\alpha} \sim 
(1+\frac{\nu_+}{k_{off}}) (1-\frac{\alpha}{k_{on}}) \sim 
1+ \frac{\nu_+}{k_{off}} - \frac{\alpha}{k_{on}}$. Thus
$(f(\alpha))_{1\to 1}= 1$ for $\alpha \sim \frac{k_{on}}{k_{off}}\nu_+$.



\section{Example 2}  \label{section:example2}


 We enlarge the network of Example 1 by including a  reaction $X_0\to X_1$ proceeding through an 
intermediate $X_4$,
\BEQ  X_0 \overset{k'_4}{\to} X_4, \qquad X_4\overset{k_4}{\to}
X_1  \label{eq:6.1} \EEQ
The above may modelize a hidden catalytic reaction in presence of an enzyme $E$ in excess, such as 
$X_0 + E \to X_4 \to X_1 + E'$, where $E'=E$, or $E'\not=E$, if e.g. a cofactor bound to $E$ has 
been taken away by $X_0$.

\Medskip There are many possible scale-splittings of this network, 
including (extending the scale-splitting, case 3' of Example 1)

\medskip
\begin{center}
\begin{tikzpicture}[scale=0.8]
\draw[dotted](0,0)--(0,-4); \draw[dotted](1,-4)--(1,-5);
\draw[dotted](2,0)--(2,-5);
\draw[<-, ultra thick](0,0)--(2,0);  \draw(-0.3,0) node {$0$};
\draw(1,0.3) node {$k_{off}$};
\draw(2.3,0) node {$1$};  \draw(2,0) node {\textbullet};
\draw(5,0) node {$n_1$};

\draw[dashed](-1,-0.5)--(4,-0.5);
\draw[<-](0.2,-1)--(2,-1);
\draw[<-](0,-1)--(0.2,-1);
\draw(1,-0.7) node {$\nu_+$};
\draw(5,-1) node {$n_2$};

\draw[dashed](-1,-1.5)--(4,-1.5);
\draw[->, ultra thick](0,-2)--(2,-2); \draw(-0,-2) node {\textbullet};
\draw(1,-1.7) node {$k_{on}$};
\draw(5,-2) node {$n_3$};

\draw[dashed](-1,-2.5)--(4,-2.5);

\draw[dashed](-1,-4.5)--(4,-4.5);
\draw[->](0,-4)--(1,-4);  \draw[->, ultra thick](1,-5)--(2,-5);
\draw(1.5,-4.7) node {$k_4$}; \draw(0.5,-3.7) node {$k'_4$};
\draw(1,-5) node {\textbullet};
\draw(1-0.3,-5) node {$4$};
\draw(5,-4) node {$n'_4 = n_{0\to 4}$}; \draw(5,-5) node {$n_4$};
\end{tikzpicture}
\end{center}
\medskip
with $n'_4 = \lfloor log_{b} k'_4 \rfloor, n_4 = \lfloor
\log_{b} k_4 \rfloor$.  We prefer to use the notation 
$n_{0\to 4}$ instead of $n'_4$ because it is not a vertex 
scale, but the scale of the transition from $X_0$ to $X_4$.

\Medskip The analysis below is generalized in \S \ref{subsection:threshold} and \ref{subsection:Regime}
below, where the same notations are used. Thus the reader eager to get to the results may skip some
details, and concentrate either on this particular case or read the arguments in Section \ref{section:sigma*}.  

\Medskip {\em Two regimes.} Since the Lyapunov exponent is 
upper-bounded by the largest autocatalytic rate (see Lemma \ref{lem:1}), we know that $\lambda^*\lesssim
k_{on}$, so we need only consider $\alpha\preceq k_{on}$.   There are two loops contributing
to $(f(\alpha))_{1\to 1}$:  
$\gamma_1 : 1\overset{k_{off}+2\nu_+}{\to} 0\overset{k_{on}}{\to} 1$ and $\gamma_2 :   1\overset{k_{off}+2\nu_+}{\to} 0\overset{k'_4}{\to} 4 \overset{k_4}{\to} 1 $.  Letting 
\BEQ \eps_1 := \frac{\nu_+}{k_{off}+\nu_+}\sim \frac{\nu_+}{k_{off}},
\EEQ
be the deficiency weight of $X_1$, see (\ref{eq:deficiency-weight}), we have  
\BEQ (w(\alpha))_{1\to 0} = \frac{k_{off}+2\nu_+}{k_{off} +
\nu_++\alpha} = (1+\eps_1) \frac{k_{off}+\nu_+}{k_{off} +
\nu_++\alpha} \sim  (1+\eps_1) (1-\frac{\alpha}{k_{off}}), 
\EEQ
\BEQ (w(\alpha))_{0\to 1} = \frac{k_{on}}{k_{on}+k'_4+\alpha}
\sim (1-\frac{k'_4}{k_{on}})(1-\frac{\alpha}{k_{on}}),
\EEQ
from which
\BEQ (w(\alpha))_{\gamma_1} =  (w(\alpha))_{1\to 0}(w(\alpha))_{0\to 1} \sim 
(1+\eps_1)(1-\frac{k'_4}{k_{on}})(1 - \frac{\alpha}{k_{on}});
\EEQ 

\BEQ (w(\alpha))_{0\to 4} = \frac{k'_4}{k_{on}+k'_4+\alpha} \sim
\frac{k'_4}{k_{on}} (1-\frac{\alpha}{k_{on}})
\EEQ
and
\BEQ 
 (w(\alpha))_{4\to 0} = \frac{k_4}{\alpha+k_4} \sim \begin{cases} 
 \frac{k_4}{\alpha},  \qquad \alpha \succ k_4 \\   1-\frac{\alpha}{k_4},
 \qquad \alpha\prec k_4 \end{cases}
\EEQ
By comparison, $(\tilde{w}(0))_{1\to 0}=1,\ (\tilde{w}(0))_{0\to 1} =  \frac{k_{on}}{k_{on}+k'_4}
\sim 1-\frac{k'_4}{k_{on}}, \  (\tilde{w}(0))_{0\to 4} = \frac{k'_4}{k_{on}+k'_4}\sim \frac{k'_4}{k_{on}}, \ 
(\tilde{w}(0))_{4\to 0}=1$. 
Thus we must distinguish between the two scale windows, (A) $\alpha\succ k_4$, and  (B) $\alpha\prec k_4$.
Solving the threshold equation leads to introducing scale-window dependent regimes and sub-regimes.


\Bigskip {\bf (A) Assume first that $k_4 \prec \alpha\prec k_{on}$.} Then
\BEQ (f(\alpha))_{1\to 1} \sim c_{high} - \alpha c'_{high} + 
\frac{c_1(\alpha)}{\alpha}  \label{eq:TypeI-1}
\EEQ
with
\BEQ c_{high} = (w(0))_{\gamma_1} \sim (1+\eps_1)(1-\frac{k'_4}{k_{on}}), \qquad c'_{high} = - \frac{d}{d\alpha}  (w(\alpha))_{\gamma_1} \Big|_{\alpha=0} \sim  k_{on}^{-1};
\EEQ
\BEQ c_1(\alpha) = (w(\alpha))_{1\to 0} (w(\alpha))_{0\to 4} k_4 
\sim  (1+\eps_1) (1-\frac{\alpha}{k_{on}})  \frac{k'_4}{k_{on}}
k_4\sim  \frac{k'_4}{k_{on}} k_4;
\EEQ

\BEQ   \del c_{high}  := c_{high}-1 \sim 
\eps_1 - \frac{k'_4}{k_{on}}  \label{eq:TypeI-deltac-high}.
\EEQ
The threshold value of $\alpha$ for which $(f(\alpha))_{1\to 1}=1$
is such that 
\BEQ {\mathrm{(Threshold\ equation)}} \qquad \del c_{high} - \alpha c'_{high} + \frac{c_1(\alpha)}{\alpha}=0.
\label{eq:TypeI-cc1}
\EEQ

 Depending on the sign of $\del c_{high}$, or equivalently, on
whether $\eps_1/(k'_4/k_{on})\gtrless 1$,  we get two regimes in which $\eps_1/(k'_4/k_{on})\succ 1$,
resp. $\eps_1/(k'_4/k_{on})\prec 1$, plus a delicate {\em transition regime} where $\eps_1 \sim \frac{k'_4}{k_{on}}$, which falls out of our discussion. We call  Regime (1), resp.
Regime (2), the case when $\del c_{high}>0$, resp. $\del c_{high}<0$.
Introduce the ratio
\BEQ R:= (\frac{k_{on}}{k_{off}}\nu_+)/k'_4. \label{eq:TypeI-R}
\EEQ
Then $R\gg 1$ in Regime (1), and $R\ll 1$ in Regime (2).

\Medskip
{\bf Assume first we are in Regime (1),} 
\BEQ (1)\qquad \frac{k_{on}}{k_{off}}\nu_+ \gg k'_4 \EEQ
Then  $\del c_{high}\sim \eps_1 \approx 
 \frac{\nu_+}{k_{off}} >0$. In (\ref{eq:TypeI-cc1}), we get
 two positive terms, and one negative term, therefore, two sub-regimes. Solving the 
 equation is possible if either (a) $\del c_{high} \sim \alpha c'_{high}$ and $\frac{c_1(\alpha)}{\alpha}\ll \del c_{high}$, or  (b)
 $\alpha c'_{high} \sim \frac{c_1(\alpha)}{\alpha}$ and $\del c_{high}
 \ll \frac{c_1(\alpha)}{\alpha}$. Case (b) is not possible
 since it yields $\alpha\sim \sqrt{k'_4 k_4}$, but 
 $\del c_{high} /(\frac{c_1(\alpha)}{\alpha}) \sim  R \sqrt{\frac{k'_4}{k_4}}
 \gg 1$. Case (a), on the other hand, is possible, yielding
 \BEQ \alpha \sim \frac{k_{on}}{k_{off}}\nu_+ \label{eq:2.36}
  \EEQ
(which is in the requested interval $[k_4,k_{on}]$ since
 $\frac{k_{on}}{k_{off}}\nu_+ = k'_4 R\gg k'_4>k_4$); namely, $(\frac{c_1(\alpha)}{\alpha})/ \del c_{high} = k'_4 k_4
(\frac{k_{off}}{k_{on}\nu_+})^2 = (\frac{k'_4}{k_4} R^2)^{-1}
\ll 1$. 

\Medskip  Thus (\ref{eq:2.36}) (Sub-regime (1) (a)) gives the looked-for estimate
of $\lambda^*$ when $R\gg 1$ (Regime (1)): we have found
\BEQ \lambda^* \simeq \frac{k_{on}}{k_{off}}\nu_+ \qquad {\mathrm{if}}\ k'_4\prec \frac{k_{on}}{k_{off}}\nu_+.
\EEQ

\Medskip
{\bf Assume now we are in Regime (2),} i.e. $R\ll 1$, equivalently,
$ (2)\qquad \frac{k_{on}}{k_{off}}\nu_+ \ll k'_4$,
from which $\del c_{high} \sim -\frac{k'_4}{k_{on}}<0$. 
Assuming still that $\alpha\gg k_4$, (\ref{eq:TypeI-cc1})
implies either (a) $|\del c_{high}|\sim \frac{c_1(\alpha)}{\alpha}$ 
and $\alpha c'_{high}\ll |\del c_{high}|$, or (b) $\alpha c'_{high}
\sim \frac{c_1(\alpha)}{\alpha}$ and $|\del c_{high}|\ll \alpha 
c'_{high}$. Neither case is valid, since (a) $\alpha\sim k_4$,
a border case  not compatible with the hypothesis $ \alpha\gg
k_4$; and (b) $\alpha\sim \sqrt{k_4 k'_4}$, but then 
$|\del c_{high}|/(\alpha c'_{high}) \sim
\sqrt{\frac{k'_4}{k_4}}\gg 1$. 
 The conclusion is: Regime (2) is excluded.

\Bigskip {\bf (B) Assume next that $\alpha\ll k_4$.} We already know that $R\ll 1$.  Then, instead
of (\ref{eq:TypeI-1}), we have
\BEQ (f(\alpha))_{1\to 1} \sim 1+\del c_{high}-\alpha c'_{high} \EEQ
with now 
\BEQ \del c_{high}  = -1+ (w(0))_{\gamma_1} + (w(0))_{\gamma_2} = \eps_1,
\EEQ
\BEQ
 c'_{high} =  -\frac{d}{d\alpha} ((w(\alpha))_{\gamma_1} + (w(\alpha))_{\gamma_2})\Big|_{\alpha=0} \sim  \frac{1}{k_{on}} + \frac{k'_4/k_{on}}{k_4} \sim \frac{1}{k_{on}} \frac{k'_4}{k_4}.
\EEQ
Thus $(f(\alpha))_{1\to 1}=1$ for $\alpha$ such that 
$\alpha\sim \frac{\del c_{high}}{c'_{high}} = Rk_4$. Note that 
$\alpha\ll k_4$  since
$R\ll 1$ by assumption. We have found:
\BEQ \lambda^* \simeq \frac{k_4}{k'_4} \frac{k_{on}}{k_{off}}\nu_+ 
\qquad {\mathrm{if}}\ k'_4\succ \frac{k_{on}}{k_{off}}\nu_+. \EEQ


\Bigskip We illustrate the above computations in the log-log
 plot  below with (choosing a base constant $b=10$)
  $k_{off} = 1, \nu_+ = 10^{-2}, k_{on} = 10^{-3}$, 
 $\frac{k_{on}}{k_{off}}
 \nu_+= 10^{-5}$, $k_4=10^{-8}$, and 
 $k'_4$ varying between $10^{-2} \, \times \frac{k_{on}}{k_{off}}
 \nu_+= 10^{-7}$ and $10^{2} \, \times \frac{k_{on}}{k_{off}}
 \nu_+= 10^{-3}$. In yellow, the true Lyapunov exponent determined
 by numerical diagonalization. The agreement between theory and numerics is manifest, save at the transition between Regime (1)  where $R\gg 1$
 ($k'_4<10^{-5}$, curve $\alpha \sim \frac{k_{on}}{k_{off}}\nu_+$ in red) and Regime (2) where $R\ll 1$ 
  ($k'_4>10^{-5}$, curve $\alpha\sim  Rk_4$  in blue). The region defined by $R\approx 1$ is called a 
  {\em resonance}. In Appendix, we prove that at the transition, $\lambda^*\sim\lambda^0$, with
  $\lambda^0 = \sqrt{ \frac{k_{on}}{k_{off}}\nu_+ k_4}$, a value materialized by
  the green dashed line, which is the geometric mean between
 the estimates for $\lambda^*$ obtained in Regime (1) and in Regime (2). We can furthermore get 
 the slope at the resonance on the log-log plot,  $\frac{\del \lambda^*}{\lambda^*} \sim -
 \sqrt{\frac{k'_4}{k_4}} \ \frac{\del k'_4}{k'_4}$, a large but finite number (the log-log plot looks
 non-differentiable, square-root like, around the transition value, but it is not).

\Medskip   Yellow and blue curves part again when 
  $k'_4\to k_{on}=10^{-3}$ and the scale ordering on the multi-scale diagram below
  (\ref{eq:6.1})  changes. 

\begin{figure}[!t]
\centering
\includegraphics[scale=0.6]{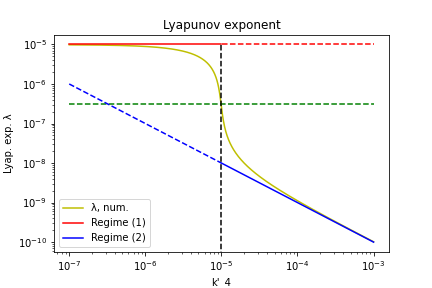}
\caption{Lyapunov exponent as a function of kinetic parameter
$k'_4$}
\label{Fig:1}
\end{figure}


\section{$\sigma^*$-dominant graphs}  \label{section:sigma*}


Dominant cycles (see Section \ref{section:multi-scale-splitting}) play a major role in our multi-scale
algorithm. Examples 1,2 above, and 3 below, are all in a subclass of graphs which is specific, in the
sense that they admit a distinguished vertex $\sigma^*$ which is necessarily involved in all
autocatalytic cycles.  {\em We assume in the following that $\lambda^*>0$.} 

\begin{Definition}[$\sigma^*$-dominant graph]  \label{def:sigma*-dominant}
 A graph $G$ is  $\sigma^*$-dominant if it admits
 a distinguished vertex, denoted $\sigma^*\in V$, such that 
  the $\sigma^*$-dominant Hypothesis below holds,

\begin{center}
($\sigma^*$-dominant Hypothesis) there exists a dominant cycle $\gamma^*:\sigma^*=\sigma_1\to 
\cdots\to \sigma_{\ell^*} = \sigma^*$, and  all dominant cycles pass through $\sigma^*$
\end{center}

\end{Definition}

 The  $\sigma^*$-dominant Hypothesis will make it possible to 'solve' the implicit equation following from the Resolvent Fromula (see Proposition \ref{prop:main-formula})
for $\lambda^*$ if $\lambda^*>0$,, namely, to get the scale of the threshold $\alpha_{thr}=\lambda^*$. The {\bf monotonicity
property} of 
$\alpha \mapsto f(\alpha)_{\sigma^*\to\sigma^*}$, see Section \ref{section:resolvent} and first lines
of Section \ref{section:multi-scale-splitting}, allows a 'trial and error' method, where one starts
with a very large scale $n_{\alpha}$ for which $f(\alpha)_{\sigma^*\to\sigma^*}<1$,  and lowers it step by 
step until  
$f(\alpha)_{\sigma^*\to\sigma^*}>1$.

\Bigskip A {\em simple excursion from $\sigma^*$ to $\sigma^*$} is an excursion $\gamma: 
\sigma^*=\sigma_1\to\cdots\to \sigma_{\ell}=\sigma^*$ s.t. all intermediate vertices are distinct, 
$(1\le i<j\le \ell-1)\Rightarrow (\sigma_i\not=\sigma_j)$. A {\em dominant excursion  from $\sigma^*$ to $\sigma^*$} is an excursion $\gamma: 
\sigma^*=\sigma_1\to\cdots\to \sigma_{\ell}=\sigma^*$ which is a dominant cycle. Let now $\sigma\not=\sigma^*$; by the $\sigma^*$-dominant Hypothesis, there exists no dominant excursion from $\sigma^*$ to $\sigma^*$. Let however 
\BEQ n_{\sigma^*\twoheadrightarrow \sigma} := \min_{\gamma:\sigma^*=\sigma_1\to\cdots\to 
\sigma_{\ell} = \sigma} \sum_{i=1}^{\ell-1} \Big\{ n(\sigma_i)-n(\sigma_i\to\sigma_{i+1})\Big\}
\label{eq:nsigma*totosigma}
\EEQ
(note the double-headed arrow, distinguishing $n_{\sigma^*\twoheadrightarrow\sigma}$ from $n_{\sigma^*\to\sigma}$, and compare with (\ref{eq:path-weight})); this scale index is approximately equal to minus the log 
of a leading path from $\sigma^*$ to $\sigma$.   Let $\lambda^*$ be the Lyapunov exponent, and 
$f^*_{\sigma^*\to\sigma} = f(\lambda^*)_{\sigma^*\to\sigma}$, $\sigma=\sigma^*$ or $\not=\sigma^*$; 
by the resolvent identity, $f^*_{\sigma^*\to\sigma^*}\equiv 1$.   We summarize the main properties implied by 
the $\sigma^*$-dominant Hypothesis, which will prove crucial in \S \ref{subsection:Lyapunov-eigenvector}:

\begin{Theorem}  \label{th:sigma*-dominant}
Assume that the $\sigma^*$-dominant Hypothesis of Definition \ref{def:sigma*-dominant} holds. Then,
if $b$ is large enough: 
\begin{itemize}
\item[(i)] there exists a {\em dominant simple excursion from $\sigma^*$ to $\sigma^*$}, namely a dominant
excursion $\gamma^*_{\sigma^*}:\sigma^*=\sigma_1\to\cdots\to \sigma_{\ell^*}=\sigma^*$ from $\sigma^*$ to $\sigma^*$ s.t. $w(\lambda^*)_{\gamma^*_{\sigma^*}}\sim 1$;
\item[(ii)] for every $\sigma\not=\sigma^*$, there exists a {\em leading simple excursion from $\sigma^*$ to $\sigma$}, namely a leading
excursion $\gamma^*_{\sigma}:\sigma^*=\sigma_1\to\cdots\to \sigma_{\ell^*}=\sigma$ from $\sigma^*$ to $\sigma$ s.t. $\sum_{i=1}^{\ell-1} \Big\{ n(\sigma_i)-n(\sigma_i\to\sigma_{i+1})\Big\} = n_{\sigma^*\twoheadrightarrow \sigma}$; furthermore, $f^*_{\sigma^*\to\sigma} \sim w(\lambda^*)_{\gamma^*_{\sigma}}$;
\item[(iii)] let $\omega^{*,\sigma} := \sum_{\ell> 1}\sum_{\gamma:\sigma=\sigma_1\to\cdots,\sigma_{\ell}=\sigma, \ 
\sigma_i\not=\sigma^*, i=1,\ldots,\ell} 
w(\lambda^*)_{\gamma}$ be the total weight of non-trivial paths from $\sigma$ to $\sigma$ 
{\em that do not pass through $\sigma^*$}. Then $\omega^{*,\sigma} \prec 1$.
\item[(iv)] let $\omega^*_n$ be the total weight $\sum_{\gamma:\sigma^*\to\cdots\to\sigma^*, \ 
\ell(\gamma)\ge n} w(\lambda^*)_{\gamma}$ of paths of length $n$ from 
$\sigma^*$ to $\sigma^*$. Then $\min_{n\ge 1} (\omega^*_n+\cdots + \omega^*_{n+\ell^*-1})\succeq 1$. 
\end{itemize}
\end{Theorem}

The proof is postponed to Appendix. Note that (see discussion of transition weights in 
Section \ref{section:multi-scale-splitting}), if $\alpha_{thr}\sim \lambda^*> 0$, then  $w(0)_{\gamma^*}\sim w(\alpha_{thr})_{\gamma^*}\sim 
w(\lambda^*)_{\gamma^*}\sim 1$ and  $w(\alpha_{thr})_{\gamma^*_{\sigma}}\sim 
w(\lambda^*)_{\gamma^*_{\sigma}} \preceq w(0)_{\gamma^*_{\sigma}}$.
The latter approximate inequality $\preceq$ is strict if $n_{\alpha_{thr}}$ is larger than 
one of the scales of the vertices along $\gamma^*_{\sigma}$.  Similarly, $f^*_{\sigma^*\to\sigma}\le f(0)_{\sigma^*\to\sigma}$. Finally, $w(\lambda^*)_{\gamma}\le w(0)_{\gamma}$ (actually, the
sum in (iv) diverges for $\lambda>0$ if one replaces $w(\lambda^*)_{\gamma}$ by $w(0)_{\gamma}$ since
$f(0)_{\sigma^*\to\sigma^*}>1$).


\subsection{Threshold equation from Resolvent Formula} \label{subsection:threshold}


\Bigskip {\em Rewriting the Resolvent Formula (\ref{eq:main-formula}).} Let
$x(\sigma) =  \frac{|A_{\sigma,\sigma}|+\kappa_{\sigma}}{
|A_{\sigma,\sigma}| + \alpha}
$
be the ratio between the denominators in 
$(w(\alpha))_{\sigma\to\sigma'} = \frac{A_{\sigma',\sigma}}{|A_{\sigma,\sigma}|+\alpha}$ and in $(\tilde{w}(0))_{\sigma\to\sigma'} = \frac{A_{\sigma',\sigma}}{|A_{\sigma,\sigma}|+\kappa_{\sigma}}$. 
Then 
\BEQ f(\alpha)_{\sigma\to\sigma} = \sum_{\gamma:\sigma\to
\sigma} w(\alpha)_{\gamma}  =  \sum_{\gamma:\sigma\to
\sigma}   \tilde{w}(0)_{\gamma}  \prod_{i=1}^{\ell(\gamma)-1} 
x(\sigma_i) =
\langle \sum_{\gamma:\sigma\to\sigma} \prod_{i=1}^{\ell(\gamma)-1} x(\sigma_i) \rangle_{\sigma}
\EEQ
where the average $\langle \ \cdot \ \rangle_{\sigma}$ is w.r. to the probability
measure $(\tilde{w}(0)_{\gamma})_{\gamma}$ over excursions
$\gamma$ from $\sigma$ to $\sigma$.

\Bigskip {\em Taylor expansion of $x(\sigma)$.}
Recall the notations $a_{\sigma}:=|A_{\sigma,\sigma}|$ for the 
diagonal entries (see 
(\ref{eq:a})) and  $\eps_{\sigma}=\frac{\kappa_{\sigma}}{a_{\sigma}} $, see (\ref{eq:deficiency-weight}), for
deficiency  weights; we have already prove that $\eps_{\sigma}\le 1$, and that $|\eps_{\sigma}|\prec 1$ 
if $\eps_{\sigma}<0$. Thus,  the product of factors $\prod_{i=1}^{\ell-1} (1+\eps_{\sigma_i})$ along an excursion $\gamma:\sigma=\sigma_1\to\cdots\to\sigma_{\ell}=\sigma$ may be rewritten as 
$1+\eps_{\gamma}$, with $\eps_{\gamma}\approx \sum_{i=1}^{\ell-1} \eps_{\sigma_i} \preceq 1$.  Then
\BEQ x(\sigma)\sim \begin{cases}   
 (1+\eps_{\sigma})
(1-\frac{\alpha}{a_{\sigma}}) \qquad \alpha\ll a_{\sigma} \\
 (1+\eps_{\sigma})
\frac{a_{\sigma}}{\alpha} \qquad  \alpha\gg a_{\sigma} \end{cases}
\label{eq:x(sigma)}
\EEQ
This defines two regimes  (small vs.
large $\alpha$) that depend on $\sigma$. 

\Medskip {\em Some further notations.}  If $\gamma:\sigma_1\to\cdots\to\sigma_{\ell}$ is a path or a union of paths (not
necessarily an excursion),
we let
\BEQ a_{\gamma}:= \min\Big(a_{\sigma_i}; \ i=1,\ldots,\ell\Big)
\label{eq:a-gamma}
\EEQ
and define
\BEQ n_{\gamma}:= \lfloor \log_{b} a_{\gamma}\rfloor
\label{eq:n-gamma}
\EEQ
to be the {\em scale} of $\gamma$; it can also be defined
as the minimum of the scales of all vertices along $\gamma$.
Then
\BEQ \Pi_{\gamma}:=  \prod_{i=1}^{\ell} a_{\sigma_i}
\label{eq:Pi-gamma}
\EEQ
is the product of all $a_{\sigma}$ along $\gamma$; 
the corresponding scale is 
\BEQ N_{\gamma}:=\lfloor 
\log_{b} \Pi_{\gamma}\rfloor. \label{eq:N-gamma}
\EEQ

\Medskip {\em Remark.}
 Assume first $\alpha\gg b^{n_1}$. Then we are in the large $\alpha$
scale regime for all $\sigma$, hence (for every $\sigma$)
$ (f(\alpha))_{\sigma\to\sigma} \sim \langle  \prod_{i=1}^{\ell(\gamma)-1} \frac{a_{\sigma_i}}{\alpha} \rangle_{\gamma} <1,
$
since the averaged quantity $\prod_{i=1}^{\ell(\gamma)-1} \frac{a_{\sigma_i}}{\alpha}$ is $o(1)$. 
 Therefore, $\lambda^*\lesssim b^{n_1}$, and we may
 assume that $n_{\alpha} <n_1$.

\Bigskip {\em Scale window.} A scale window is an integer
interval $\{n+1,\ldots,n'-1\}$, where $n=n_{\sigma},\, n'=n_{\sigma'}$
are consecutive vertex scales, i.e. $n_{\sigma}<n_{\sigma'}$
and there exists no vertex scale $n_{\sigma''}$ such that
$n_{\sigma}<n_{\sigma''}<n_{\sigma'}$.

\Bigskip {\em Fix a scale window $\{n+1,\ldots,n'-1\}$,  and a scale 
 $n_{\alpha} \in \{n+1,\ldots,n'-1\}$}.
An excursion $\gamma:\sigma^*=\sigma_1\to\cdots\to\sigma_{\ell}=\sigma^*$ is {\em high} if all $n_{\sigma_i}>n_{\alpha}$ (equivalently, if $\alpha$ is small along $\gamma$). 
All other excursions are {\em low}; more precisely,  we say
that
$\gamma$ is $k$-low if the number of indices $i=1,\ldots,\ell-1$
such that $n_{\sigma_i}<n_{\alpha}$ is equal to $k$. We split a general
excursion $\gamma:\sigma^*\to\sigma^*$ into $\gamma^{{\mathrm{high}}}\cup\gamma^{{\mathrm{low}}}$,
where 
\BEQ \gamma^{{\mathrm{high}}}= \{(\sigma_i)_{i=1,\ldots,\ell-1}\ |\ n_{\sigma_i}>n_{\alpha} \}, \qquad \gamma^{{\mathrm{low}}}= \{ (\sigma_i)_{ i=1,\ldots,\ell-1}\ |\ n_{\sigma_i}<n_{\alpha} \},
\EEQ
 considered as unions of paths
by connecting vertices of $\gamma^{{\mathrm{high}}}$ or $\gamma^{{\mathrm{low}}}$ if they are
connected along $\gamma$. Note that $\gamma$ is $k$-low
if and only if the length $\ell(\gamma^{{\mathrm{low}}})$ of  $\gamma^{{\mathrm{low}}}$ 
is equal to $k$. Then $(f(\alpha))_{\sigma\to\sigma}$ splits into the sum of two
terms, a high one, and a low one,
\BEA (f(\alpha))_{\sigma^*\to\sigma^*} &=& 
\Big(\sum_{\gamma \ {\mathrm{high}}} w(\alpha)_{\gamma} 
\Big) + 
\Big(\sum_{\gamma \ {\mathrm{low}}} w(\alpha)_{\gamma} 
\Big) \nonumber\\
&=& 
\Big(1+ \del c_{high} - \alpha c'_{high} \Big) +
\Big(   \sum_{k\ge 1} 
c_{k}(\alpha) \alpha^{-k} \Big) \label{eq:cc}
\EEA
with
\BEQ \del c_{high} = -1+ \sum_{\gamma: \sigma^*\to \sigma^*\ {\mathrm{high}}}
(1+\eps_{\gamma}) \, \tilde{w}(0)_{\gamma};
\EEQ
\BEA c'_{high} &=& -\frac{d}{d\alpha} \Big(\sum_{\gamma: \sigma^*\to \sigma^*\ {\mathrm{high}}}
 \, w(\alpha)_{\gamma}\Big)|_{\alpha=0}
\nonumber\\
&\sim & \sum_{\gamma: \sigma^*\to \sigma^*\ {\mathrm{high}}}
(1+\eps_{\gamma})\  \tilde{w}(0)_{\gamma}\  (\sum_{i=1}^{\ell-1} 
\frac{1}{a_{\sigma_i}})
\sim   \sum_{\gamma: \sigma^*\to \sigma^*\ {\mathrm{high}}} 
\frac{ \tilde{w}(0)_{\gamma}}{a_{\gamma}}
\EEA
where $a_{\gamma}$ is as in (\ref{eq:a-gamma}); and, for 
$k\ge 1$, (leaving out factors of order 1 such as $1+\eps_{\gamma}\sim 1$ and  the product of factors $1-\frac{\alpha}{a_{\sigma}}\sim 1$ along $\gamma^{{\mathrm{high}}}$)
$ c_{k}(\alpha) \sim \sum_{\gamma: \sigma^*\to \sigma^*\  k-{\mathrm{low}}}   \, \tilde{w}(0)_{\gamma} \ 
\Pi_{\gamma^{{\mathrm{low}}}},
$
where $\Pi_{\gamma^{{\mathrm{low}}}}$ is as in (\ref{eq:Pi-gamma}).

\Medskip The threshold value of $\alpha$ for which 
$(f(\alpha))_{1\to 1}=1$ is obtained by solving an equation
for the threshold similar to  (\ref{eq:TypeI-cc1}) for  Example 2, namely,

\begin{center}  
\BEQ {\mathbf{(Threshold\ equation)}}\ \ 
\del c_{high} - \alpha c'_{high} +   \sum_{k\ge 1} 
c_{k}(\alpha) \alpha^{-k}=0  \label{eq:threshold} 
\EEQ
\end{center}

\noindent
Therefore,  we  get a priori two regimes:

\begin{itemize}
\item[\textbullet] In Regime (1), called {\em low $\alpha$-regime}, 
\BEQ  (1) \qquad \del c_{high}>0  \EEQ
and the threshold value can be expected to be obtained by
compensating the largest term among $\del c_{high}$ and
$(\frac{c_k(\alpha)}{\alpha^k})_{k\ge 1}$ by $\alpha c'_{high}$.
More precisely, in case (a), $\del c_{high}\sim \alpha c'_{high}$
and $\max_{k\ge 1} \Big(c_k(\alpha) \alpha^{-k}\Big) \ll 
\del c_{high}$. In case $(b_k)$ $(k\ge 1)$,  
$c_k(\alpha)\alpha^{-k} \sim \alpha c'_{high}$ and 
$\del c_{high},\max_{j\not=k} \Big(c_j(\alpha) \alpha^{-j}\Big) \ll c_k(\alpha)\alpha^{-k}$.

\item[\textbullet]    On the other hand, in Regime (2),
called {\em high $\alpha$-regime}, 
\BEQ  (1) \qquad \del c_{high}<0  \EEQ
and the threshold value can be expected to be obtained by
compensating\\
  $\max(|\del c_{high}|,\alpha c'_{high})$ by the
maximum of the 
$(\frac{c_k(\alpha)}{\alpha^k})_{k\ge 1}$. More precisely, in 
case $(a_k)$, $k\ge 1$,  $|\del c_{high}|\sim c_k(\alpha) \alpha^{-k}$ and $\alpha c'_{high}, \max_{j\not=k} \Big(c_j(\alpha) \alpha^{-j}\Big) \ll c_k(\alpha)\alpha^{-k}$.
In 
case $(b_k)$, $k\ge 1$,  $\alpha c'_{high}\sim c_k(\alpha) \alpha^{-k}$ and $\del c_{high}, \max_{j\not=k} \Big(c_j(\alpha) \alpha^{-j}\Big) \ll c_k(\alpha)\alpha^{-k}$.

\end{itemize}

The reason why Regime (1) is called low $\alpha$-regime is
that when $\alpha=0$, all excursions are high, and the autocatalytic assumption $f(0)_{\sigma^*\to
\sigma^*}>1$ 
implies $\del c_{high}>0$. On the contrary, when $\alpha\to +\infty$, all excursions are low, and $\del c_{high}=-1<0$.

\Bigskip In the next subsection, we shall prove the following result:

\begin{Theorem}[threshold parameter and approximate Lyapunov exponent]  \label{th:thr}
Away from transition regimes, $\lambda^* \simeq \alpha_{thr}$, with
\BEQ  \alpha_{thr} = \del c_{high}/c'_{high}. \label{eq:th:thr} \EEQ  
\end{Theorem}


\subsection{Only Regime (1) (a) is relevant} \label{subsection:Regime}


We  prove here that only Regime (1) (a) is possible outside
transition regimes. Thus low path contributions $c_k(\alpha)$
are irrelevant. This implies directly Theorem \ref{th:thr}. 

\Medskip {\em  Simplified threshold equation.} We
solve a simplified version of  the threshold 
equation (\ref{eq:threshold}) $\del c_{high}-\alpha c'_{high}
+\sum_{\gamma \ {\mathrm{low}}} w(\alpha)_{\gamma}=0$, in which the contribution $\sum_{\gamma \ {\mathrm{low}}} w(\alpha)_{\gamma}\sim \sum_{k\ge 1} c_k(\alpha) \alpha^{-k}$ of
low-lying paths has been eliminated,

\begin{center}
\BEQ
{\mathbf{(Simplified\ threshold\ equation)}} \qquad  \del c_{high} - \alpha c'_{high} = 0
\EEQ
\end{center}

\noindent
The solution is trivially
\BEQ \alpha_{thr} = \del c_{high}/c'_{high}, \EEQ
that is, (\ref{eq:th:thr}), 
and it is accepted if and only if $\alpha_{thr}$ is in the scale window
$\{n+1,\ldots,n'-1\}$. 
Decreasing $\alpha$ (which means changing scale window as
one passes vertex scales), the function $\alpha\mapsto\del c_{high}-\alpha c'_{high}$ increases
by the monotonicity argument from a strictly negative 
value for $n_{\alpha} \succ n_1$ to a strictly positive value
as $\alpha\to 0$.  Thus the simplified threshold equation
has a unique solution, fixing in particular a 
threshold scale window  $\{n_{thr}+1,\ldots,n'_{thr}-1\}$ for which $\alpha_{thr} = \del c_{high}/c'_{high}$. In particular, 
$\del c_{high}>0$ in this window. 

\Medskip {\em  Transition regime.}  Since 
$\sum_{\gamma} \tilde{w}(0)_{\gamma}=1$, we have
$\del c_{high} \sim  \sum_{\gamma \ {\mathrm{high}}} 
\eps_{\gamma} \tilde{w}(0)_{\gamma} -\sum_{\gamma \ {\mathrm{low}}} 
\tilde{w}(0)_{\gamma};$ a particular case of this decomposition is (\ref{eq:TypeI-deltac-high}). 
The transition regime, which falls out of our discussion, is defined by the parameter constraint (or {\em resonance}) $ \sum_{\gamma \ {\mathrm{high}}} 
\eps_{\gamma} \tilde{w}(0)_{\gamma} \sim \sum_{\gamma \ {\mathrm{low}}} 
\tilde{w}(0)_{\gamma}$. Assuming we are out of this transition regime, we either have $\del c_{high} \sim 
\sum_{\gamma \ {\mathrm{high}}} 
\eps_{\gamma} \tilde{w}(0)_{\gamma} \succ \sum_{\gamma \ {\mathrm{low}}} 
\tilde{w}(0)_{\gamma} >0$, or $\del c_{high}<0$. 

\Medskip {\em Comparison to original threshold equation.}
Let $\alpha\ll \alpha_{thr}$ be smaller than $\alpha$ but
still in the same scale window. Then   $\del c_{high}-\alpha c'_{high}>0$, whence obviously, $\del c_{high}-\alpha c'_{high}
+ \sum_{\gamma \ {\mathrm{low}}} w(\alpha)_{\gamma}>0$. Thus 
$\lambda^*>\alpha$.

\noindent Let now $\alpha\gg \alpha_{thr}$ be larger than $\alpha$ but
still in the same scale window. Since $\del c_{high}>0$ and we are not in the transition regime, 
the solution of the rectified equation $\del c_{high}-\alpha'_{thr} c'_{high} + \sum_{\gamma \ {\mathrm{low}}}
w(0)_{\gamma}=0$ is $\alpha'_{thr}\approx \alpha_{thr}$, hence $\alpha\gg \alpha'_{thr}$. But then,
 $\del c_{high}-\alpha c'_{high} + \sum_{\gamma \ {\mathrm{low}}}
w(\alpha)_{\gamma}< \del c_{high}-\alpha_{thr} c'_{high} + \sum_{\gamma \ {\mathrm{low}}}
w(0)_{\gamma}=0 $. 
Thus $\lambda^*<\alpha$.


\subsection{Threshold equation in terms of integer scales} \label{subsection:threshold-integer}


{\em We assume in the rest of the section that there is no degradation, so that all deficiency
rates $\kappa_{\sigma}$ are $\ge 0$.}  
 We shall now show that the threshold equation can be solved in terms of the integer
{\em scales} $n_{\sigma}, n_{\sigma\to\sigma'}$ and  the 
deficiency scales
\BEQ d_{\sigma}\sim \lfloor \log_{b} (\kappa_{\sigma})\rfloor \in \Z\cup\{-\infty\}
\EEQ
 ($d_{\sigma}=-\infty$ if $\kappa_{\sigma}=0$). Namely,
if $\gamma$ is an excursion of length $\ell$,
\BEQ n_{\gamma}= \min(n_{\sigma_i}; i=1,\ldots,\ell); 
\qquad N_{\gamma} \sim \sum_{i=1}^{\ell} n_{\sigma_i}; \EEQ
\BEQ -\log_{b}\eps_{\gamma}\sim  \min\Big( 
 n_{\sigma_i}-d_{\sigma_i}, \ 1\le i\le \ell-1\Big)  \label{eq:log-eps(gamma)}
 \EEQ
 
Then, $-\log_{b} c'_{high}$ may be replaced by 
\BEQ -n'_{high} := \min\Big(D_{\gamma}+n_{\gamma}\ |\ \gamma:\sigma\to\sigma
 \ {\mathrm{high}} \Big), \label{eq:n'high}
 \EEQ
  where
\BEQ D_{\gamma}= \sum_{i=1}^{\ell-1} (n_{\sigma_i} - 
n_{\sigma_i\to\sigma_{i+1}}) \sim -\lfloor \log_{b} 
\tilde{w}(0)_{\gamma}\rfloor \ge 0  \label{eq:D(gamma)}
\EEQ 
is the {\bf depth} of the excursion $\gamma$.

\Medskip Let us finally deal with
 $\del c_{high}$, which can be rewritten using 
 the normalization condition as  $-\del c_{high,1} + \del c_{high,2}$, with
\BEQ \del c_{high,1} = \sum_{\gamma:\sigma^*\to\sigma^* \ {\mathrm{low}}} (\tilde{w}(0))_{\gamma} \sim b^{-q_1}, \qquad q_1:= \min(D_{\gamma} \ |\ 
\gamma:\sigma\to\sigma \ {\mathrm{low}}) 
\EEQ
and (using (\ref{eq:log-eps(gamma)}) and (\ref{eq:D(gamma)}))
\BEQ \del c_{high,2} = \sum_{\gamma:\sigma^*\to\sigma^* \ {\mathrm{high}}}
 \eps_{\gamma} (\tilde{w}(0))_{\gamma} \sim  
 b^{-q_2}, \EEQ
 \BEQ  \qquad q_2:=\min\Big( D_{\gamma} + 
 \min(n_{\sigma_i}-p_{\sigma_i},\ 1\le i\le\ell-1)  \ |\   \gamma:\sigma\to\sigma \ {\mathrm{high}}\Big)  \label{eq:q2}
 \EEQ
Thus (leaving out equality cases), we are in Regime (1) 
if  $q_2<q_1$, implying $\del c_{high}\sim b^{-q_2}$.

\Medskip Choose then a scale $n_{\alpha} = \log_{b}\alpha$. 
"High" and "low" are notions which depend on a given
integer interval $I_{\alpha}$ for $n_{\alpha}$.
  Since only Regime (1), case (a) is relevant, we get for $\alpha= \alpha_{thr}$ the 
  final scaling form for the treshold equation:
\BEQ {\mathbf{(Threshold\ equation,\ scaling\ form)}} \qquad  n_{\alpha}\sim
-q_2 - n'_{high} .  \label{eq:thr-scaling-form}
\EEQ
 It should first be checked whether
 this value is contained in $I_{\alpha}$. If so, we must also
 check the conditions     $n_{k-low}-kn_{\alpha} <
 -q_2$ for all $k\ge 1$.


\subsection{Determination of the Lyapunov eigenvector} \label{subsection:Lyapunov-eigenvector}


We estimate here the Lyapunov weights $\pi^*_{\sigma} = (|A_{\sigma,\sigma}|+\lambda^*) v^*_{\sigma} $, 
see (\ref{eq:pi*}) for irreducible graphs, and then generalize to the
case of non-irreducible ones.  Recall (see \S \ref{subsection:resolvent-formula}) that     
$f(\alpha)_{\sigma^*\to\sigma}$ is the total
weight of
excursions from $\sigma^*$ to $\sigma$; we want to prove that 
$\pi^*_{\sigma}\simeq f(\alpha)_{\sigma^*\to\sigma}$ for any $\alpha\simeq \lambda^*
\simeq \alpha_{thr}$, where $\alpha_{thr}$ is the estimate of $\lambda^*$ obtained  indifferently from (\ref{eq:th:thr})
or (\ref{eq:thr-scaling-form}). Though $f(\alpha)_{\sigma^*\to\sigma}$ cannot
be computed exactly, it can be estimated using the leading
simple excursions $\gamma^*_{\sigma^*},\gamma^*_{\sigma}$ defined in  Theorem \ref{th:sigma*-dominant} 
(i), (ii).  Namely,
  $f(\alpha)_{\sigma^*\to\sigma^*}\sim 1$.  
If $\sigma\not=\sigma^*$, then 
$f(\alpha)_{\sigma^*\to\sigma}$  can be estimated by the weight of a leading  simple
excursion $\gamma^*_{\sigma}$ from $\sigma^*$ to $\sigma$. We formulate
our result in the following way:

\begin{Lemma}[Lyapunov weights for irreducible graphs] \label{lem:Lya-weight}
Normalize $\pi^*$ so that it is a probability measure. Then 
\BEQ \pi^*_{\sigma} \sim w(\alpha_{thr})_{\gamma^*_{\sigma}}, \qquad 
\sigma\in V  \label{eq:Lya-weight}
\EEQ
In particular, $\pi^*_{\sigma^*}\sim 1$.
\end{Lemma}

\noindent {\bf Proof.} Define $P := W(\lambda^*) \simeq W(\alpha_{thr})$ to be the discrete-time 
transition matrix with coefficients $P_{\sigma,\sigma'} = w(\lambda^*)_{\sigma\to\sigma'}$, and  
\BEQ G_N^{\sigma^*}(\sigma) \equiv G_N(\sigma^*\to\sigma) := \sum_{n=0}^{N-1} (P^n)_{\sigma^*,\sigma}\qquad (N\ge 1).  \label{eq:Green}
\EEQ
 The kernel 
$G_N=(G_N(\sigma^*\to\sigma))_{\sigma^*,\sigma}$ is a generalization of a truncated version of the  Green kernel, defined for non-defective, transient Markov chains as the
limit of the sum when $N\to\infty$ (the sum diverges for null-recurrent
chains); it is equal by definition to the average number of times 
a random trajectory with length $N$ started from $\sigma^*$ visits 
$\sigma$. The usual Green kernel is invariant by multiplication to
the right by the transition matrix; here, we get approximate invariance,
\BEQ (G_N^{\sigma^*} P)_{\sigma'} = \sum_{\sigma} G_N^{\sigma^*}(\sigma)
P_{\sigma,\sigma'} = \sum_{n=1}^N (P^n)_{\sigma^*,\sigma'} = 
G_N(\sigma^*\to\sigma') + (P^N)_{\sigma^*,\sigma'} - \del_{\sigma^*,\sigma'}. 
\EEQ  
Thus $\frac{1}{N} G^{\sigma^*}_N P \sim_{N\to\infty} \frac{1}{N} G^{\sigma^*}_N$, which (provided convergence holds) implies :
$\frac{1}{N} G_N(\sigma^*\to\sigma) \to c_0\pi^*_{\sigma}$ for some constant
$c_0> 0$.  Thus we need only prove that there exist two
 constants $c_1,c_2>0$ such that,  for $N$ large enough,
$c_1  w^*_{\gamma^*_{\sigma}} < \frac{1}{N}
 G_N(\sigma^*\to\sigma) < c_2  w^*_{\gamma^*_{\sigma}} $.

\Bigskip We discuss first a lower bound. Consider all paths $\gamma$ of length
$\le N$ such that $\gamma=\gamma_1\cup\gamma_2$ (concatenation
of paths), where $\gamma_2 = \gamma^*_{\sigma}$ (see Theorem \ref{th:sigma*-dominant} (ii)) is a leading
simple excursion from $\sigma^*$ to $\sigma$, 
and $\gamma_1$ is any path of length $N_1\le N-\ell(\gamma^*_{\sigma})$ from 
$\sigma^*$ to $\sigma^*$. By definition, $G_N(\sigma^*\to\sigma)$ is bounded from below by  the total weight of such paths. Then (summing over $N_1$), we
can use Theorem \ref{th:sigma*-dominant} (iv), and conclude that   $G_N(\sigma^*\to \sigma)> c_1 N w^*_{\gamma^*_{\sigma}}$ for some constant $c_1$ of order 1. 

\Medskip The upper bound is very simple: paths giving a further contribution
to $G_N(\sigma^*\to\sigma)$ are continuations of paths $\gamma$ 
of the form discussed in the lower bound that loop around $\sigma$, 
in other words, repeated excursions from $\sigma$ to $\sigma$. 
But Theorem \ref{th:sigma*-dominant} (iv) implies that the sum over 
such contributions converges to a finite value of small order $\prec 1$.  \hfill\eop

\Bigskip {\bf Extension to non-irreducible graphs.} Let 
${\cal G}^*$ be the SCC (strongly connected component) containing $\sigma^*$,  
and ${\cal G}'$ be the union of all other SCCs; thus 
$A$  is in the block form

\medskip
\begin{center}
\begin{tikzpicture}
\draw(0,0.8)--(-0.3,0.8)--(-0.3,-0.8)--(0,-0.8);
\draw(-2,0) node {$A=$};
\draw(0,1.2) node {${\cal G}^*$};
\draw(1,1.2) node {${\cal G}'$};
\draw(-0.7,0.5) node {${\cal G}^*$};
\draw(-0.7,-0.5) node {${\cal G}'$};
\draw(0,0.5) node {$A^*$};
\draw(1,0.5) node {$0$};
\draw(0,-0.5) node {$K$};
\draw(1,-0.5) node {$A'$};
\draw(1,0.8)--(1.25,0.8)--(1.25,-0.8)--(1,-0.8);

\end{tikzpicture}
\end{center}

\Medskip {\em We assume that the Lyapunov exponent $\lambda^*_{A'}$ of $A'$ is 
$<\lambda^*$.} Then $\lambda^* = \lambda^*(A^*)$. Start from
a (unique up to normalization) distribution vector $v^*_{{\cal G}^*}$ on ${\cal G}^*$  such that $(A^*-\lambda^*\Id) v^*_{{\cal G}^*} = 0$, 
then there exists a unique extension of $v^*_{{\cal G}^*}$ into 
a  quasi-stationary distribution vector $v^*$ of the block form
$v^* = \left[\begin{array}{c} v^*_{{\cal G}^*} \\ v^*_{{\cal G}'}
\end{array}\right]$; furthermore, 
$(A'-\lambda^*\Id) v^*_{{\cal G}'} = -Kv^*_{{\cal G}^*},$
 solved as 
\BEQ v^*_{{\cal G}'} = \int_0^{+\infty} dt\, e^{t(A'-\lambda^*)} 
\ Kv^*_{{\cal G}^*}
\EEQ
The integral converges because $\lambda^*_{A'}<\lambda^*$. 
One recognizes in the last expressions the resolvent of $A'$; the
solution (see (\ref{eq:traj})) may be expressed as a sum over paths,
\BEQ v^*_{\sigma'} = \sum_{\sigma''\in {\cal G}', \tau^* \in 
{\cal G}^*}  \Big(  \sum_{\gamma':\sigma''\to 
\sigma'} \frac{w(\lambda^*)_{\gamma'}}{|A_{\sigma',\sigma'}|+\lambda^*}  \Big) \  K_{\sigma'',\tau^*} \  v^*_{\tau^*}
\EEQ
which we can also rewritten as 
\BEQ \pi^*_{\sigma'}
= \sum_{\tau^* \in 
{\cal G}^*}   \pi^*_{\tau^*}\ \Big( \sum_{\gamma: \tau^*\to {\cal G}'\to \sigma'}  w(\lambda^*)_{\gamma}   \Big)
 \label{eq:Lya-weight-bis}
\EEQ
where
\BEQ \pi^*_{\sigma'} = (k_{\sigma'}+\lambda^*)v^*_{\sigma'}, \qquad
\pi^*_{\tau^*} = (k_{\tau^*}+\lambda^*)v^*_{\tau^*}
\EEQ
In the first expression, $\gamma'$ ranges over the set of paths
in ${\cal G}'$ from $\sigma''$\ to $\sigma'$. The product
$\frac{  w(\lambda)^*_{\gamma'} K_{\sigma'',\tau^*}}{k_{\tau^*}+\lambda^*} =
\frac{ k_{\tau^*\to \sigma''}\  w(\lambda)^*_{\gamma'}}{k_{\tau^*}+\lambda^*}$ may be rewritten as $w(\lambda^*)_{\gamma}$,
where $\gamma$ is the concatenation of the edge $\tau^*\to \sigma''$
with $\gamma'$. This gives the second expression, where the sum 
ranges over paths $\gamma$ from $\tau^*$ to $\sigma'$
whose first step is an edge $\tau^*\to \sigma''$ from $\tau^*$ 
to ${\cal G}'$.  Eq. (\ref{eq:Lya-weight-bis}) makes sense even
if $\sigma'$ is an absorbing state, provided $\lambda^*>0$ 
(namely, $k_{\sigma'}=0$ but the total transition rate $\lambda^*$ 
out of $\sigma'$ involves the degradation reaction $\sigma'\overset{\lambda^*}{\to} \emptyset$). 

\Medskip Note that (\ref{eq:Lya-weight-bis}) is very similar 
to (\ref{eq:Lya-weight}). By the same arguments as those used in 
the introductive Remarks to the section, or in the proof of Lemma
\ref{lem:Lya-weight}, it is enough to choose in the term between
parentheses in (\ref{eq:Lya-weight-bis}) a dominant, simple path,
say, $\gamma^{\tau^*}_{\sigma'}$. Then $\pi^*_{\sigma'} \sim 
\max\Big(\pi^*_{\tau^*} \,  w(\lambda^*)_{\gamma^{\tau^*}_{\sigma'}}, \ 
\tau^*\in {\cal G}^*\Big)$ can be estimated by finding  a 
starting vertex $\tau^*$ in the SCC ${\cal G}^*$ maximizing the
product of the weight $(\pi^*_{{\cal G}^*})_{\tau^*}$ (found using Lemma 
\ref{lem:Lya-weight}) by the weight of the corresponding dominant
path from $\tau^*$ to $\sigma'$.


\subsection{Directed acyclic graph (DAG) structure}


We now show that one may estimate of Lyapunov weights of a $\sigma^*$-dominant graph  $G$
following  Lemma \ref{lem:Lya-weight}  by multiplying the weights along the edges 
of a DAG (directed acyclic graph) extracted from $G$. Recall that a DAG is a directed graph 
$\T=(V,{\cal E})$
containing no cycle. Writing $x<x'$ if there exists a path from $x$ to $x'$ in $\T$ defines
a partial ordering on $\T$; minimal elements for this ordering are called roots.   An example of set endowed with a natural DAG structure is the set of 
communication classes of a Markov chain; minimal classes (in the usual definition) are also minimal
for the above ordering. 

\Medskip Let $\{\gamma^*_{\sigma,1},\ldots,\gamma^*_{\sigma,p_{\sigma}}\}$
be the set of all leading simple excursions from $\sigma^*$ to $\sigma\in V$; ${\cal E}$ the
(non-necessarily disjoint) union
of the edges of these for every $\sigma\in V$; and $\T$ be the 
graph with vertex set $V$ and edge set ${\cal E}$.   Recall that the integer $D_{\gamma}$, defined in (\ref{eq:D(gamma)}), 
is essentially minus the log-weight of the path $\gamma$.

\begin{Lemma} \label{lem:DAG} 
 $\T$ is a DAG rooted in $\sigma^*$. Furthermore, if $\gamma,\gamma'$ are two paths from 
 $x$ to $y$, then $D_{\gamma} = D_{\gamma'}$. 
\end{Lemma} 

\noindent {\bf Proof.}
We note the following properties.  We say here that $\gamma$ is of 
higher weight than $\gamma'$ if $D_{\gamma}<D_{\gamma'}$. 

\Medskip (1) If $x\in \gamma^*_{\sigma_i}\cap \gamma^*_{\sigma',i'}$ with $(\sigma,i)\not=
(\sigma',i')$, then $D_{\sigma^* \overset{\gamma^*_{\sigma,i}}{\to} x} =  D_{\sigma^*\overset{\gamma^*_{\sigma',i'}}{\to} x}$. Namely, assume by absurd that 
$D_{\sigma^* \overset{\gamma^*_{\sigma,i}}{\to} x} > D_{\sigma^*\overset{\gamma^*_{\sigma',i'}}{\to} x}$. Then the mixed path $\sigma^*\overset{\gamma^*_{\sigma',i'}}{\to} x \overset{\gamma^*_{\sigma,i}}{\to} \sigma$ is of higher  weight than $\gamma^*_{\sigma,i}$, which
is contradictory with the fact that $\gamma^*_{\sigma,i}$ is leading. Thus one may define
unambiguously  $D_x$ as the common value of all $D_{\sigma^*\overset{\gamma^*_{\sigma,i}}{\to} x}$. 

\Medskip (2) Generalizing, let  $\gamma$ be a path from 
 $\sigma^*$ to $x$. It is obtained by concatenating pieces of leading simple excursions 
 $\gamma_{\sigma_1,i_1}\big|_{\sigma^*\to x_1},\ldots, \gamma_{\sigma_q,i_q}\big|_{x_{q-1}\to x_q}$, with $x_q=x$. By (1), $D_{x_1}= D_{\sigma^*\overset{\gamma^*_{\sigma_1,i_1}}{\to} x_1} =
 D_{\sigma^*\overset{\gamma^*_{\sigma_2,i_2}}{\to} x_1}$. Thus one may replace 
 $\gamma_{\sigma_1,i_1}\big|_{\sigma^*\to x_1} \uplus \gamma_{\sigma_2,i_2}\big|_{x_1\to x_2}$ by $\gamma_{\sigma_2,i_2}\big|_{\sigma^*\to x_2}$ without changing the depth. Continuing inductively
 along the path, we get $D_{\gamma} = D_{\sigma^*\overset{\gamma_{\sigma_q,i_q}}{\to} x} = D_x$. 
 
 \Medskip (3) Generalizing further, let $\gamma$ be a path from $x$ to $y$. As in (2), 
 it is obtained by concatenating pieces of leading simple excursions 
 $\gamma_{\sigma_1,i_1}\big|_{x_0\to x_1},\ldots, \gamma_{\sigma_q,i_q}\big|_{x_{q-1}\to x_q}$, with $x_0=x, x_q=y$. Concatenating it with $\gamma_0:=\gamma_{\sigma_1,i_1}\big|_{\sigma^*\to x_1}$, we 
 get: $D_{\gamma_0\uplus \gamma} = D_x$. Hence $D_{\gamma}= D_x - D_{\gamma_0} = D_x-D_y$ is 
 independent of the choice of path from $x$ to $y$.

\Medskip (4) Assume (by absurd) that $\T$ contains a cycle $x_1\to x_2\to\cdots x_p\to x_1$. 
Remove from the list $x_1,\ldots,x_p$ all "intermediate" species $x_i$ such that  
$x_{i-1}\to x_i$ and $x_i\to x_{i+1}$ are edges along the same leading simple excursion, we 
get (after renumbering) a sequence of paths $x_1\overset{\gamma_1}{\to} x_2,\cdots, 
x_{q-1}\overset{\gamma_{q-1}}{\to} x_q, x_q \overset{\gamma_q}{\to} x_1$, with 
$\gamma_j \subset \gamma^*_{\sigma_j,q_j}$. Since the depth is strictly increasing along a path,
$D(x_1)<D(x_2)<\ldots<D(x_1)$: contradiction.  \hfill\eop


\subsection{Examples}  \label{subsection:multiscale-examples}


Let us illustrate the computations of this section for Example 1 and Example 2. Apart from presenting the
results of Sections \ref{section:example1} and \ref{section:example2} as part of a systematic algorithm,
the interest is to derive the estimates for the Lyapunov eigenvectors. The threshold parameter 
$\alpha_{thr}$ is estimated from the scaling form (\ref{eq:thr-scaling-form}), and then, the Lyapunov
eigenvector $v^*$ from  (\ref{eq:Lya-weight}).


\subsubsection{Example 1}


We present here in full detail the computations in the case
when $n_1>n_2>n_3$, with $n_1=\lfloor \log_{b} k_{on}\rfloor$,
$n_1=\lfloor \log_{b} k_{off}\rfloor$, 
$n_1=\lfloor \log_{b} \nu_+\rfloor$  (see case 2 in  Fig. \ref{fig:example1}), and $n_1=\lfloor \log_{b} k_{off}\rfloor$,
$n_1=\lfloor \log_{b} k_{on}\rfloor$, 
$n_1=\lfloor \log_{b} \nu_+\rfloor$  (case 2'). Since
$\lambda^*\le \nu_+$, $n_{\alpha}\le n_3$ in both cases, so
that both excursions from $\sigma^*=1$ to 1,
\BEQ \gamma_1: 1 \overset{k_{off}}{\to} 0 \to 1, \qquad \gamma_2: 1 \overset{\nu_+}{\twoheadrightarrow 0} \to 1 \EEQ
are high, so we are in regime (1), case (a). Furthermore, the scale of $\gamma_1,\gamma_2$   is $n_{\gamma_1}=n_{\gamma_2}=n_2$, and the depth of $\gamma_1$ is
$D_{\gamma_1}=0$. Also, $\eps_1=\frac{\nu_+}{\nu_++k_{off}}\sim \frac{\nu_+}{k_{off}}$. 

\begin{itemize}
\item[(1)] (see Fig. 2) Here, $-\log_{b}(\eps_1)\sim D_{\gamma_2}\sim  n_2-n_3$. Thus
\BEQ -n'_{high}= \min(D_{\gamma_1}+n_{\gamma_1}, D_{\gamma_2}
+n_{\gamma_2})= D_{\gamma_1}+n_{\gamma_1} =  n_2; 
\label{eq:5.25}
 \EEQ
\BEQ q_2 = \min(D_{\gamma_1}+(n_2-n_3),D_{\gamma_2}+(n_2-n_3) )
 = D_{\gamma_1}+(n_2-n_3) = n_2-n_3
 \label{eq:5.26}
\EEQ
hence  $n_{\alpha}\sim -q_2-n'_{high} = n_3\sim \log_{b}(\nu_+)$.

Note that the unique path minimizing (\ref{eq:n'high}), see (\ref{eq:5.25}), is 
$\gamma_1$,  and similarly for (\ref{eq:q2}), see
(\ref{eq:5.26}). Thus 
\BEQ v^* \sim a^{-1} \sim \left(\begin{array}{c} k_{on}^{-1} \\
k_{off}^{-1} \end{array}\right)  \label{eq:5.29}
\EEQ

\item[(2)] (see Fig. 2') Here, $-\log_{b}(\eps_1)\sim D_{\gamma_2} \sim n_1-n_3$. Thus $-n'_{high}=n_2$, 
$q_2=n_1-n_3$ instead of $n_2-n_3$. Hence
 $n_{\alpha}\sim -q_2-n'_{high} = n_2-n_1+n_3 \sim \log_{b}
 (\frac{k_{on}}{k_{off}}\nu_+)$. As in (1), the unique 
 minimizing path is $\gamma_1$. Hence (\ref{eq:5.29}) also
 holds in this case.
\end{itemize}

\Medskip Let us now present briefly the results in the other 
cases. {\em In case 1 or 1',} ($\nu_+\gg k_{on},k_{off}$), 
$\lambda^*\approx k_{on}$,  $-\log_{b}(\eps_1)\sim 0$, 
$n_{\gamma_1}=n_{\gamma_2} = n_3$ (case 1) or $n_2$ (case 1'). Then $D_{\gamma_1}\sim n_1-n_2$ (case 1) or $n_1-n_3$ (case 1'),
 $D_{\gamma_2}\sim 0$, $q_2\sim D_{\gamma_2}\sim 0$, and 
 $-n'_{high} = D_{\gamma_2}+n_{\gamma_2} = n_3$ (case 1), 
 $n_2$ (case 1'). The unique minimizing path is $\gamma_2$ in
 both cases. Therefore $v^* \sim a^{-1}\sim \left(\begin{array}{c}
 k_{on}^{-1} \\ \nu_+^{-1} \end{array}\right)$. 
 {\em In case 3} ($k_{on}\gg \nu_+\gg k_{off}$), $\lambda^*\approx
 \nu_+$, $n_{\gamma_1}=n_{\gamma_2}=n_2$, $D_{\gamma_1}\sim 
 n_2-n_3, D_{\gamma_2}\sim 0$, and $-\log_{b}(\eps_1)\sim 0$.
 Therefore the unique minimizing path is $\gamma_2$,  
  $q_2\sim 0, -n'_{high} \sim n_2$, and $v^* \sim a^{-1}\sim \left(\begin{array}{c}
 k_{on}^{-1} \\ \nu_+^{-1} \end{array}\right)$.  Finally,
 {\em in case 3'} ($k_{off}\gg \nu_+\gg k_{on}$), $\lambda^* \approx \frac{k_{on}}{k_{off}}\nu_+$, $-\log_{b}(\eps_1)\sim
 n_1-n_2$, $n_{\gamma_1}=n_{\gamma_2}=n_3$, $D_{\gamma_1}\sim 0$,
 $D_{\gamma_2}\sim n_1-n_2$, $q_2\sim n_1-n_2$, $-n'_{high} \sim
 n_3$. The unique minimizing path is $\gamma_1$ , therefore $v^*\sim a^{-1}\sim \left(\begin{array}{c}
 k_{on}^{-1} \\ k_{off}^{-1} \end{array}\right)$.


\subsubsection{Example 2}


We take $\sigma^*=1$, so that
\BEQ {\mathrm{(vertex\ scales)}} \qquad
n_{\gamma_1} = n_3, \ n_{\gamma_2}=n_4, \ p_1=n_2 \EEQ
\BEQ  {\mathrm{(path\ depths)}} \qquad  D_{\gamma_1}=0, \ D_{\gamma_2} = n_3 -n_{0\to 4} \EEQ

\Bigskip {\bf Assume first (A)} : $n_4\prec n_{\alpha}\prec n_3$, 
namely, $k_4\ll \alpha \ll k_{on}$.   Then $\gamma_1 : 1\to 0\to 1$ is high,
$\gamma_2: 1\to 0\to 4\to 1$ is low, $\gamma_2^{low}=(4)$, 
$N_{\gamma_2^{low}}\sim n_4$. Further,
\BEQ n'_{high} = -n_{\gamma_1}=-n_3; \EEQ
\BEQ n_{1-low} = -D_{\gamma_2} + N_{\gamma_2^{low}} \sim
-n_3+n_{0\to 4} + n_4 ;\EEQ
and
\BEQ \begin{cases} q_1 = D_{\gamma_2} = n_3-n_{0\to 4}
\sim \log_{b}(k_{on}/k'_4) \\
q_2 = D_{\gamma_1} + n_1-d_1 = n_1-n_2 \sim \log_{b}(k_{off}/\nu_+) \end{cases}
\EEQ

Regime (1) is defined by $q_2\prec q_1$, which amounts to 
$R:=(\frac{k_{on}}{k_{off}}\nu_+)/k'_4\gg 1$.  Then
in case (a), we get 
\BEQ n_{\alpha} \sim -q_2-n'_{high} \sim -n_1+n_2+n_3 \sim \log_{b}(\frac{k_{on}}{k_{off}}\nu_+) \EEQ
and the condition $n_{1-low}-n_{\alpha}\prec -q_2$ is verified
since $n_{1-low}-n_{\alpha}+q_2 = 2(n_1-n_2-n_3) + n_{0\to 4}
+ n_4 \sim -2\log_{b}R - \log(k'_4/k_4)<0$. 

\Medskip Regime (2) is defined by $q_1\prec q_2$. In case $(a_1)$, 
we get $n_{\alpha} \sim n_{1-low} + q_1= n_4$, which is rejected
since $n_4\prec n_{\alpha}$ by hypothesis.

\Medskip Both $\gamma_1,\gamma_2$ play a role in the determination of main indices $n'_{high},n_{1-low},q_1,q_2$. The 
weights $w(\alpha)_{\gamma_1}, w(\alpha)_{\gamma_2}$ are
$\sim 1$ and $\sim \frac{k'_4}{k_{on}} \, \times\, 
\frac{k_4}{k_4+\alpha} \sim \frac{k'_4}{k_{on}} \frac{k_4}{\alpha}
\sim \frac{k'_4 k_4 k_{off}}{k_{on}^2 \nu_+}$; taking logs
and changing signs, we get depth indices $D(\alpha)_{\gamma_1}\sim 0, D(\alpha)_{\gamma_2} \sim -n_1+ n_2+2n_3 - n_{0\to 4}-n_4>0$.

\Medskip We now estimate the Lyapunov weights by using Lemma
\ref{lem:Lya-weight}. The dominant simple excursion from $1$ to $0$,
resp. $4$, may be chosen as $\gamma^*_1 : 1\to 0$, resp. $\gamma^*_4: g
1\to 0 \to 4$. Since $0,1$ are high, we readily find $w(\alpha)_{\gamma^*_1} \sim 1, w(\alpha)_{\gamma^*_4} \sim 
\frac{k'_4}{k_{on}}$; thus  
\BEQ \pi^*_0 \sim \pi^*_1 \sim 1, \qquad \pi^*_4 \sim \frac{k'_4}{k_{on}}
\EEQ
from which, dividing by $k_{\sigma}+\alpha\sim \begin{cases} 
k_{on} \\ k_{off} \\ \alpha \end{cases}$, 
 $ v^* \propto\left(\begin{array}{c} 1 \\ \frac{k_{on}}{k_{off}} \\
\frac{k'_4 k_{off}}{k_{on} \nu_+} \end{array}\right)
$
or equivalently,
\BEQ -\log_{b} v^* \sim  \left(\begin{array}{c} 0 \\ 
n_1-n_3 \\ (n_2+n_3) - (n_1+n_{0\to 4}) \end{array}\right). \EEQ

 
\Bigskip {\bf Assume now  (B)}:  $n_{\alpha}\prec n_4$, namely,
$\alpha\ll k_4$. Then both $\gamma_1,\gamma_2$ are high, no
we are automatically in Regime (1), case (a), and $q_1$ needs
not be introduced. Then
\BEQ n'_{high} = -\min(D_{\gamma_1}+n_{\gamma_1},D_{\gamma_2}+
n_{\gamma_2}) = -(D_{\gamma_2}+n_{\gamma_2}) = n_{0\to 4}-n_3-n_4 
\sim -\log_{b}(k_{on}/k'_4k_4) 
\EEQ
and
\BEQ q_2= \min(D_{\gamma_1}+n_1-d_1,D_{\gamma_2}+n_1-d_1)
=n_1-n_2 \sim \log_{b}(k_{off}/\nu_+)
\EEQ
so then, 
\BEQ n_{\alpha} \sim -q_2-n'_{high} \sim \log_{b}(Rk_4)
\EEQ

Here also, both paths $\gamma_1,\gamma_2$ play a role in the
determination of the indices, but now (since $n_{\alpha}$ is
the smallest scale), $D(\alpha)_{\gamma_1}\sim D_{\gamma_1}\sim 0,
D(\alpha)_{\gamma_2}\sim D_{\gamma_2}\sim n_3-n_{0\to 4}$. 

\Medskip We again estimate the Lyapunov weights using $\gamma^*_1,
\gamma^*_4$ and 
 get similarly (with $k_4+\alpha\sim k_4$ this time)
$v^* \sim 
\left(\begin{array}{c} k_{on}^{-1} \\ k_{off}^{-1} \\
\frac{k'_4}{k_4 k_{on}} \end{array}\right)$. 
The largest component is $v^*_4 \sim \frac{k'_4}{k_4 k_{on}}$. 
Dividing by $v^*_4$ yields  the normalized formula,
\BEQ v^* \sim \left(\begin{array}{c} \frac{k_4}{k'_4} \\
\frac{k_{on} k_4}{k_{off} k'_4} \\ 1 \end{array}\right)
\EEQ
or equivalently,
\BEQ -\log_{b} v^* \sim \left(\begin{array}{c} n_{0\to 4}-n_4 \\
(n_1+n_{0\to 4})-(n_3+n_4) \\ 0 \end{array}\right)
\EEQ

\Bigskip {\bf Lyapunov eigenvector.}  Normalize
$v^*$ so that its top coefficient is $\sim 1$, then (from Lemma \ref{lem:Lya-weight})
\BEQ  v^* \sim \begin{cases} v^{{\mathbf{A}}} \equiv  \left(\begin{array}{c} 1 \\ \frac{k_{on}}{k_{off}} \\
\frac{k'_4 k_{off}}{k_{on} \nu_+} \end{array}\right), \qquad 
R\gg 1 \ {\mathbf{(A)}} \\ 
v^{\mathbf{B}} \equiv \left(\begin{array}{c} \frac{k_4}{k'_4} \\
\frac{k_{on} k_4}{k_{off} k'_4} \\ 1 \end{array}\right), \qquad 
R\ll 1   \ {\mathbf{(B)}} \end{cases}
\label{eq:vAvB}
\EEQ

\begin{figure}[!t]
\centering
\includegraphics[scale=0.6]{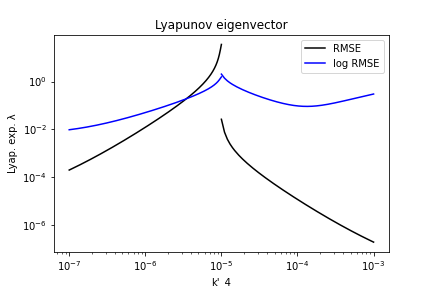}
\caption{RMSE errors on Lyapunov eigenvector as a function of kinetic parameter
$k'_4$}
\label{Fig:2}
\end{figure}

We computed the root mean-squared errors (RMSE), in black on Figure \ref{Fig:2},  $ \Big(\sum_{i=0,1,4}
(v^*_i-v_i^{\mathbf{A}})^2\Big)^{1/2}$  ($R\gg 1$),   $ \Big(\sum_{i=0,1,4}
(v^*_i-v_i^{\mathbf{B}})^2\Big)^{1/2}$  ($R\ll 1$), and also
the logarithmic  RMSE (in blue),  $\Big( \sum_{i=0,1,4}
(\log_{b}(v^*_i/v_i^{\mathbf{A}}))^2\Big)^{1/2}$  ($R\gg 1$),   $
\Big( \sum_{i=0,1,4}
(\log_{b}(v^*_i/v_i^{\mathbf{B}}))^2\Big)^{1/2}$  ($R\ll 1$).  Figure \ref{Fig:2}
reproduces the results; both RMSE are small, except around
transition scales. The accurateness of (\ref{eq:vAvB}) is corroborated by the smallness of the logarithmic RMSE (the RMSE fails to give any quantitative information about the smallest component of $v^*$). 

\section{Example 3. Formose}  \label{section:formose}


The formose network consists of non-branching oses (sugars) $X_i$ with 
$i\ge 2$ carbon atoms (formaldehyde $X_1$ is assumed to
be chemostated), related by addition reactions $X_j+X_{i-j}\to X_i$ called aldolizations,
and reverse fragmentation reactions, $X_i\to X_j+X_{i-j}$, called retroaldolizations.  Species $X_2,X_3,X_4,\dots,X_n$ are labeled by the number of 
carbon atoms (keto-enol tautomers are considered as one species,
since tautomerization is a fast reaction, and branching oses are left out).  Main
reactions are of 3 types,
\begin{itemize}
\item[(i)] aldolizations with formaldehyde, $X_i \to X_{i+1}, \qquad i\ge 2$;
\item[(ii)]  their reverse counterparts, $X_{i+1} \to X_i, \qquad i\ge 2$;
\item[(iii)] general retroaldolizations, $X_i \to X_j + X_{i-j}, \qquad j,i-j\ge 2$. 
\end{itemize}

\Medskip It is reasonable to assume that kinetic constants
of reactions within a given type (i), (ii) or (iii) are 
of the same order. We restrict to oses with $\le n=5$ carbons, both because extension to larger $n$ is straightforward and
can be done by the reader, and because larger oses tend to
 cyclize, requiring a more detailed analysis for relevance. Within our framework, this means that there
are three relevant scales (as in Example 1, with an obvious
identification by letting 
kinetic rates $\sim k_{on}$ for (i), $\sim k_{off}$ for (ii),
$\sim \nu_+$ for (iii)),  and six possible
multi-scale splittings, in particular (compare with subcases 2,2'
in Fig. \ref{fig:example1})

\bigskip
\begin{center}
\begin{tikzpicture}[scale = 0.8]

\draw(2.5,1.5) node {$2$};

\draw(0,0.4) node {$2$};
\draw(1.25,0.4) node {$3$};
\draw(2.75,0.4) node {$4$};
\draw[dotted](0,0)--(0,-3.05);
\draw[dotted](1.25,0)--(1.25,-2.95);
\draw[dotted](2.75,0)--(2.75,-2.5);
\draw[dotted](4,0)--(4,-3.05);

\draw(4,0.4) node {$5$};
\draw(0,0) node {\textbullet};
\draw[->, ultra thick](0,0)--(1,0); 
\draw(1.25,0) node {\textbullet};
\draw[->, ultra thick](1.5,0)--(2.5,0); 
 
 \draw(2.75,0) node {\textbullet};
\draw[->, ultra thick](3,0)--(4,0); 

\draw(6,0) node {$n_1$};

\draw[dashed](-1,-0.5)--(6,-0.5); 

\draw(4,-1) node {\textbullet};
\draw[->, ultra thick](4,-1)--(2.75,-1);
\draw[->](2.75,-1.5)--(1.25,-1.5);
\draw[->](1.25,-1)--(0,-1);
\draw(6,-1.25) node {$n_2$};

\draw[dashed](-1,-2)--(6,-2); 

\draw[->](2.75,-2.5)--(0.2,-2.5);
\draw[->](0.2,-2.5)--(0,-2.5);

\draw(4,-2.95)--(1.45,-2.95);
\draw[->](1.45,-2.95)--(1.25,-2.95);

\draw(4,-3.05)--(0.2,-3.05);
\draw[->](1.45,-3.05)--(0,-3.05);

\draw(6,-2.75) node {$n_3$};


\begin{scope}[shift={(10,0)}]
\draw(2.5,1.5) node {$2'$};

\draw(0,0.4) node {$2$};
\draw(1.25,0.4) node {$3$};
\draw(2.75,0.4) node {$4$};
\draw[dotted](0,0)--(0,-3.05);
\draw[dotted](1.25,0)--(1.25,-2.95);
\draw[dotted](2.75,0)--(2.75,-2.5);
\draw[dotted](4,0)--(4,-3.05);

\draw(4,0.4) node {$5$};
\draw(4,0) node {\textbullet};
\draw[<-, ultra thick](0,0)--(1,0); 
\draw(1.25,0) node {\textbullet};
\draw[<-, ultra thick](1.5,0)--(2.5,0); 
 
 \draw(2.75,0) node {\textbullet};
\draw[<-, ultra thick](3,0)--(4,0); 

\draw(6,0) node {$n_1$};

\draw[dashed](-1,-0.5)--(6,-0.5); 

\draw(0,-1) node {\textbullet};
\draw[<-](4,-1)--(2.75,-1);
\draw[<-](2.75,-1.5)--(1.25,-1.5);
\draw[<-, ultra thick](1.25,-1)--(0,-1);
\draw(6,-1.25) node {$n_2$};

\draw[dashed](-1,-2)--(6,-2); 

\draw[->](2.75,-2.5)--(0.2,-2.5);
\draw[->](0.2,-2.5)--(0,-2.5);

\draw(4,-2.95)--(1.45,-2.95);
\draw[->](1.45,-2.95)--(1.25,-2.95);

\draw(4,-3.05)--(0.2,-3.05);
\draw[->](1.45,-3.05)--(0,-3.05);

\draw(6,-2.75) node {$n_3$};

\end{scope}
\end{tikzpicture}
\end{center}

\Bigskip   For both splittings (Case 2, Case 2'), all excursions are high, since
$n_{\alpha}\le n_3<\min_{\sigma} n_{\sigma}$. Thus we
are in regime (1), case (a). We choose $\sigma^*=5$ in Case 2, and
$\sigma^*=2$ in Case 2'. Note that it is not possible to choose
e.g. $\sigma^*=2$ in Case 2, because the $\sigma^*$-dominant
Hypothesis would not be verified (the cycle $4\rightleftarrows 5$ 
is dominant in that case). We use bold arrows for dominant edges.

\Medskip 
{\bf Case 2.}  Vertex scales are $n_1$ for $\sigma=2,3,4$ and
$n_2$ for $\sigma=5$.  Main excursions are 
\BEQ \gamma_1: 5 \pmb{\to} 4 \pmb{\to} 5, \ \gamma_2: 5 \pmb{\to} 4\to 3\pmb{\to} 4 \pmb{\to} 5,
\  \gamma_3: 5 \pmb{\to} 4\to 3\to 2 \pmb{\to} 3 \pmb{\to} 4 \pmb{\to} 5 \EEQ

\begin{tikzpicture}
\draw(-1.75,0) node {$\gamma_4: 5$};
\draw[->](-1,0)--(0,0); 
\draw(1.3,0) node {$3\pmb{\to} 4\pmb{\to} 5,$};

\begin{scope}[shift={(6.6,0)}] 
\draw(-2,0) node {$\gamma_5 =5$};
\draw[->](-1.25,0)--(-0.25,0); 
\draw(1.3,0) node {$2\pmb{\to} 3\pmb{\to} 4 \pmb{\to} 5$};
\end{scope}

\draw(0.5,-1) node {$\gamma'_6: 5 \pmb{\to} 4$};
\draw[->] (1.5,-1)--(2.3,-1);
\draw[->](2.3,-1)--(2.5,-1);
\draw(4,-1) node {$2\pmb{\to} 3\pmb{\to} 4\pmb{\to} 5$};
\end{tikzpicture}

We have $-\log_{b}(\eps_4)\sim n_1-n_3 \succ -\log_{b}(\eps_5) \sim
n_2-n_3$, whence  $-\log_{b} \eps_{\gamma}  \sim n_2-n_3$ for 
all $\gamma=\gamma_1,\ldots,\gamma_6$. All excursion scales are equal: $n_{\gamma}=n_2$
for all $\gamma$. 
Excursion depths are
\BEQ D_{\gamma_1}=0,\ D_{\gamma_2}=n_1-n_2, \
 D_{\gamma_3}=2(n_1-n_2); \qquad 
 D_{\gamma_4} =  D_{\gamma_5} = n_1-n_3; \qquad
 D_{\gamma'_3} =  n_1-n_3
\EEQ

We have $-\log_{b}(\eps_4)\sim n_1-n_3 \succ -\log_{b}(\eps_5) \sim
n_2-n_3$, whence  $-\log_{b} \eps_{\gamma}  \sim n_1-n_3$ for 
$\gamma=\gamma_1,\ldots,\gamma_6$.  Then
\BEQ -n'_{high} = n_2 +  \min_{\gamma} D_{\gamma} = n_2 \EEQ
\BEQ q_2 = \min_{\gamma} (D_{\gamma} + (n_2-n_3)) = n_2-n_3 \EEQ 
Thus
\BEQ n_{\alpha} \sim -q_2-n'_{high} \sim n_3 \sim \log_{b} (\nu_+) 
\EEQ
Thus, the Lyapunov exponent is $\sim \nu_+$: retroaldolizations $X_4\to 2X_2, X_5\to X_2+X_3$
are the kinetically limiting steps.
 
\Medskip   Let us now discuss the Lyapunov eigenvector using 
Lemma \ref{lem:Lya-weight}. The dominant simple excursion to $4$ is  $\gamma^*_4 : 5\to 4$, with weight $w^*_{\gamma^*_4}\sim 1$.  Things
are not so clear for excursions to $3$ or $2$. Let  $\gamma^{*,1}_3 : 5 \to 4 \to 3,\ \gamma^{*,2}_3:
5\longrightarrow 3$, and $\gamma^{*,1}_2: 5\longrightarrow 2, \ \gamma^{*,2}_2: 
5\to 4\longrightarrow 2, \ \gamma^{*,3}_3:5\to 4\to 3\to 2$. The
respective  weights are  $w^*_{\gamma^{*,1}_3}\sim \frac{k_{off}}{k_{on}}, 
w^*_{\gamma^{*,2}_3} \sim \frac{\nu_+}{k_{off}}$, and
 $w^*_{\gamma^{*,1}_2}
\sim \frac{\nu_+}{k_{off}} \succ w^*_{\gamma^{*,2}_2}
\sim \frac{\nu_+}{k_{on}}, \ w^*_{\gamma^{*,3}_2}
\sim (\frac{k_{off}}{k_{on}})^2 $. 
Thus,  we get three cases:
\begin{enumerate}
\item $(\frac{\nu_+}{k_{off}}\ge \frac{k_{off}}{k_{on}})$ 
$\pi^* \sim \left(\begin{array}{c} \frac{\nu_+}{k_{off}} \\ \frac{\nu_+}{k_{off}} \\ 1
 \\  1 \end{array}\right)$;
\item $(\frac{k_{off}}{k_{on}})^2 \le \frac{\nu_+}{k_{off}} \le  \frac{k_{off}}{k_{on}})$ 
$\pi^* \sim \left(\begin{array}{c} \frac{\nu_+}{k_{off}} \\ \frac{k_{off}}{k_{on}} \\ 1
 \\  1 \end{array}\right)$;
\item $(\frac{\nu_+}{k_{off}}\le (\frac{k_{off}}{k_{on}})^2)$
$\pi^* \sim \left(\begin{array}{c}  (\frac{k_{off}}{k_{on}})^2  \\ \frac{k_{off}}{k_{on}} \\ 1
 \\  1 \end{array}\right)$.
\end{enumerate} 
 In the end, one must divide
 by $k_{\sigma}+\lambda^* \sim \begin{cases} k_{on} \\ k_{on} \\
 k_{on} \\ k_{off} \end{cases}$ to get $v^*$.

\Bigskip 
{\bf Case 2'.}  Vertex scales are $n_1$ for $\sigma=3,4,5$ and
$n_2$ for $\sigma=2$.  Main excursions are 
\BEQ \gamma_1: 2 \pmb{\to} 3\pmb{\to} 2, \ \gamma_2: 2\pmb{\to} 3\to 4\pmb{\to} 3\pmb{\to} 2,
\  \gamma_3: 2\pmb{\to} 3\to 4\to 5\pmb{\to} 4\pmb{\to} 3\pmb{\to} 2 \EEQ

\begin{tikzpicture}
\draw(0,0) node {$\gamma_4: 2 \pmb{\to} 3\to 4$};
\draw[->] (1.5,0)--(2.3,0);
\draw[->](2.3,0)--(2.5,0);
\draw(5.3,0) node {$2,\ \gamma'_4 : 2\pmb{\to} 3\to 4\to 5\pmb{\to} 4 $};
\begin{scope}[shift={(6.6,0)}]
\draw[->] (1.5,0)--(2.3,0);
\draw[->](2.3,0)--(2.5,0);
\end{scope}
\draw(9.3,0) node {$2$};

\draw(0,-1) node {$\gamma_5: 2\pmb{\to} 3\to 4\to 5$};
\draw (2,-1)--(2.8,-1);
\draw[->](2.8,-1)--(3,-1);
\draw(3.3,-1) node {$2,$};

\begin{scope}[shift={(6,1)}]
\draw(0,-2) node {$\gamma_6: 2\pmb{\to} 3\to 4\to 5$};
\draw (2,-2)--(2.8,-2);
\draw[->](2.8,-2)--(3,-2);
\draw(3.8,-2) node {$3\pmb{\to} 2$};
\end{scope}
\end{tikzpicture}

\Medskip We have $-\log_{b}(\eps_4) \sim
-\log_{b}(\eps_5) \sim n_1-n_3$, whence (excluding $\gamma_1$
since it does not go through $X_4$ or $X_5$),  $-\log_{b} (\eps_{\gamma}) \sim n_1-n_3$
for all $\gamma\not=\gamma_1$.  All excursion scales are
equal: $n_{\gamma}=n_2$ for all $\gamma$.
Then, excursion depths are
\BEQ D_{\gamma_1}=0,\ D_{\gamma_2}=n_1-n_2, \
 D_{\gamma_3}=2(n_1-n_2)
\EEQ
\BEQ D_{\gamma_4}= (n_1-n_2)+(n_1-n_3), \ 
D_{\gamma'_4} = 2(n_1-n_2)+  (n_1-n_3) \EEQ
\BEQ D_{\gamma_5} = D_{\gamma_6} = 2(n_1-n_2)+(n_1-n_3)
\EEQ

Then
\BEQ -n'_{high} = n_2 + \min_{\gamma} D_{\gamma} =n_2
\EEQ 

\BEQ q_2 = \min_{\gamma\not=\gamma_1} (D_{\gamma} + (n_1-n_3)) =
(n_1-n_2) + (n_1-n_3)
\EEQ

Thus
\BEQ n_{\alpha} \sim -q_2-n'_{high} \sim  
2n_2-2n_1 + n_3 \sim \log_{b} \Big( (\frac{k_{on}}{k_{off}})^2 \nu_+
\Big)
\EEQ  
as compared to $\log_{b} \Big( (\frac{k_{on}}{k_{off}}) \nu_+
\Big)$ for Example 1. The ratio $\frac{k_{on}}{k_{off}}$ is
squared roughly because it takes two steps to go from $X_2$ to $X_4$
(reactant of the 1-2 reaction $X_4\to 2X_2$), and each
would cost (at equilibrium, letting $\nu_+=0$)  $\sim 
\frac{k_{on}}{k_{off}}$. 

\Medskip Let us now discuss the Lyapunov eigenvector using 
Lemma \ref{lem:Lya-weight}. Dominant simple excursions are
$\gamma^*_3 : 2\to 3,  \ \gamma^*_4 : 2 \to 3 \to 4,\ \gamma^*_5:
2\to 3\to 4\to 5$, with weights $w^*_{\gamma^*_3}\sim 1, 
w^*_{\gamma^*_4} \sim \frac{k_{on}}{k_{off}}, \ w^*_{\gamma^*_5}
\sim (\frac{k_{on}}{k_{off}})^2$. Thus 
$\pi^* \sim \left(\begin{array}{c} 1 \\ 1 \\  \frac{k_{on}}{k_{off}}
 \\  (\frac{k_{on}}{k_{off}})^2 \end{array}\right)$. Dividing
 by $k_{\sigma}+\lambda^* \sim \begin{cases} k_{on} \\ k_{off} \\
 k_{off} \\ k_{off} \end{cases}$, we get 
\BEQ v^* \sim k_{on}^{-1} \ \left(\begin{array}{c} 1   \\  \frac{k_{on}}{k_{off}}
 \\  (\frac{k_{on}}{k_{off}})^2 \\   (\frac{k_{on}}{k_{off}})^3
  \end{array}\right)
\EEQ



\section{Appendix}  \label{section:appendix} 


\subsection{Proof of Theorem \ref{thm:homogeneous}}


We first prove bounds valid for an arbitrary chemical network, possibly with degradation, see 
{\bf A.}. Applying these bounds in {\bf B.} for degradationless networks yields Theorem 
\ref{thm:homogeneous}.   We use the notation
\BEQ D := \max_{\sigma\in V} \Big(\frac{\max(\kappa_{\sigma},0)}{|A_{\sigma,\sigma}|} \Big). 
\label{eq:D}
 \EEQ

\Medskip {\bf A.} {\em We start by proving a general double inequality for $\lambda^*$, see (\ref{eq:double-inequality-lambda*}). } Choose
some uniform degradation rate $\alpha\ge 0$. Let
$f(\alpha)_{\sigma\to \sigma'}$, resp. $\tilde{f}(\alpha)_{\sigma\to \sigma'}$ be the weight of excursions from $\sigma$ to $\sigma'$ for the discrete-time  generalized Markov chain associated to the generator $A$, resp. for the discrete-time
Markov chain associated to $\tilde{A}$. In particular, if it can be proven that $\max_{\sigma}(f(\alpha)_{\sigma\to
\sigma})>1$, then the network is autocatalytic, implying
the lower bound $\lambda^*>\alpha$.  Fix  a species index $\sigma^*\in V$.  A first-step analysis yields the  following linear system
for the vector $(\tilde{f}(\alpha)_{\sigma\to\sigma^*}))_{\sigma}$, 
\BEQ  \begin{cases} \tilde{f}(\alpha)_{\sigma\to \sigma^*} = \sum_{\sigma'\not=\sigma} \tilde{w}(\alpha)_{\sigma\to \sigma'}
\tilde{f}(\alpha)_{\sigma'\to\sigma^*} + \tilde{w}(\alpha)_{\sigma\to \sigma^*} \qquad \sigma\not=\sigma^* \\
\tilde{f}(\alpha)_{\sigma^*\to \sigma^*} =\sum_{\sigma'\not=\sigma^*} \tilde{w}(\alpha)_{\sigma^*\to \sigma'}
\tilde{f}(\alpha)_{\sigma'\to\sigma^*} 
\end{cases}
\EEQ
In particular,
\BEQ \sum_{\sigma'\not=\sigma^*} \tilde{w}(0)_{\sigma^*\to
\sigma'} \tilde{f}(0)_{\sigma'\to\sigma^*} = 
\tilde{f}(0)_{\sigma^*\to\sigma^*}=1 \label{eq:proba-cons}
\EEQ
by probability conservation.

\Medskip
The solution $f=f(w)$ of the general system 
$\begin{cases} f_{\sigma\to\sigma^*} = \sum_{\sigma'\not=\sigma}
w_{\sigma\to\sigma'} f_{\sigma'\to\sigma^*} + w_{\sigma\to
\sigma^*} \qquad \sigma\not=\sigma^* \\
f_{\sigma^*\to\sigma^*} = \sum_{\sigma'\not=\sigma^*}
w_{\sigma^*\to\sigma'} f_{\sigma'\to\sigma^*} \end{cases}$
with $w\ge 0$  may be obtained by iterating the system of equations\\
$\begin{cases} f^{(n+1)}_{\sigma\to\sigma^*} = \sum_{\sigma'\not=\sigma}
w_{\sigma\to\sigma'} f^{(n)}_{\sigma'\to\sigma^*} + w_{\sigma
\to \sigma^*} \qquad \sigma\not=\sigma^* \\
f^{(n+1)}_{\sigma^*\to\sigma^*} = \sum_{\sigma'\not=\sigma^*}
w_{\sigma^*\to\sigma'} f^{(n)}_{\sigma'\to\sigma^*} \end{cases}$
with initial condition $f^{(0)}=0$ 
and taking the limit $n\to\infty$. Therefore it is monotonous
in $w$, i.e. $\Big(w'_{\sigma\to\sigma'}\ge w_{\sigma\to\sigma'},
\ \sigma\not=\sigma'\in V\Big)\Rightarrow \Big(f(w')_{\sigma\to\sigma^*}
\ge f(w)_{\sigma\to\sigma^*}, \ \sigma\in V\Big)$.

We let 
\BEQ m:=\min_{\sigma}| \tilde{A}_{\sigma,\sigma}|, \qquad 
M:=\max_{\sigma}| \tilde{A}_{\sigma,\sigma}| \EEQ
and $c := \alpha/m$, $C:=\alpha/M$ $(c>C>0)$. 
Then
\BEQ \frac{1}{1+c} \tilde{w}(0)_{\sigma\to\sigma'} \le  \tilde{w}(\alpha)_{\sigma\to\sigma'}= \frac{\tilde{A}_{\sigma',\sigma}}{|\tilde{A}_{\sigma,\sigma}|+\alpha} 
\le  \frac{1}{1+C} \tilde{w}(0)_{\sigma\to\sigma'}
\EEQ
hence (by monotonicity in $w$)
\BEQ  \frac{1}{1+c}
\tilde{f}(0)_{\sigma\to\sigma^*}\le \tilde{f}(\alpha)_{\sigma\to\sigma^*} \le  \frac{1}{1+C}
\tilde{f}(0)_{\sigma\to\sigma^*}, \qquad \sigma\in \Sigma
\label{eq:tildefalphatildef0}
\EEQ
Also, monotonicity in $w$ and (\ref{eq:w>=wtilde}) imply 
\BEQ f(\alpha)_{\sigma\to\sigma^*}\ge \tilde{f}(\alpha)_{\sigma\to
\sigma^*}, \qquad \sigma\in \Sigma
\label{eq:f>=ftilde}
\EEQ
Conversely, 
\BEQ w(\alpha)_{\sigma\to\sigma'} = \frac{A_{\sigma',\sigma}}{
|A_{\sigma,\sigma}|+\alpha} = \frac{|\tilde{A}_{\sigma,\sigma}|+
\alpha}{|A_{\sigma,\sigma}|+\alpha}\  \tilde{w}(\alpha)_{\sigma\to\sigma'} \EEQ
hence (see (\ref{eq:D}))
\BEQ \frac{w(\alpha)_{\sigma\to\sigma'}}{\tilde{w}(\alpha)_{\sigma\to\sigma'}} \le \max_{\sigma} 
\frac{|\tilde{A}_{\sigma,\sigma}|+
\alpha}{|A_{\sigma,\sigma}|+\alpha} = 1+ \max_{\sigma}
\frac{\kappa_{\sigma}}{|A_{\sigma,\sigma}|+\alpha} \le 
1+D
\EEQ
from which we get by monotonicity: 
\BEQ f(\alpha)_{\sigma\to\sigma^*} \le (1+D)\tilde{f}(\alpha)_{\sigma\to\sigma^*} \label{eq:f<=Dftilde}
\EEQ

\Bigskip
Now,
\BEQ f(\alpha)_{\sigma^*\to\sigma^*} = \sum_{\sigma'\not=\sigma^*} w(\alpha)_{\sigma^*\to \sigma'}
f(\alpha)_{\sigma'\to\sigma^*}  
\EEQ
Using
\BEQ w(\alpha)_{\sigma^*\to \sigma'}
= \tilde{w}(0)_{\sigma^*\to \sigma'} \ \times\ 
\frac{|\tilde{A}_{\sigma^*,\sigma^*}|}{|\tilde{A}_{\sigma^*,\sigma^*}| - \kappa_{\sigma^*}+\alpha}. \label{eq:wwAkappa}
\EEQ
and (\ref{eq:tildefalphatildef0}), (\ref{eq:f>=ftilde}),
(\ref{eq:proba-cons}), we get
\BEA
f(\alpha)_{\sigma^*\to\sigma^*} &=& \frac{|\tilde{A}_{\sigma^*,\sigma^*}|}{|\tilde{A}_{\sigma^*,\sigma^*}| - \kappa_{\sigma^*}+\alpha}  \ \times\ \sum_{\sigma'\not=\sigma} \tilde{w}(0)_{\sigma^*\to\sigma'} f(\alpha)_{\sigma'\to\sigma^*} \\
&\ge & \frac{|\tilde{A}_{\sigma^*,\sigma^*}|}{|\tilde{A}_{\sigma^*,\sigma^*}| - \kappa_{\sigma^*}+\alpha}  \ \times\ \sum_{\sigma'\not=\sigma} \tilde{w}(0)_{\sigma^*\to\sigma'} \tilde{f}(\alpha)_{\sigma'\to\sigma^*} \nonumber\\
&\ge & h(m,\alpha|\sigma^*)  \sum_{\sigma'\not=\sigma} \tilde{w}(0)_{\sigma^*\to\sigma'} \tilde{f}(0)_{\sigma'\to\sigma^*} = h(m,\alpha|\sigma^*)
\EEA
where 
\BEQ h(k,\alpha|\sigma^*):= \frac{k}{k+\alpha} \frac{|\tilde{A}_{\sigma^*,\sigma^*}|}{|\tilde{A}_{\sigma^*,\sigma^*}| - \kappa_{\sigma^*}+\alpha}.
\EEQ 
Conversely, combining (\ref{eq:f<=Dftilde}),  (\ref{eq:tildefalphatildef0}) and  (\ref{eq:proba-cons}), we get
\BEA f(\alpha)_{\sigma^*\to\sigma^*} &=& \frac{|\tilde{A}_{\sigma^*,\sigma^*}|}{|\tilde{A}_{\sigma^*,\sigma^*}| - \kappa_{\sigma^*}+\alpha}  \ \times\ \sum_{\sigma'\not=\sigma} \tilde{w}(0)_{\sigma^*\to\sigma'} f(\alpha)_{\sigma'\to\sigma^*} \\
&\le & \frac{|\tilde{A}_{\sigma^*,\sigma^*}|}{|\tilde{A}_{\sigma^*,\sigma^*}| - \kappa_{\sigma^*}+\alpha}  \ \times\
\frac{1+D}{1+C}
\sum_{\sigma'\not=\sigma} \tilde{w}(0)_{\sigma^*\to\sigma'} \tilde{f}(0)_{\sigma'\to\sigma^*} \nonumber\\
& = & (1+D)h(M,\alpha|\sigma^*)
\EEA

\Medskip
If there exists any $\sigma^*$ for which $h(m,\alpha|\sigma^*)>1$,
then $\lambda^*>\alpha$. Conversely, if  there exists any $\sigma^*$ for which $(1+D)h(M,\alpha|\sigma^*)<1$,
then $\lambda^*<\alpha$. Therefore  

\BEQ  \sup\{\alpha\ge 0 \ |\ \max_{\sigma^*}h(m,\alpha|\sigma^*)>1\}\le \lambda^* \le \inf\{\alpha \ge 0 \ |\ \min_{\sigma^*} (1+D)h(M,\alpha|\sigma^*)<1\}.
\EEQ
The lower bound is inexistent if $\max_{\sigma^*}h(m,\alpha|\sigma^*)\le 1$ for all $\alpha\ge 0$,
in particular of course, if $\lambda^*<0$. 
 
\Medskip Introduce the quantities
\BEQ  y(m|\sigma^*)^2 = 4m\kappa_{\sigma^*}, \qquad  y(M|\sigma^*)^2 = 4M (D |\tilde{A}_{\sigma^*,\sigma^*}| + \kappa_{\sigma^*})
\EEQ
Note that $y(m|\sigma^*)^2<0$ if $\kappa_{\sigma^*}<0$; but then, 
$x(m|\sigma^*)^2 + y(m|\sigma^*)^2 = (|\tilde{A}_{\sigma^*,\sigma^*}| + |\kappa_{\sigma^*}| +m  )^2 
-4 m |\kappa_{\sigma^*}| \ge (|\kappa_{\sigma^*}| +m )^2 
-4 m |\kappa_{\sigma^*}| \ge 0$; this positivity property also holds with $M$ instead of $m$.  
Now, $h(m,\alpha|\sigma^*)>1$ if and only if 
$\alpha^2 + (|\tilde{A}_{\sigma^*,\sigma^*}| - \kappa_{\sigma^*}
+m)\alpha - m\kappa_{\sigma^*}<0$, equivalently,  if
\BEQ \alpha<\alpha_{thr}(m|\sigma^*):= 
\frac{-x(m|\sigma^*)+\sqrt{x(m|\sigma^*)^2+y(m|\sigma^*)^2}}{2},
\label{eq:alphathrmsqrt}
\EEQ
where
\BEQ x(m|\sigma^*)= |\tilde{A}_{\sigma^*,\sigma^*}| - \kappa_{\sigma^*}
+m = |A_{\sigma^*,\sigma^*}|+m, \EEQ

Conversely, $(1+D)h(M,\alpha|\sigma^*)<1$ if and only if
$\alpha^2 + (|\tilde{A}_{\sigma^*,\sigma^*}| - \kappa_{\sigma^*}
+M)\alpha - M(D|\tilde{A}_{\sigma^*,\sigma^*}|+ \kappa_{\sigma^*})<0$, equivalently,  if
\BEQ \alpha>\alpha_{thr}(M|\sigma^*):= 
\frac{-x(M|\sigma^*)+\sqrt{x(M|\sigma^*)^2+y(M|\sigma^*)^2}}{2},
\label{eq:alphathrMsqrt}
\EEQ 
Thus we get our key formula,
\BEQ \max_{\sigma^*}\alpha_{thr}(m|\sigma^*)\le \lambda^* \le  \min_{\sigma^*}\alpha_{thr}(M|\sigma^*).
\label{eq:double-inequality-lambda*}
\EEQ


\Bigskip  {\bf B.}  {\em In the following, we prove the estimates  for $\lambda^*$ given in the three scaling regimes (weakly autocatalytic/strongly autocatalytic/intermediate) of 
Theorem \ref{thm:homogeneous}
under  
homogeneity conditions for the underlying 1-1 reaction network.}  Since there is no degradation by assumption, $\kappa_{\sigma}\ge 0$ for all $\sigma\in V$ now. 
  We assume that all non 1-1 reactions $\rho$ are of order $p(\rho)=2$, so that 
\BEQ |A_{\sigma,\sigma}| = |A^{(0)}_{\sigma,\sigma}|+\kappa_{\sigma}, \qquad   |\tilde{A}_{\sigma,\sigma}| = |A^{(0)}_{\sigma,\sigma}|+2\kappa_{\sigma}
\EEQ
where $|A^{(0)}|_{\sigma,\sigma}$ is the sum of rates of 
1-1 reactions leaving $\sigma$. This hypothesis is not
needed (it would be enough to assume that 
$\max_{\rho\in R} p(\rho)<\infty$ to get the same
final estimates), however it is very natural when one thinks
of biochemical reactions. Note that $D\le 1$ under this
condition.
 The important (homogeneity) hypothesis is the
following: 

\begin{center} {\bf  (H) \qquad All $|A^{(0)}_{\sigma,\sigma}|\approx k$ are  of the same order}
\end{center}

In particular, we then have: 
\BEQ m\approx k + \min_{\sigma} \kappa_{\sigma}, \qquad 
 M\approx k + \max_{\sigma} \kappa_{\sigma}
\label{eq:m}
\EEQ
\BEQ D\approx \max_{\sigma} \Big(\frac{\kappa_{\sigma}}{k+\kappa_{\sigma}} \Big) \approx \min\Big(1, \frac{\max_{\sigma} \kappa_{\sigma}}{k}\Big)
\EEQ

\Bigskip {\em We start with the lower bound.} Let
\BEQ  g(\sigma^*):=(\frac{x(m|\sigma^*)^2}{y(m|\sigma^*)^2} = 
\frac{(|A_{\sigma^*,\sigma^*}| + 
m)^2}{4m\kappa_{\sigma^*}}   \label{eq:g}
\EEQ
then under our assumptions,   $g(\sigma^*)\approx \frac{(2m+\kappa_{\sigma^*})^2}{4m\kappa_{\sigma^*}} = \frac{m}{\kappa_{\sigma}^*} + 1 +
\frac{\kappa_{\sigma^*}}{4m}$. This expression, as a function 
of $\kappa_{\sigma^*}$,  has a minimum
at $\kappa_{\sigma^*}=2m$, namely, $2$. If $x\gg y>0$, then
(by Taylor expansion)
$\frac{-x+\sqrt{x^2+y^2}}{2}\sim \frac{y^2}{4x}$, but more
generally, $\frac{-x+\sqrt{x^2+y^2}}{2}\approx \frac{y^2}{x}$
in the whole domain $x>y>0$.  The r.-h.s. of
(\ref{eq:alphathrmsqrt}) depends critically on this ratio.
We get 4 cases.

\Medskip 
 (1) If $\kappa_{\sigma^*}\ll k$ {\em (weak autocatalytic rate)},
which implies in particular $m\approx k$, 
then $g(\sigma^*)\gg 1$, $x(m|\sigma^*)\gg y(m|\sigma^*)$ and
$\alpha_{thr}(\sigma^*)\sim 
\frac{y(m|\sigma^*)^2}{4x(m|\sigma^*)} = \frac{m\kappa_{\sigma^*}}{|A_{\sigma^*,\sigma^*}|+m} \sim \frac{\kappa_{\sigma^*}}{2}
$ increases linearly in $\kappa_{\sigma^*}$ (perturbative
regime).

\Medskip If, on the other hand, $\kappa_{\sigma^*}\gtrsim k$  {\em (strong autocatalytic rate)},
then $\kappa_{\sigma^*}\gtrsim m$, but we do not necessarily have
$\kappa_{\sigma^*}\gg m$, and the Taylor expansion does not
 give a correct lower bound for $\alpha_{thr}(\sigma^*)$ if the
 latter condition is not satisfied. We distinguish three cases:
 
 \Medskip(2)\  if  $\kappa_{\sigma^*}\gg m$ (equivalently,
 $\kappa_{\sigma^*}\gg k$ and $\kappa_{\sigma^*}\gg \min_{\sigma}
 \kappa_{\sigma}$, so that autocatalytic rates differ widely
 in order of magnitude), then the Taylor expansion is correct,
 and  one gets
$\alpha_{thr}(\sigma^*)\approx \frac{y(m|\sigma^*)^2}{4x(m|\sigma^*)} = 
\frac{m\kappa_{\sigma^*}}{|A_{\sigma^*,\sigma^*}|+m} \approx
\frac{m\kappa_{\sigma^*}}{k+\kappa_{\sigma^*}+m}\approx m$
saturates as $\max_{\sigma^*}\to\infty$ while $\min_{\sigma}
\kappa_{\sigma}$ is fixed. 

\Medskip(3)\ if $\kappa_{\sigma^*}\approx k$, then $m\approx k$,
$|A_{\sigma^*,\sigma^*}|\approx k+\kappa_{\sigma^*}\approx k$,
and $x(m|\sigma^*),y(m|\sigma^*)\approx k$ too. Then 
$\alpha_{thr}(\sigma^*)\approx k\approx \kappa_{\sigma^*}$, which
is the relevant scale of the problem. 

\Medskip(4)\ finally, if $\kappa_{\sigma^*}\gg k$ but
$\kappa_{\sigma^*}\approx \min_{\sigma}\kappa_{\sigma}$ (namely,
if all autocatalytic rates dominate w.r. to Markov rates 
which are of order $k$), 
then $m,y(m|\sigma^*),x(m|\sigma^*)\approx \kappa_{\sigma^*}$, 
hence $\alpha_{thr}(\sigma^*)\approx \kappa_{\sigma^*}$. 

\Medskip Summarizing, we have three scales in this problem:
$k$; $\min_{\sigma}\kappa_{\sigma}$ and $\kappa_{\sigma^*}$. 
Apart from the perturbative regime (1) where $\alpha_{thr}(\sigma^*)\approx \kappa_{\sigma^*}$ is dictated by the
kinetically limiting step, in all other cases, 
 $\alpha_{thr}(\sigma^*)$ is subdominant, i.e. it is  of the same order as the
intermediate scale: 
(2) if  $\kappa_{\sigma^*}\gg  k,\min_{\sigma}\kappa_{\sigma}$, 
then $\alpha_{thr}(\sigma^*) \approx m\approx \max(k,\min_{\sigma}\kappa_{\sigma})$; (3) if $k\approx \kappa_{\sigma^*} \gtrsim 
\min_{\sigma}\kappa_{\sigma}$, then $\alpha_{thr}(\sigma^*)
\approx k$; (4) if $\kappa_{\sigma^*},\min_{\sigma}\kappa_{\sigma}
\gg k$, then $\alpha_{thr}(\sigma^*)\approx \kappa_{\sigma^*}$. 

\Medskip Maximizing in $\sigma^*$, we get the estimates 
$\max_{\sigma^*}\alpha_{thr}(m|\sigma^*) \approx \begin{cases}
\max_{\sigma}\kappa_{\sigma} \\ \min_{\sigma} \kappa_{\sigma} \\
k\end{cases}$ depending on the scaling regime (weakly autocatalytic/strongly autocatalytic/intermediate), see statement of Theorem \ref{thm:homogeneous}.


\Bigskip  {\em Let us now consider the upper bound.} We first remark that
$M\approx \max(k, \max_{\sigma}\kappa_{\sigma})$, $x(M|\sigma^*)
\approx M$, $D\approx \frac{\max_{\sigma}\kappa_{\sigma}}{
k+ \max_{\sigma}\kappa_{\sigma}} \approx \frac{\max_{\sigma}
\kappa_{\sigma}}{M}$. Then  $D|\tilde{A}_{\sigma^*,\sigma^*}| \ge 
\frac{\kappa_{\sigma^*}}{|A_{\sigma^*,\sigma^*}|} \ \times\ 
|A_{\sigma^*,\sigma^*}| = \kappa_{\sigma^*}$, whence $y(M|\sigma^*)\approx 
\sqrt{MD|\tilde{A}_{\sigma^*,\sigma^*}|}$. Finally, redefining
$g(\sigma^*)$ as $(\frac{x(M|\sigma^*)}{y(M|\sigma^*)})^2$, 
we have
\BEQ  g(\sigma^*) \approx \frac{M}{D|\tilde{A}_{\sigma^*,\sigma^*}|} \gtrsim \frac{M}{D\max_{\sigma}|\tilde{A}_{\sigma^*,\sigma^*}|} \approx 1 \EEQ
Thus the estimates on $\alpha_{thr}(M|\sigma^*)$ can be
obtained in the same way as those for $\alpha_{thr}(m|\sigma)$,   see paragraph below (\ref{eq:g}).

\Medskip The two scaling parameters here are $k$ and 
$\max_{\sigma} \kappa_{\sigma}$. We find 
3 possibilities:

\Medskip (1) {\em (weakly autocatalytic regime)} if  $\max_{\sigma}(\kappa_{\sigma})\ll k$, then $M\approx k$, 
$D|\tilde{A}_{\sigma^*,\sigma^*}|\approx \frac{\max_{\sigma}(\kappa_{\sigma})}{M} \, \times\, k \approx \max(\kappa_{\sigma})$,
$g(\sigma^*)\gg 1$ and $\alpha_{thr}(M|\sigma^*)\sim 
\frac{y(M|\sigma^*)^2}{4x(M|\sigma^*)} \sim \frac{\kappa_{\sigma}^*}{2}$, exactly as for the lower bound obtained in  case (1)
above.

\Medskip (2) If $\kappa_{\sigma^*}\gtrsim k$  {\em (strong autocatalytic rate)}, then $\max_{\sigma}\kappa_{\sigma}\gtrsim k$,  
whence $M\sim \max_{\sigma}(\kappa_{\sigma})$, $D\approx 1$, 
$D|\tilde{A}_{\sigma^*,\sigma^*}| \approx \kappa_{\sigma^*}$,
and $\alpha_{thr}(M|\sigma^*)\approx \frac{y(M|\sigma^*)^2}{x(M|\sigma^*)} \approx \frac{MD|\tilde{A}_{\sigma^*,\sigma^*}|}{\max_{\sigma} \kappa_{\sigma}} \approx \kappa_{\sigma^*}$. 

\Medskip (3) If $\max_{\sigma}\kappa_{\sigma} \gtrsim k \gtrsim 
\kappa_{\sigma^*}$, then $M\sim \max_{\sigma}(\kappa_{\sigma})$, $D\approx 1$ as before, but now $D|\tilde{A}_{\sigma^*,\sigma^*}| \approx k$, and $\alpha_{thr}(M|\sigma^*)\approx \frac{y(M|\sigma^*)^2}{x(M|\sigma^*)} \approx \frac{MD|\tilde{A}_{\sigma^*,\sigma^*}|}{\max_{\sigma} \kappa_{\sigma}} \approx k$.

\Medskip Minimizing in $\sigma^*$, we get the estimates 
$\min_{\sigma^*}\alpha_{thr}(M|\sigma^*) \approx \begin{cases}
\max_{\sigma}\kappa_{\sigma} \\ \min_{\sigma} \kappa_{\sigma} \\
k\end{cases}$ depending on  the scaling regime (weakly autocatalytic/strongly autocatalytic/intermediate) 
as above.

\vskip 1 cm
\noindent Combining the two bounds, we have finally proved Theorem \ref{thm:homogeneous}. \hfill\eop


\subsection{Resonances}


\noindent {\em Resonance regime of Example 2} (see \S \ref{section:example2}).  The adjoint generator is
$A = \left[\begin{array}{ccc} -k_{on}-k'_4 & k_{off} + 2\nu_+ & 0 \\ k_{on} & -k_{off}-\nu_+ & k_4 \\
k'_4 & 0 & -k_4 \end{array}\right]$. We find

\BEQ -\det(A-\lambda) = (\lambda+k_{on}+k'_4)(\lambda+k_{off}+\nu_+)(\lambda+k_4) - k_4(k_{off}+2\nu_+)(k_{on}+k'_4) - \lambda (k_{off}+2\nu_+)k_{on}.
\EEQ
The determinant vanishes for $\lambda=\lambda^*(k'_4)$. 
The transition regime is defined by $k'_4\sim \frac{k_{on}\nu_+}{k_{off}}$; in this regime, 
$\lambda^*\preceq \frac{k_{on}}{k_{off}}\nu_+$. Expanding the above determinant yields to leading orders
\BEQ \lambda^2 k_{off} + \lambda\Big\{  \Big[(k_{off}+\nu_+)k'_4-\nu_+ k_{on} \Big]
 + k_{off} k_4 \Big\} \sim k_{on} \nu_+ k_4 
\EEQ
Note that  $k_{off}k_4 \ll (k_{off}+\nu_+)k'_4 \sim  \nu_+ k_{on}$ in the transition regime. The term between curly brackets vanishes for some value $(k'_4)^0$ of $k'_4$ located in the transition regime. Let $\del k'_4 := k'_4 - (k'_4)^0$, and $\eta := \frac{\del k'_4}{k'_4}$ the logarithmic 
discrete derivative; the determinant vanishes when $\lambda=\lambda^*$ solves the second degree equation
$ \lambda^2 k_{off} + \lambda\eta\nu_+ k_{on} \sim \nu_+ k_{on} k_4.  $
 When $\eta=0$, we get 
 \BEQ \lambda^* \equiv \lambda^0 \sim \sqrt{\nu_+ \frac{k_{on}}{k_{off}} k_4}. \EEQ
 Letting $\eta$ vary, the discriminant $\Del$ is such that $\sqrt{\Del} \sim 2 \sqrt{\nu_+ k_{off}k_{on} k_4}
+ \frac{(\nu_+ k_{on}\eta)^2}{4\sqrt{\nu_+ k_{off}k_{on} k_4}}$ whenever $ |\eta|\preceq \sqrt{k_4/k'_4}$, yielding under this condition (with $\del \lambda:=\lambda^*-
\lambda^0$)
\BEQ \frac{\del \lambda}{\lambda}\sim -\frac{\nu_+ k_{on}}{2k_{off}\lambda^0} \frac{\del k'_4}{k'_4} 
\sim - \sqrt{\frac{k'_4}{k_4}} \, \frac{\del k'_4}{k'_4} \EEQ
When $ |\eta|\sim\sqrt{k_4/k'_4}$, we get $|\frac{\del \lambda}{\lambda}|\sim 1$.


\subsection{Proof of Theorem \ref{th:sigma*-dominant}}


\noindent We assume here that the $\sigma^*$-dominant Hypothesis of Section \ref{section:sigma*} holds.
 We split the proof of Theorem \ref{th:sigma*-dominant} into several points.
 
\begin{enumerate}
\item {\em (dominant simple excursion from $\sigma^*$ to $\sigma^*$)} Start from the dominant cycle
$\gamma^*:\sigma^*=\sigma_1\to\cdots\to \sigma_{\ell^*}=\sigma^*$.  Assume it contains an inner loop, i.e. that
there exist $1\le i<j\le \ell^*$ (excluding the trivial case $(i,j)=(1,\ell^*)$) such that $\sigma_i=\sigma_j$; then the contracted
cycle $\sigma_1\to\cdots\to \sigma_i\to\sigma_{j+1}\to\cdots\to\sigma_{\ell^*}$ is still a dominant cycle. 
Iterating, we get in the end a dominant cycle from $\sigma^*$ to $\sigma^*$ without inner loop, i.e. a dominant simple excursion from $\sigma^*$ to $\sigma^*$.
 Since there are only finitely many simple excursions
from $\sigma^*$ to $\sigma^*$ (because $V$ is finite), there exists a dominant
   simple excursion  (denoted $\gamma^*_{\sigma^*}$)
 with  maximum length  $\ell_{\sigma^*}-1\le |V|$.  This implies (i). 
 
\item {\em (leading simple excursion from $\sigma^*$ to $\sigma$)} Let $\sigma\not=\sigma^*$. 
Start from any excursion $\gamma:\sigma^*=\sigma_1\to\cdots\sigma_{\ell}=\sigma$ from $\sigma^*$ to $\sigma$ 
minimizing the sum in (\ref{eq:nsigma*totosigma}). Assume as in 1; that it contains an inner loop, then
contracting it yields a shorter excursion minimizing (\ref{eq:nsigma*totosigma}), since all terms in
the sum are $\ge 0$. Iterating, we get in the end a leading simple excursion $\gamma^*_{\sigma}$. This proves the first
part of (ii). 

\item {\em (exponential decrease for long excursions)} {\em  The total weight of 
excursions $\gamma:\sigma^*\to\sigma^*$ of length $\ell>|V| $ decreases exponentially with $\ell$, 
\BEQ f^*_{\ell}:= \sum_{\gamma:\sigma^*\to \sigma^*\ {\mathrm{excursion}}, \ \ell(\gamma)=\ell} w(\lambda^*)_{\gamma}  \le c^{\ell/|V|} 
\EEQ
with $c \prec 1$. }

Namely, let $\ell>|V|$; then any excursion $\gamma:\sigma^*\to\sigma^*$ of length $\ell$ must contain
an inner loop $\gamma':\sigma_i\to\cdots\to\sigma_j$ with   $\sigma_i=\sigma_j\equiv\sigma\not=\sigma^*$ and $1<i<j<\ell$. 
Choosing $i$ minimal, and then $j$ minimal (following the loop-erasing construction \cite{LERW}), 
we actually have $j-1\le |V|$, and $\gamma'$ is a simple cycle from $\sigma$ to $\sigma$. By the 
$\sigma^*$-dominant Hypothesis,  $\gamma'$ is non-dominant, whence $w(\lambda^*)_{\gamma'}\prec 1$. Summing
over all possibilities for $\sigma_2,\ldots,\sigma_j$, we get $f^*_{\ell} \le Cb \max_{1\le i\le |V|}
f^*_{\ell-i}$, where $b$ is the base constant, and $C$ is some finite, $b$-independent constant. Thus (by an 
easy induction)
$f^*_{\ell} \le (Cb)^{\ell/|V|-1} f^*_{|V|}$, implying exponential decrease.  

\item 
Similarly, {\em let $\omega_{\ge m}^*$, $m\ge |V|$,  be the total weight
of paths $\gamma:\sigma^*=\sigma_1\to\cdots\to\sigma_n$ $(n\ge m)$ such that 
$\sigma_i\not=\sigma^*$ for all $i>1$ ("incomplete" excursions from
$\sigma^*$ to $\sigma^*$) ; then 
$f_{\ge m}^* \le c^{m/|V|}$ with $c\prec 1$.}  The proof (adapted from 2.) is left to the reader. 

\item Similarly again, let $\omega^{*,\sigma}_{\ell}, \ell> |V|$, be the total weight of paths 
$\gamma:\sigma=\sigma_1\to\cdots\to \sigma_{\ell} = \sigma$ {\em not passing through} $\sigma^*$. Then
$\omega^{*,\sigma}_{\ell}\le (Cb)^{\ell/|V|-1} \omega^{*,\sigma}_{|V|}$. Then, the $\sigma^*$-dominant
hypothesis implies $\omega^{*,\sigma}_{|V|}\prec 1$. Summing over $\ell>|V|$, we get (iii).

\item Let $\omega^*_n = \sum_{\gamma:\sigma^*\to\sigma^* \, {\mathrm{cycle}},\ \ell(\gamma)=n+1}
w(\lambda^*)_{\gamma}$ be the total weight of paths of length $n$ from 
$\sigma^*$ to $\sigma^*$, and $\Omega^*_n := \sum_{i=\max(n,\ell^*)}^{n+\ell^*-1} \omega^*_i$. We already know from 1. that $\Omega^*_{\ell^*} \ge \cdots\ge \Omega^*_1\ge 
w(\lambda^*)_{\gamma^*}\succeq 1$. 
For $n\ge \ell^*$, we have the recurrence relation
\BEA  \Omega^*_n & =& \omega^*_n + \ldots+ \omega^*_{n+\ell^*-1} \nonumber\\
  &\ge& \sum_{\gamma:\sigma^*\to\sigma^* \ {\mathrm{excursion}},\  2\le \ell(\gamma)<n} w(\lambda^*)_{\gamma} 
\Omega^*_{n-\ell(\gamma)} - O(f_{\ge n}^*) \nonumber\\
&\ge & (1-c^n)  \min_{1\le n'<n} \Omega^*_n - O(c^n)
\EEA
for some $c\prec 1$ by 2. and 3. We have used the resolvent identity: $\sum_{\gamma:\sigma^*\to\sigma^*\  {\mathrm{excursion}}}
w(\lambda^*)_{\gamma} \equiv 1$.  From this we deduce by an easy induction that 
$\min_{n\ge 1}\succeq 1$. This implies (iv).

\item {\em (leading simple excursion from $\sigma^*$ to $\sigma$, cont'd)} Reasoning as in 5. 
but with the total weight $f^*_{\ell}(\sigma^*\to\sigma)$ of excursions $\gamma:\sigma^*\to\sigma$ of
length $\ell$, we finally prove: $ f^*_{\sigma^*\to\sigma} \sim w(\lambda^*)_{\gamma^*_{\sigma}}
$, which concludes (ii).

\end{enumerate}


\section*{Index of definitions and notations} \label{section:index}


\begin{tabular}{cc}
dominant edge, \S \ref{section:multi-scale-splitting} & dominant SCC, \S \ref{subsection:homogeneous} \\
  $\sigma^*$-dominant graph, 
\S \ref{section:sigma*} & 
leading path, \S \ref{section:multi-scale-splitting} \\
merging, \S \ref{subsection:homogeneous} & rewiring, \S \ref{subsection:homogeneous} \\
\\  $a_{\gamma}$, \S \ref{subsection:threshold} & 
$a_{\sigma} = |A_{\sigma,\sigma}|$, \S \ref{section:multi-scale-splitting} \\ 
$A$, adjoint generator \S \ref{subsection:Markov} & $\tilde{A}$, proba. preserving generator
\S \ref{subsection:Markov} \\
$c_k(\alpha)$,  \S \ref{subsection:threshold} & $c'_{high}$, \S \ref{subsection:threshold} \\
$d_{\sigma}$, deficiency scale \S \ref{subsection:threshold-integer}  & $D_{\gamma}$, excursion
depth \S \ref{subsection:threshold-integer} \\
$E$, edge set of Markov graph, \S \ref{subsection:Markov} &
$f(\alpha)_{\sigma\to\sigma'}$, excursion weight \S \ref{subsection:resolvent-formula} \\
$G$, Markov graph\S \ref{subsection:Markov}  & 
$J^{\rho},J_{lin}^{\rho}$, reaction currents   \S \ref{subsection:maths} \\
$k_{\sigma}$, total outgoing rate \S \ref{section:multi-scale-splitting} & $k_{\sigma\to\sigma'}$, 
transition rate \S \ref{subsection:Markov} \\
$n_i$, scales \S \ref{subsection:notations} &
$n(k)$, reaction scale \S \ref{subsection:aim}, \ref{subsection:notations} \\ 
$n_{\gamma}$, \S \ref{subsection:threshold} &
$n_{\sigma}$, vertex scale \S \ref{section:multi-scale-splitting} \\
  $n_{\sigma\to\sigma'}$,
edge scale  \S \ref{section:multi-scale-splitting} &  $n_{\sigma^*\twoheadrightarrow \sigma}$, 
\S \ref{section:sigma*} \\
$n'_{high}$, \S \ref{subsection:threshold-integer} \\
$N_{\gamma}$, \ref{subsection:Markov} & 
$p(\rho)$, reaction order \S \ref{subsection:Markov} \\
$q_{1,2}$, \S \ref{subsection:threshold-integer} & $R_p$, order $p$ reaction set
\S \ref{subsection:Markov} \\
$R(\alpha)$, resolvent \S \ref{subsection:resolvent-formula} &  ${\mathbb S}$, stoechiometric
matrix \S \ref{subsection:maths} \\
  $V$, species set \S \ref{subsection:Markov} & 
 $v^*$, Lyapunov eigenvector,  \S \ref{subsection:Markov}  \\
$w_{\sigma\to\sigma'}, \tilde{w}_{\sigma\to\sigma'}$, transition weights \S \ref{subsection:resolvent-formula} & $w(\alpha)_{\sigma\to\sigma'}, \tilde{w}(\alpha)_{\sigma\to\sigma'}$, transition weights \S \ref{subsection:resolvent-formula} \\
\\
$\alpha$, resolvent parameter \S \ref{subsection:resolvent-formula} &  $\alpha_{thr}$, threshold
parameter \S \ref{subsection:resolvent-formula}, (\ref{eq:th:thr})   \\
$\beta_{\sigma}$, degradation rate  \S \ref{subsection:maths} & 
$\gamma$, path \S \ref{section:resolvent}  \\ 
$\gamma^{high}, \gamma^{low}$, \S \ref{subsection:threshold}  & 
$\del c_{high}$, $\del c_{high,i}$,  \S \ref{subsection:threshold} \\
$\eps_{\sigma}$, deficiency weight \S \ref{section:multi-scale-splitting} & 
$\kappa_{\sigma}$, deficiency rate \S  \ref{subsection:Markov}, \ref{section:multi-scale-splitting}  \\ 
$\lambda^*$, Lyapunov exponent \S \ref{subsection:Markov} & $\tilde{\mu}$, 
stationary measure \S \ref{subsection:Markov}   \\
$\tilde{\pi}$, stationary (discrete-time) measure \S \ref{subsection:resolvent-formula} & 
$\Pi_{\gamma}$,  \S \ref{subsection:threshold} \\
$\rho$, reaction index & $\sigma$, species index  \\
$\sim,\simeq$, \S \ref{subsection:notations} & $\prec,\succ$, \S \ref{subsection:notations}  \\
\end{tabular}

\newpage

\Medskip {\em Data availability statement.} There are no data supporting the results.
   
\Medskip {\em Conflict of interest statement.} There is no conflict of interest.


\end{document}